%
%
%
\newif\ifsect\newif\iffinal
\secttrue\finaltrue
\def\smallsect #1. #2\par{\bigbreak\noindent{\bf #1.}\enspace{\bf #2}\par
	\global\parano=#1\global\eqnumbo=1\global\thmno=1
	\nobreak\smallskip\nobreak\noindent\message{#2}}
\def\thm #1: #2{\medbreak\noindent{\bf #1:}\if(#2\thmp\else\thmn#2\fi}
\def\thmp #1) { (#1)\thmn{}}
\def\thmn#1#2\par{\enspace{\sl #1#2}\par
        \ifdim\lastskip<\medskipamount \removelastskip\penalty 55\medskip\fi}
\def\square{{\msam\char"03}}
\def\qedn{\thinspace\null\nobreak\hfill\square\par\medbreak}
\def\pf{\ifdim\lastskip<\smallskipamount \removelastskip\smallskip\fi
        \noindent{\sl Proof\/}:\enspace}

\newcount\parano
\newcount\eqnumbo
\newcount\thmno
\newcount\versiono
\def\neweqt#1$${\xdef #1{(\number\parano.\number\eqnumbo)}
	\eqno #1$$
	\iffinal\else\rsimb#1\fi
	\global \advance \eqnumbo by 1}
\def\newthmt#1 #2: #3{\xdef #2{\number\parano.\number\thmno}
	\global \advance \thmno by 1
	\medbreak\noindent
	\iffinal\else\lsimb#2\fi
	{\bf #1 #2:}\if(#3\thmp\else\thmn#3\fi}
\def\neweqf#1$${\xdef #1{(\number\eqnumbo)}
	\eqno #1$$
	\iffinal\else\rlap{$\smash{\hbox{\hfilneg\string#1\hfilneg}}$}\fi
	\global \advance \eqnumbo by 1}
\def\newthmf#1 #2: #3{\xdef #2{\number\thmno}
	\global \advance \thmno by 1
	\medbreak\noindent
	\iffinal\else\llap{$\smash{\hbox{\hfilneg\string#1\hfilneg}}$}\fi
	{\bf #1 #2:}\if(#3\thmp\else\thmn#3\fi}
\def\inizia{\ifsect\let\neweq=\neweqt\else\let\neweq=\neweqf\fi
\ifsect\let\newthm=\newthmt\else\let\newthm=\newthmf\fi}
\def\bititolo{\empty}
\gdef\begin #1 #2\par{\xdef\titolo{#2}
\ifsect\let\neweq=\neweqt\else\let\neweq=\neweqf\fi
\ifsect\let\newthm=\newthmt\else\let\newthm=\newthmf\fi
\centerline{\titlefont\titolo}
\if\bititolo\empty\else\medskip\centerline{\titlefont\bititolo}
\xdef\titolo{\titolo\ \bititolo}\fi
\bigskip
\centerline{\bigfont \autore}
\if\istituto!\else\bigskip
\centerline{\istituto}
\centerline{\indirizzo}
\centerline{\email}\fi
\medskip
\centerline{#1~\anno}
\bigskip\bigskip
\ifsect\else\global\thmno=1\global\eqnumbo=1\fi}
\def\istituto{Dipartimento di Matematica, Universit\`a di Roma ``Tor Vergata''}
\def\indirizzo{Via della Ricerca Scientifica, 00133 Roma, Italy}
\def\email{E-mail: abate@mat.uniroma2.it}
\def\anno{\the\year}
\font\titlefont=cmssbx10 scaled \magstep1
\font\bigfont=cmr12 
\font\eightrm=cmr8
\font\sc=cmcsc10
\font\bbr=msbm10
\font\sbbr=msbm7 
\font\ssbbr=msbm5
\font\msam=msam10
\def\ca #1{{\cal #1}}
\nopagenumbers
\binoppenalty=10000
\relpenalty=10000
\newfam\amsfam
\textfont\amsfam=\bbr \scriptfont\amsfam=\sbbr \scriptscriptfont\amsfam=\ssbbr
\let\de=\partial
\def\eps{\varepsilon}
\def\phe{\varphi}

\def\Aut{\mathop{\rm Aut}\nolimits}

\def\End{\mathop{\hbox{\rm End}}}
\def\Ker{\mathop{\rm Ker}\nolimits}

\def\Re{\mathop{\rm Re}\nolimits}

\def\sp{\mathop{\rm sp}\nolimits}
\def\id{\mathop{\rm id}\nolimits}
\mathchardef\void="083F
\def\C{{\mathchoice{\hbox{\bbr C}}{\hbox{\bbr C}}{\hbox{\sbbr C}}
{\hbox{\sbbr C}}}}
\def\N{{\mathchoice{\hbox{\bbr N}}{\hbox{\bbr N}}{\hbox{\sbbr N}}
{\hbox{\sbbr N}}}}
\def\P{{\mathchoice{\hbox{\bbr P}}{\hbox{\bbr P}}{\hbox{\sbbr P}}
{\hbox{\sbbr P}}}}
\newcount\notitle
\notitle=1
\headline={\ifodd\pageno\rhead\else\lhead\fi}
\def\rhead{\ifnum\pageno=\notitle\iffinal\hfill\else\hfill\tt Version 
\the\versiono; \the\day/\the\month/\the\year\fi\else\hfill\eightrm\titolo\hfill
\folio\fi}
\def\lhead{\ifnum\pageno=\notitle\hfill\else\eightrm\folio\hfill\autore\hfill
\fi}
\newbox\bibliobox
\def\setref #1{\setbox\bibliobox=\hbox{[#1]\enspace}
	\parindent=\wd\bibliobox}
\def\biblap#1{\noindent\hang\rlap{[#1]\enspace}\indent\ignorespaces}
\def\art#1 #2: #3! #4! #5 #6 #7-#8 \par{\biblap{#1}#2: {\sl #3\/}.
	#4 {\bf #5} (#6)\if.#7\else, \hbox{#7--#8}\fi.\par\smallskip}
\def\book#1 #2: #3! #4 \par{\biblap{#1}#2: {\bf #3.} #4.\par\smallskip}
\def\coll#1 #2: #3! #4! #5 \par{\biblap{#1}#2: {\sl #3\/}. In {\bf #4,} 
#5.\par\smallskip}
\def\pre#1 #2: #3! #4! #5 \par{\biblap{#1}#2: {\sl #3\/}. #4, #5.\par\smallskip}
\versiono=4
\def\autore{Marco Abate\footnote{${}^1$}{\rm Partially supported by Progetto MURST di
Rilevante Interesse Nazionale {\it Propriet\`a geometriche delle variet\`a reali e
complesse.}}}
%
%
%
%
%
%
\begin {April} Diagonalization of non-diagonalizable discrete holomorphic dynamical systems

{\narrower {\sc Abstract.} We shall describe a canonical procedure to associate to any 
(germ of) holomorphic self-map~$F$ of~$\C^n$ fixing the origin so that $dF_O$ is
invertible and non-diagonalizable an~$n$-dimensional complex manifold~$M$, a holomorphic
map~$\pi\colon M\to\C^n$, a point ${\bf e}\in M$ and a (germ of) holomorphic self-map~$\tilde
F$ of~$M$ such that: $\pi$ restricted to $M\setminus\pi^{-1}(O)$ is a biholomorphism between
$M\setminus\pi^{-1}(O)$ and $\C^n\setminus\{O\}$; $\pi\circ\tilde F=F\circ\pi$; and
$\bf e$ is a fixed point of~$\tilde F$ such that $d\tilde F_{{\bf e}}$ is
diagonalizable. Furthermore, we shall use this construction to describe the local dynamics of
such an~$F$ nearby the origin when~$\sp(dF_O)=\{1\}$.\footnote{}{{\it 1991 Mathematical
Subjects Classification:} Primary 32H50, 32H02, 58F23}

}

\smallsect 0. Introduction

In passing from one to several variables, possibly the first new phenomenon one has to deal
with is the existence of non-diagonalizable linear maps. Roughly speaking, one can think
of them as some sort of singularity in the space of all linear maps; indeed, a generic
linear endomorphism is diagonalizable. It would be interesting to have a device to
``resolve'' the singularity, similarly to what happens in algebraic geometry for
singularities of complex spaces. 

In this paper we shall describe exactly such a device, in a more general holomorphic
setting. Let $F\in\End(\C^n,O)$ be a (germ of) holomorphic self-map of $\C^n$ keeping the
origin fixed and such that $dF_O$ is invertible and non-diagonalizable. We shall build in a
canonical way (depending only on the block structure of the Jordan form of $dF_O$) a new
holomorphic map
$\tilde F$ semi-conjugate to~$F$ (and actually conjugate to~$F$ outside the origin) with a
canonical fixed point~{\bf e} such that $d\tilde F_{{\bf e}}$ is diagonalizable; the price to
pay is that we have to change the base manifold. We shall in fact prove the
following result (see Theorem~2.4):\def\autore{Marco Abate}

\newthm Theorem \zuno: (Diagonalization Theorem)
Let $F\in\End(\C^n,O)$ be such that $dF_O$ is invertible and non-diagonalizable. Then there
exist a complex $n$-dimensional manifold~$M$, a holomorphic projection $\pi\colon
M\to\C^n$, a canonical point ${\bf e}\in M$ and a (germ at $\pi^{-1}(O)$ of)
holomorphic self-map $\tilde F\colon M\to M$ such that:
\smallskip
\item{\rm(i)}$\pi$ restricted to $M\setminus\pi^{-1}(O)$ is a biholomorphism between
$M\setminus\pi^{-1}(O)$ and $\C^n\setminus\{O\}$;
\item{\rm(ii)}$\pi\circ\tilde F=F\circ\pi$;
\item{\rm(iii)}$\bf e$ is a fixed point of~$\tilde F$, and $d\tilde F_{{\bf e}}$ is
diagonalizable.
\smallskip
\noindent More precisely, if the Jordan canonical form of~$dF_O$ contains $\rho\ge1$ blocks
of length $\mu_1\ge\mu_2\ge\cdots\ge\mu_\rho\ge1$ corresponding respectively to eigenvalues
$\lambda_1$,~$\lambda_2,\ldots,\lambda_\rho\in\C$, then $d\tilde F_{\bf e}$ has eigenvalues
$\tilde\lambda_1$,~$1$,~$\lambda_2/\lambda_1,\ldots,\lambda_\rho/\lambda_1$ of multiplicity
respectively~$1$,~$\mu_1-1$,~$\mu_2,\ldots,\mu_\rho$, where $\tilde\lambda_1=\lambda_1$ if
$\mu_1>\mu_2$, and $\tilde\lambda_1=\lambda_1^2/\lambda_2$ if~$\mu_1=\mu_2$.

One subtle point must be stressed here. If the only aim is to diagonalize the differential,
one can choose among several different constructions; but most of them are useless for the
dynamical applications we have in mind. For instance, the standard way to resolve
singularities in algebraic geometry is by blowing up points. One could do the same here:
$M$ could be obtained by $\C^n$ blowing up a suitable sequence of points, and then
there is a unique way to lift~$F$ to a self-map~$\tilde F$ of~$M$ enjoying some of the
properties we are looking for. Unfortunately, this naive approach is too rough: the
manifold~$M$ constructed in this way is so large that many properties of the original
map~$F$ will be hidden inside the singular divisor~$\pi^{-1}(O)$.

To give an idea why this is the case (see Remark~3.3 for a more precise explanation), let
us discuss what is known about the local dynamics of~$F$ nearby the fixed point~$O$. In the
hyperbolic case (that is, when $dF_O$ has no eigenvalues of modulus one) the stable manifold
theorem (see, e.g., [Wu] for the statement in the complex case; see also~[S] and [R1, 2] for
the attracting case) describes completely the situation: there are two local
\hbox{$F$-invariant} manifolds, the stable one~$W^s$ and the unstable one~$W^u$, intersecting
transversally at the origin, such that~$(F|_{W^s})^k\to O$ and $(F|_{W^u})^{-k}\to O$ as
$k\to+\infty$, uniformly on compact sets. More generally, the local dynamics is
topologically conjugated to the dynamics induced by the differential~$dF_O$, with~$W^s$
corresponding to the direct sum of the generalized eigenspaces associated to eigenvalues
with modulus less than one, and $W^u$ corresponding to the direct sum of the
generalized eigenspaces associated to eigenvalues with modulus greater than one.

In the non-hyperbolic case, the theory at present is far less complete. One can recover a
good generalization of the classical one-variable Fatou-Leau theorem in the semi-attractive
case, when $dF_O$ has $1$ as eigenvalue of multiplicity one, and the others eigenvalues have
absolute value less than~1. In this case (studied first by Fatou~[F], and later by Ueda~[U1,
2] and Hakim~[H1]) either $F$ admits a holomorphic curve of fixed points passing through the
origin or there exists a basin of attraction to the origin, formed by~$k-1$ petals, where
$k\ge 2$ is the multiplicity of the origin as fixed point of~$F$; furthermore, Nishimura~[N]
has a description of the dynamics when there is a curve of fixed points.

Another situation that has been studied is when $dF_O=\id$, that is when $F$ is tangent to
the identity. In this case Hakim~[H2, 3] (see also Weickert~[W]) has proved that for $F$
generic there exists an
$F$-invariant stable (i.e., attracted to the origin) holomorphic curve with the origin {\it
in its boundary;} furthermore, there are estimates on the rate of approach of stable orbits
to the origin (see Section~3 for a precise statement of Hakim's results). Notice that, in
general, it is not possible to extend such a stable curve holomorphically through the
origin. It should also be mentioned that Rivi~[Ri] combined Hakim's results on maps tangent
to the identity with results on the semiattractive case to obtain a
description of the dynamics when there is a $dF_O$-invariant decomposition $\C^n=V_1\oplus
V_2$, with $dF_O|_{V_1}=\id$ and $\sp(dF_O|_{V_2})\subset\{|\lambda|<1\}$.

One feature that Hakim's and Weickert's works made clear is that one has to study orbits
converging to the origin tangentially to a given direction~$v\in\C^n$. It is easy to see
that such a~$v$ must be an eigenvector of~$dF_O$. Of course, not all the eigenvectors are
tangent to an orbit; but nevertheless this observation points out that, from a
dynamical point of view, the eigenvectors of $dF_O$ should be treated differently from the
non-eigenvectors. 

Now we can go back to our discussion of the manifold~$M$
in Theorem~\zuno. Blowing up points one deals with all the tangent directions in the
same way; and the previous discussion suggests that this should not be the case. The correct
replacement is blowing up {\it submanifolds;} in this way we are able to keep track of the
different status of the different tangent directions --- and we shall then be able to recover
easily informations about the local dynamics of~$F$ from informations about the dynamics
of~$\tilde F$ (see, e.g., Corollary~3.2).

In Section~1 we describe the canonical procedure for building
the manifold~$M$. It depends only on the Jordan block structure of the
differential~$dF_O$, and is obtained by blowing up a sequence of at most~$\mu_1+1$
submanifolds, where~$\mu_1$ is the dimension of the largest Jordan block in~$dF_O$. In
Section~2 we describe how to lift the map~$F$ to the blow-ups, and we give the
proof of Theorem~\zuno. It should be remarked that the construction is completely explicit;
for instance, it is possible to compute the local power series expansion of the lifted
map~$\tilde F$ in terms of the local power series expansion of~$F$, and this is essential
for the applications.

In Section~3 we apply the Diagonalization Theorem to dynamics. Since the eigenvalues
of~$d\tilde F_{\bf e}$ are quotients of the eigenvalues of~$dF_O$, this is really meaningful
only when all the eigenvalues of~$dF_O$ have modulus one. We shall concentrate on the case
$\sp(dF_O)=\{1\}$, because then $\tilde F$ is tangent to the identity. It turns
out that, for generic $F$, one and exactly one of the $\tilde F$-stable holomorphic
curves whose existence is guaranteed by Hakim's results is contained in
$M\setminus\pi^{-1}(O)$; its projection under~$\pi$ is then an $F$-stable
holomorphic curve, with the origin in its boundary (Corollary~3.2).

Thus we can apply Hakim's theory to generic maps~$F$ whose differential is
non-diagonalizable and such that $\sp(dF_O)=\{1\}$. Actually, our technique is flexible
enough to be used even for some classes of non-generic maps (see Section~3 for the
definition of ``generic'' in this context). For instance, we have fairly complete results in
the bi-dimensional case (Corollary~3.3), showing among other things that the dynamics might
depend strongly on the third degree terms of the map~$F$ even when the quadratic part is not
identically zero. Furthermore, we get yet another version of the Fatou-Bieberbach phenomenon
(Remark~3.7). 

A priori, one might suspect that other $\tilde F$-stable holomorphic curves might give rise 
at least to some other $F$-orbits converging to the origin, if not to $F$-stable holomorphic
curves. In the last section of this paper we shall show that, under some mild assumption on
the rate of convergence to zero of the orbit, if $dF_O$ is the canonical Jordan block~$J_n$
of order $n$ associated to~1 then this is not the case: roughly speaking, then, for such maps
the stable dynamics nearby the origin is described by Corollary~3.2.

I would like to end this introduction quoting a few lines from~[F, p.~135--137]: ``Ce cas
[that is, $n=2$ and $\sp(dF_O)=\{1\}$], tr\`es important au point de vue des applications
aux \'equations de la dynamique, exigerait de longues et difficiles recherches pour \^etre
\'elucid\'e compl\`etement. (...) Prenons par example (...) le cas limite
$$
\cases{x_1=x+\alpha y,\cr  y_1=y+a'x^2,\cr}
$$
substitution birationnelle que nous \'etudierons plus en detail dans la second partie de ce
M\'emoire''. Unfortunately, the promised second part never appeared; but now, after
seventy-five years, we are at last able to describe the dynamics of Fatou's example. 

\smallsect 1. The blow-up sequence

As described in the introduction, to diagonalize a non-diagonalizable dynamical system we
shall replace $\C^n$ by a suitable complex manifold obtained blowing-up a specific sequence
of submanifolds, depending on the Jordan block structure of the differential of the map
generating the dynamical system. In this section we introduce the general machinery
needed.

First of all we fix a number of notations. Given $0\le r<n$, a {\sl
splitting}~$\ca P$ of {\sl weight}~$r$ of~$n$ is a subdivision of $\{1,\ldots,n\}$ as a
disjoint union $\{1,\ldots,n\}=\ca P'\cup\ca P''$, where $\hbox{\rm card}\,\ca P'=r$ e
$\hbox{\rm card}\,\ca P''=n-r$. The {\sl standard splitting} of weight~$r$ is
$\{1,\ldots,r\}\cup\{r+1,\ldots,n\}$. If
$z=(z_1,\ldots,z_n)\in\C^n$ and $\ca P$ is a splitting of weight~$r>0$ with $\ca
P'=\{i_i,\ldots,i_r\}$ and $\ca P''=\{i_{r+1},\ldots,i_n\}$ (where $i_1<\cdots<i_r$ and
$i_{r+1}<\cdots<i_n$), we shall write
$z'=(z_{i_1},\ldots,z_{i_r})$ and $z''=(z_{i_{r+1}},\ldots,z_{i_n})$; if $r=0$ we
set $z''=z$, and $z'$ is empty. Finally, if $V$ is any vector space and $v\in
V\setminus\{O\}$, we denote by $[v]$ the projection of~$v$ in~$\P(V)$.

Let $M$ be a complex manifold of dimension $n\ge 2$, and $X\subset M$ a closed complex
submanifold of dimension $r\ge 0$. Let $N_{X/M}$ denote the normal bundle
of~$X$ in~$M$, and let $E_X=\P(N_{X/M})$ be the projective normal bundle, whose fiber
over $p\in X$ is $E_p=\P(T_pM/T_pX)$. The {\sl blow-up of~$M$ along~$X$} is the set
$$
\tilde M_X=(M\setminus X)\cup E_X\;,
$$
endowed with the manifold structure we shall presently describe, together with the
projection $\sigma\colon\tilde M_X\to M$ given by $\sigma|_{M\setminus X}
=\id_{M\setminus X}$ and $\sigma|_{E_p}\equiv\{p\}$ for $p\in X$. The set
$E_X=\sigma^{-1}(X)$ is the {\sl exceptional divisor} of the blow-up.

A chart $\phe=(z_1,\ldots,z_n)\colon V\to\C^n$ is {\sl adapted to~$X$} if there is a
splitting $\ca P$ of weight $r=\dim X$ such that $V\cap X=\{z''=0\}$. Choose a chart
$(V,\phe)$ adapted to~$X$, and for $j\in\ca P''$ and $q\in V\cap X$ set
$X_j=\{z_j=0\}\subset V$,
$L_{j,q}=\P\bigl(\Ker(dz_j)_q /T_qX\bigr)\subset E_q$, $L_j=\bigcup_{q\in V\cap X}L_{j,q}$,
$E_{V\cap X}=\sigma^{-1}(V\cap X)$ and $V_j=(V\setminus X_j)\cup (E_{V\cap X}\setminus L_j)$.
Define $\chi_j\colon V_j\to\C^n$ by
$$
\chi_j(q)_h=\cases{\phe(q)_h& if $h\in\ca P'$,\cr
	z_h(q)/z_j(q)& if $h\in\ca P''\setminus\{j\}$,\cr 
	z_j(q)&if $h=j$,\cr}
$$
if $q\in V\setminus X_j$, and by
$$
\chi_j([v])_h=\cases{\phe\bigl(\sigma([v])\bigr)_h&if $h\in\ca P'$,\cr
	d(z_h)_{\sigma([v])}(v)/d(z_j)_{\sigma([v])}(v)&if $h\in\ca P''\setminus\{j\}$,\cr
	0&if $h=j$,\cr}
$$	
if $[v]\in E_{V\cap X}\setminus L_j$. Then it is not difficult to check that the
charts $(V_j,\chi_j)$, together with an atlas of~$M\setminus X$, endow $\tilde M_X$ with a
structure of $n$-dimensional complex manifold, as claimed, such that the projection~$\sigma$
is holomorphic everywhere. For future reference, we record here that
$$
\phe\circ\sigma\circ\chi_j^{-1}(w)_h=\cases{w_h&if $h\in\ca P'\cup\{j\}$,\cr
	w_jw_h&if $h\in\ca P''\setminus\{j\}$.\cr}
\neweq\eqdzero
$$

The fiber~$E_p$ of the exceptional divisor over a point~$p\in X$ is a projective space; so
the choice of an adapted chart yields an explicit isomorphism with~$\P^{n-r-1}(\C)$ that we
shall denote by $\iota_{p,\phe}\colon E_p\to\P^{n-r-1}(\C)$. Finally, if
$Y\subseteq M$ is a submanifold of~$M$, then the {\sl proper transform} of~$Y$ is
$\tilde Y=\overline{\sigma^{-1}(Y\setminus X)}\subset\tilde M_X$. 

To describe the sequence of blow-ups we need some more notations. Given $\rho\ge 1$, a {\sl
$\rho$-partition of~$n$} is a set $\ca M=\{\mu_1,\ldots,\mu_\rho\}\subset\N$ with
$\mu_1\ge\cdots\ge\mu_\rho\ge 1$ and $\mu_1+\cdots+\mu_\rho=n$. The {\sl length} $\ell(\ca
M)$ of~$\ca M$ is~$\mu_1$ if~$\mu_2<\mu_1$, and~$\mu_1+1$ if $\mu_2=\mu_1$. 

To a $\rho$-partition~$\ca M$ we can associate several objects. First of all, we define
$\nu_1,\ldots,\nu_\rho\in\N$ by setting $\nu_1=0$ and $\nu_j=\nu_{j-1}+\mu_{j-1}$ for
$j=2,\ldots,\rho$. Then we define sets $\ca P'_{kl}\subset\{1,\ldots,n\}$ for $0\le
k\le\mu_1-1$ and $1\le l\le\rho$ by setting
$$
\ca P'_{kl}=\cases{\void&if $k=0$,\cr
	\{\nu_l+1,\ldots,\nu_l+\min(k,\mu_l)\}& if $1\le k\le\mu_1-1$.\cr}
$$
If $\mu_2=\mu_1$, we also define $\ca P'_{\mu_1,l}$ for $1\le l\le\rho$ by
$$
\ca P'_{\mu_1,\l}=\cases{\{\nu_l+1,\ldots,\nu_l+\mu_l\}& if $l\neq 2$,\cr
	\{\nu_2+1,\ldots,\nu_2+\mu_2-1\}&if $l=2$;\cr}
$$
we also set $\ca P'_{\mu_1+1,1}=\{1,\ldots,\mu_1,\nu_2+\mu_2\}$.

Then we get $\ell(\ca M)$ splittings~$\ca P_k$ of~$n$ by setting $\ca P'_k=\bigcup_{l=1}^\rho
\ca P_{kl}$ and $\ca P''_k=\{1,\ldots,n\}\setminus\ca P'_k$. Furthermore, we also get a
sequence of linear subspaces $\void=Y^0\subset Y^1\subset\cdots\subset
Y^{\ell(\ca M)-1}\subset\P^{n-1}(\C)$ by letting $Y^k$ to be the subspace generated by
$\{[e_h]\mid h\in\ca P'_k\}$, where $\{e_1,\ldots,e_n\}$ is the canonical basis of~$\C^n$.

We are now ready to associate a sequence of $\ell(\ca M)$ blow-ups to any $\rho$-partition
$\ca M$ of~$n$. Set $M^0=\C^n$, $\chi_0=\id_{\C^n}$, ${\bf
e}_0=O$ and~$X^0=\{O\}$. We start by blowing up the origin, taking $M^1=\tilde M^0_{X^0}$
and $\pi_1=\sigma_1\colon M^1\to M^0$. Since $M^0=\C^n$ has a canonical chart adapted to
$X^0$ (that is, centered at the origin), the exceptional divisor~$E^1=\pi_1^{-1}(X^0)$ is
canonically isomorphic to
$\P^{n-1}(\C)$. This allows us to define a distinguished point~${\bf e}_1\in E^1$,
corresponding to~$[e_1]\in\P^{n-1}(\C)$, and also distinguished linear subspaces $Y^k\subset
E^1$ for $k=1,\ldots,\ell(\ca M)-1$, corresponding to the previously defined linear
subspaces of~$\P^{n-1}(\C)$ associated to~$\ca M$.

Now put $X^1=Y^1$ and set $M^2=\tilde M^1_{X^1}$. Let $X^2\subset M^2$ be the
proper transform of~$Y^2$, and set~$M^3=\tilde M^2_{X^2}$. Next, let $X^3\subset M^3$ be the
proper transform (with respect to $\sigma_3\colon M^3\to M^2$) of the proper transform (with
respect to $\sigma_2\colon M^2\to M^1$) of~$Y^3$, and put $M^4=\tilde M^3_{X^3}$. Proceeding
in this way, we define for $k=2,\ldots,\ell(\ca M)-1$ the manifold~$M^{k+1}$ as the blow-up
of~$M^k$ along the iterated proper transform~$X^k$ of~$Y^k$; we denote by
$\sigma_{k+1}\colon M^{k+1}\to M^k$ the associated projection, and by
$E^{k+1}=\sigma^{-1}_{k+1}(X^k)\subset M^{k+1}$ the exceptional divisor. For
$k=1,\ldots,\ell(\ca M)$ we also put $\pi_k=\sigma_1\circ\cdots\circ\sigma_k\colon M^k\to
M^0$; the set $\pi_k^{-1}(X^0)$ will be called the {\sl singular divisor} of~$M^k$.

At each stage of this construction there are canonical charts adapted to the submanifolds
involved:

\newthm Lemma \quno:
For $1\le k\le \ell(\ca M)$ we can find a distinguished point ${\bf e}_k\in M^k$ and a
canonical chart $(V_k,\chi_k)$ centered in~${\bf e}_k$ such that:
$$
V_k\cap X^k=\chi_k^{-1}\left(\{w_1=0\}\cap\bigcap_{h\in\ca P''_k}\{w_h=0\}\right);
\neweq\eqduno
$$
$$
V_k\cap\pi_k^{-1}(X^0)=\chi_k^{-1}\left(\bigcup_{h\in\ca P'_{k1}}\{w_h=0\}\right)\supset
V_k\cap X^k;
\neweq\eqddue
$$
and such that for $h=k+1,\ldots,\ell(\ca M)-1$ the intersection of $V_k$ with the iterated
proper transform of~$Y^h$ is
$$
\chi_k^{-1}\left(\{w_1=0\}\cap\bigcap_{h\in\ca P''_h}\{w_h=0\}\right).
$$
Furthermore, $\chi_0\circ\sigma_1\circ\chi_1^{-1}(w)=(w_1,w_1w_2,\ldots,w_1w_n)$, 
$\chi_{\mu_1}\circ\sigma_{\mu_1+1}\circ\chi_{\mu_1+1}^{-1}(w)=(w_1w_{\nu_2+\mu_2},w_2,
\ldots,w_n)$, and for $2\le k\le\mu_1$
$$
\chi_{k-1}\circ\sigma_k\circ\chi_k^{-1}(w)_h=\cases{w_h&if $h\in(\ca P'_{k-1}\setminus\{1\})
	\cup\{k\}$,\cr
	w_kw_h&if $h\in\{1\}\cup(\ca P''_{k-1}\setminus\{k\})$.\cr}
\neweq\eqdtre
$$

\pf
For $k=1$, the existence of a canonical chart adapted to~$X^0$ yields a canonical
chart~$(V_1,\chi_1)$ centered at~${\bf e}_1$ and adapted to~$X^1$; in turn this yields
a canonical basis $\{\de/\de w_1,\ldots,\de/\de w_n\}$ of~$T_{{\bf e}_1}M^1$. Furthermore, it
is easy to check that
$$
V_1\cap E^1=\chi_1^{-1}(\{w_1=0\})=V_1\cap\pi_1^{-1}(X^0)\supset
V_1\cap X^1=\chi_1^{-1}\left(\{w_1=0\}\cap\bigcap_{h\in\ca P''_1}\{w_h=0\}\right),
$$
and that
$$
\chi_0\circ\sigma_1\circ\chi_1^{-1}(w)=(w_1,w_1w_2,\ldots,w_1w_n).
$$
So the lemma is proved for $k=1$.

Assume, by induction, that the lemma holds for $k-1$. In particular, we have a
distinguished point~${\bf e}_{k-1}$ and a canonical chart $(V_{k-1},\chi_{k-1})$ centered
at~${\bf e}_{k-1}$ and adapted to~$X^{k-1}$. We thus
have a canonical basis $\{\de/\de w_1,\ldots,\de/\de w_n\}$ of $T_{{\bf e}_{k-1}}M^{k-1}$
such that $\{\de/\de w_h\mid h\in\ca P'_{k-1}\setminus\{1\}\}$ spans
$T_{{\bf e}_{k-1}}X^{k-1}$. Put
$$
{\bf e}_k=\left[{\de\over\de w_k}+T_{{\bf e}_{k-1}}X^{k-1}\right]\in\sigma_k^{-1}({\bf
		e}_{k-1})\;,
$$
(or ${\bf e}_k=\left[{\de\over\de w_{\nu_2+\mu_2}}+T_{{\bf
e}_{\mu_1}}X^{\mu_1}\right]\in\sigma_{\mu_1}^{-1}({\bf e}_{\mu_1})$ if $k=\mu_1+1$),
and let $(V_k,\chi_k)$ be the canonical chart centered in~${\bf e}_k$ constructed, as before,
via~$(V_{k-1},\chi_{k-1})$. Then it is not too difficult to check using the inductive
hypothesis that $(V_k,\chi_k)$ is as desired.\qedn

We end this section by remarking that it is easy to prove by induction that if we fix $1\le
k\le\ell(\ca M)$ and write $z=\chi_0\circ\pi_k\circ\chi_k^{-1}(w)$ then
$$
\eqalign{z_j&=\cases{w_1\prod\limits_{h=2}^j (w_h)^2\,\prod\limits_{h=j+1}^k w_h& if
	$j\in\ca P'_{k1}$,\cr
	w_1\prod\limits_{h=2}^{j-\nu_l} (w_h)^2\left(\prod\limits_{h=j-\nu_l+1}^k w_h\right)
	w_j& if $j\in\ca P'_{kl}$, $2\le l\le\rho$;\cr
	w_1\prod\limits_{h=2}^k (w_h)^2\, w_j&if $j\in\ca P''_k$;\cr}\qquad
	\hbox{if $1\le k\le\mu_1$;}\cr
z_j&=\cases{w_1\prod\limits_{h=2}^j (w_h)^2\left(\prod\limits_{h=j+1}^{\mu_1}w_h\right)
	w_{\nu_2+\mu_2}& if $j\in\ca P'_{\mu_1,1}$,\cr
	w_1\prod\limits_{h=2}^{j-\nu_l} (w_h)^2\left(\prod\limits_{h=j-\nu_l+1}^{\mu_1} w_h
	\right)w_jw_{\nu_2+\mu_2}& if
	$j\in\ca P'_{\mu_1,l}$, $2\le l\le\rho$;\cr
	w_1\prod\limits_{h=2}^{\mu_1}(w_h)^2\, (w_{\mu_2+\mu_2})^2&if
		$j\in\ca P''_{\mu_1}$;\cr}\qquad
	\hbox{if $k=\mu_1+1$.}\cr}
\neweq\equcin
$$
Furthermore, if $z_1,\ldots,z_k\neq 0$ then
$$
\eqalign{w_j=\cases{(z_1)^2/z_k& if $j=1$,\cr
	z_j/z_{j-1}& if $j\in\ca P'_{k1}\setminus\{1\}$,\cr
	z_j/z_{j-\nu_l}& if $j\in\ca P'_{kl}$, $2\le l\le\rho$,\cr
	z_j/z_k& if $j\in\ca P''_k$;\cr}\qquad\quad&\hbox{if $1\le k\le\mu_1$;}\cr
w_j=\cases{(z_1)^2/z_{\nu_2+\mu_2}& if $j=1$,\cr
	z_j/z_{j-1}& if $j\in(\ca P'_{\mu_1,1}\setminus\{1\})$,\cr
	z_j/z_{j-\nu_l}& if $j\in\ca P'_{\mu_1,l}$, $2\le l\le\rho$,\cr
	z_j/z_{\mu_1}& if $j\in\ca P''_{\mu_1}$.\cr}\qquad\quad&\hbox{if $k=\mu_1+1$.}\cr}
\neweq\equsei
$$

\vfill\eject
\smallsect 2. The diagonalization theorem

We shall denote by $\End(\C^n,O)$ the set of germs of holomorphic self-maps of $\C^n$
sending the origin~$O$ to itself; more generally, if $X$ is a closed set of a complex
manifold~$M$, we shall denote by $\End(M,X)$ the set of germs at $X$ of holomorphic
self-maps of~$M$ sending $X$ into itself. Every germ $F\in\End(\C^n,O)$ has a 
homogeneous expansion of the form
$$
F(z)=\sum_{j=1}^\infty P_j(z),
$$
where $z=(z_1,\ldots,z_n)\in\C^n$, and the $P_j$'s are $n$-uples of homogeneous polynomials
of degree~$j$ in~$z_1,\ldots,z_n$.

Let $M$ be a complex manifold of dimension~$n$, and $X$ a closed submanifold of
dimension~$r\ge 0$. We are interested to see when a germ
$F\in\End(M,X)$ can be lifted to the blow-up~$\tilde M_X$ as a germ $\tilde
F\in\End(M_X,E_X)$. Take $p\in X$, and choose charts $(V,\phe)$ and
$(\tilde V,\tilde\phe)$ adapted to~$X$ so that $p\in V$ and $F(p)\in\tilde V$. In a
neighbourhood of~$p$ we can write the homogeneous expansion of $G=\tilde\phe\circ
F\circ\phe^{-1}$ as
$$
G(z)=\sum_{l\ge 0}P_{l,z'}(z'')\;,
$$
where $P_{l,z'}$ is a $n$-uple of $l$-homogeneous polynomials with coefficients holomorphic
in~$z'$. The condition $F(X)\subseteq X$ then translates to
$$
(P_{0,z'})''\equiv 0\;.
$$
The {\sl order of~$F$ at~$p$ along~$X$} is 
$$
\nu_X(F,p)=\min\{l\mid (P_{l,\phe(p)'})''\not\equiv 0\}\ge 1\;;
$$
it is easily checked that $\nu_X(F,p)$ does not depend on the adapted charts chosen. The {\sl
order of~$F$ along~$X$} is then given by
$$
\nu_X(F)=\min\{\nu_X(F,p)\mid p\in X\}\;.
$$
Clearly the set $\{p\in X\mid \nu_X(F,p)=\nu_X(F)\}$ is open in~$X$.

We shall say that $F$ is {\sl non-degenerate at~$p$ along~$X$} if 
\smallskip

\item{(i)} $F^{-1}(p)\subseteq X$,
\item{(ii)} $\nu_X(F,p)=\nu_X(F)$, and
\item{(iii)} $\bigl(P_{l_0,\phe(p)'}(v)\bigr)''=0$ iff $v=O\in\C^{n-r}$, where
$l_0=\nu_X(F)$. 

\medskip
\noindent If $F$ is non-degenerate along~$X$ at all points of~$X$ we shall say that $F$ is
{\sl non-degenerate along~$X$.} 

\newthm Proposition \tuno:
Let $M$ be a complex manifold of dimension $n$, and $X\subset M$ a closed submanifold of
dimension~$r\ge 0$. Let $F\in\End(M,X)$ be non-degenerate along $X$. Then there exists a
unique $\tilde F\in\End(\tilde M_X,E_X)$ such that $F\circ\sigma=\sigma\circ\tilde F$.
Furthermore, if $p\in X$ and $(V,\phe)$, $(\tilde V,\tilde\phe)$ are charts adapted to~$X$
with~$p\in V$ and~$F(p)\in\tilde V$, then
$$
\tilde F\bigl([v]\bigr)=(\iota_{F(p),\tilde\phe})^{-1}\left(\left[P_{l_0,\phe(p)'}
	\bigl(\iota_{p,\phe}([v])\bigr)''\right]\right)
\neweq\eqtuno
$$
for all $[v]\in E_p$, where $l_0=\nu_X(F)$.

\pf
Since $F^{-1}(X)\subseteq X$, if $q$ does not belong to~$X$ we can safely set $\tilde
F(q)=F(q)$; we are left to define $\tilde F$ on the exceptional divisor.

Choose $p\in X$, and the charts as in the statement of the theorem; without loss of
generality, we can assume that for both charts the associated splitting is the standard one.
For
$[v]\in E_p$ choose
$r+1\le j\le n$ so that $[v]\in V_j$; if $\tilde F$ exists, we must have
$$
F\circ\sigma\circ\chi_j^{-1}=\sigma\circ\tilde F\circ\chi_j^{-1}\;.
$$
If $[v]=(\iota_{p,\phe})^{-1}[v_{r+1}:\ldots:v_n]$, we have
$$
[v]=\lim_{\zeta\to 0}\chi_j^{-1}\left(\phe(p)',{v_{r+1}\over v_j},\ldots,\zeta,\ldots,
	{v_n\over v_j}\right)\;,
$$
and so, setting again $G=\tilde\phe\circ F\circ\phe^{-1}$, 
$$
\tilde F([v])=\lim_{\zeta\to 0}\sigma^{-1}
\left(\tilde\phe^{-1}\left(G\left(\phe(p)',{\zeta\over v_j}v\right)\right)\right)\;,
$$
where with a slight abuse of notation we have put $v=(v_{r+1},\ldots,v_n)\in\C^{n-r}$. 

Now, given a sequence $\{q_k\}\subset M\setminus X$ converging
to~$q\in X$, the sequence~$\{\sigma^{-1}(q_k)\}$ converges in~$\tilde M
\setminus X$ iff $\{[\tilde\phe(q_k)'']\}$ converges in $\P^{n-r-1}(\C)$, and then
$$
\lim_{k\to\infty}\sigma^{-1}(q_k)=\iota_{q,\tilde\phe}^{-1}\left(\lim_{k\to\infty}[\tilde
	\phe(q_k)'']\right).
$$
In our case we have
$$
G\left(\phe(p)',{\zeta\over v_j}v\right)''=\sum_{l\ge l_0}P_{l,\phe(p)'}\left(
	{\zeta\over v_j}v\right)''=\left(\zeta\over v_j\right)^{l_0}\bigl(P_{l_0,\phe(p)'}(v)''+
	\zeta Q(\zeta)\bigr)\;,
$$
for a suitable holomorphic map $Q$. Therefore $[G(\phe(p)',\zeta
v/v_j)'']\to[P_{l_0,\phe(p)'}(v)'']$, and thus if $\tilde F$ exists it is given by~\eqtuno\
on the exceptional divisor.

To finish the proof we must show that an $\tilde F$ defined by \eqtuno\ on the exceptional
divisor and by $F$ elsewhere is holomorphic. Take $[v]\in E_p$, and choose $r+1\le h,k\le n$
so that $[v]\in V_h$ and $\tilde F([v])\in\tilde V_k$; we must show that $\chi_k\circ\tilde
F\circ\chi_h^{-1}$ is holomorphic. We know that
$$
G\circ(\phe\circ\sigma\circ\chi_h^{-1})=(\phe\circ\sigma\circ\chi_k^{-1})\circ
	(\chi_k\circ\tilde F\circ\chi_h^{-1})\;;
$$
so putting $\chi_k\circ\tilde F\circ\chi_h^{-1}=(\tilde f_1,\ldots,\tilde f_n)$ and
recalling \eqdzero\ we must have
$$
G(w',w_hw_{r+1},\ldots,w_h,\ldots,w_hw_n)
	=\bigl(\tilde f_1(w),\ldots,\tilde f_r(w),
	\tilde f_k(w)\tilde f_{r+1}(w),\ldots,\tilde f_k(w),\ldots,\tilde f_k(w)\tilde f_n(w)
	\bigr)\;.
$$
Writing $G=(g_1,\ldots,g_n)$ we find that if $w_h\neq 0$ then
$$
\tilde f_i(w)=\cases{g_i(w',w_hw_{r+1},\ldots,w_h,\ldots,w_hw_n)& if $1\le i\le r$ or
	$i=k$,\cr
	\displaystyle {g_i(w',w_hw_{r+1},\ldots,w_h,\ldots,w_hw_n)\over
	g_k(w',w_hw_{r+1},\ldots,w_h,\ldots,w_hw_n)}& if $r+1\le i\neq k\le n$.\cr}
\neweq\eqtdue
$$
Since the $g_i$'s are holomorphic and $\{w_h=0\}$ has codimension~1 in $\chi_h(V_h)$, to
end the proof it suffices to show that the quotients in \eqtdue\ have a limit when
$w\to\chi_h([v])$. 

Write again $\iota_{p,\phe}([v])=[v_{r+1}:\ldots:v_n]$ and $v=(v_{r+1},\ldots,v_n)$,
and assume then that $w\to\chi_h([v])$. This means that $w'\to\phe(p)'$,
$w_h\to 0$ and $(w_{r+1},\ldots,1,\ldots,w_n)\to v_h^{-1}v$. Now,
$$
G(w',w_hw_{r+1},\ldots,w_h,\ldots,w_hw_n)''
	=\sum_{l\ge l_0}\left({w_h\over v_h}\right)^l P_{l,w'}(w_{r+1}v_h,\ldots,v_h,\ldots,
	w_n v_h)''\;.
$$
Since $\tilde F([v])\in\tilde V_k$ we have $P_{l_0,\phe(p)'}(v)_k\neq 0$; therefore
$$
{g_i(w',w_hw_{r+1},\ldots,w_h,\ldots,w_hw_n)\over
	g_k(w',w_hw_{r+1},\ldots,w_h,\ldots,w_hw_n)}\to {P_{l_0,\phe(p)'}(v)_i\over
	P_{l_0,\phe(p)'}(v)_k}\;,
$$
and we are done.\qedn

Now, our construction involves iterated blow-ups; thus we are
interested to know when the map $\tilde F$ is still non-degenerate along suitable
submanifolds of~$\tilde M_X$. We shall limit ourselves to two special cases, which are enough
for our aims. 

\newthm Proposition \tdue:
Let $M$ be a complex manifold of dimension $n$, and $X\subset M$ a closed submanifold of
dimension~$r\ge 0$. Let $F\in\End(M,X)$ be non-degenerate along $X$, and $\tilde
F\in\End(\tilde M_X,E_X)$ its lifting. Let~$Y\subseteq M$ be a submanifold of~$M$ of
dimension~$r+s$ (with $s\ge1$), and $\tilde Y\subseteq\tilde M$ its proper transform.
Assume that
\medskip
{
\item{\rm (i)} $Y$ contains properly $X$;
\item{\rm (ii)} $F(Y)\subseteq Y$ and $F^{-1}(Y)\subseteq Y$;
\item{\rm (iii)} $dF_q$ is invertible for all~$q\in Y$.
\medskip
\noindent Then $\tilde F$ is non-degenerate along $\tilde Y$, and $d\tilde F_{\tilde q}$ is
invertible for all~$\tilde q\in\tilde Y$.}

\pf
First of all, notice that if $p\in X$ then $\tilde Y\cap E_p=\P(T_pY/T_pX)$, and that $\tilde
F|_{E_p}$ is induced by $dF_p$. Since, by construction, $\tilde F(\tilde Y)\subseteq\tilde Y$
and $\tilde F^{-1}(\tilde Y\setminus E_X)\subseteq\tilde Y\setminus E_X$, it suffices to
prove that $d\tilde F_{[v]}$ is invertible for all $[v]\in\tilde Y\cap E_X$. 

Fix $p\in X$ and $[v]\in\tilde Y\cap E_p$, and choose two charts $(V,\phe)$ and $(\tilde
V,\tilde\phe)$, centered in~$p$, respectively in~$F(p)$, such that
$V\cap X=\{z_{r+1}=\cdots=z_n=0\}$, $V\cap Y=\{z_{r+s+1}=\cdots=z_n=0\}$,
and analogously for $\tilde V$. In particular,
$$
\iota_{p,\phe}(\tilde Y\cap E_p)=\iota_{F(p),\tilde\phe}(\tilde Y\cap E_{F(p)})=
	\{v_{r+s+1}=\cdots=v_n=0\}\;,
$$
and we can also assume that $\iota_{p,\phe}([v])=\iota_{F(p),\tilde\phe}\bigl(\tilde F([v])
\bigr)=[1:0:\cdots:0]$. Then the charts $(V_{r+1},\chi_{r+1})$ and $(\tilde V_{r+1},
\tilde\chi_{r+1})$ are centered in~$[v]$, respectively in~$\tilde F([v])$, and adapted
to~$\tilde Y$.

Set $G=\tilde\phe\circ F\circ\phe^{-1}=(g_1,\ldots,g_n)$ and $\tilde G=\tilde\chi_{r+1}
\circ\tilde F\circ\chi_{r+1}^{-1}=(\tilde f_1,\ldots,\tilde f_n)$; the relation between
the~$g_i$'s and the $\tilde f_j$'s is given by~\eqtdue. 
Since $F(X)\subseteq X$ and $F(Y)\subseteq Y$, the jacobian matrix of $G$ at the origin is
of the form
$$
\ca A=\left|\vcenter{\offinterlineskip
\halign{&\hfil$#$\hfil\cr
\multispan3\strut&\omit\hskip1.5mm\vrule width 1pt\hskip 1.5mm&\cr
\multispan3\strut\hskip1.5mm\hfil $A$\hfil&\omit\hskip1.5mm\vrule width 1pt
	\hskip 1.5mm&\multispan3\hfil$*$\hskip1.5mm\hfil\cr
\multispan3\strut&\omit\hskip1.5mm\vrule width 1pt\hskip 1.5mm&\cr
\noalign{\hrule height 1pt}
\multispan3&\omit\hskip1.5mm\vrule height 1pt width 1pt\hskip 1.5mm&\multispan3\cr
\strut&&&\omit\hskip1.5mm\vrule width 1pt\hskip 1.5mm&B&\omit\hskip1.5mm
	\vrule\hskip1.5mm&*\hskip1.5mm\cr
\multispan3\hskip1.5mm\quad\hfil$\smash{\lower 3pt\vbox{\smash{$O$}}}$\quad
	\hfil&\omit\hskip1.5mm\vrule width 1pt&\multispan3\hskip-1mm\hrulefill\cr
\multispan3&\omit\hskip1.5mm\vrule height 1pt width 1pt\hskip 1.5mm&\omit&\omit
\hskip1.5mm\vrule height 1pt\hskip1.5mm&\omit\cr
\strut&&&\omit\hskip1.5mm\vrule width 1pt\hskip 1.5mm&O&\omit\hskip1.5mm\vrule
	\hskip1.5mm&C\hskip1.5mm\cr}}
\right|\;,
$$
with $A\in M_{r,r}(\C)$, $B\in M_{s,s}(\C)$ and $C\in M_{n-r-s,n-r-s}(\C)$. Since, by
assumption, $dF_p$ is invertible, we have
$$
\det(\ca A)=\det(A)\det(B)\det(C)\neq 0\;.
$$
Finally, $\tilde F([v])\in\tilde V_{r+1}$ translates in
$$
\lambda={\de g_{r+1}\over\de z_{r+1}}(O)\neq 0\;.
$$
Our aim is to compute $\de\tilde f_i/\de w_j$ at $w=O$. This is easy when $1\le i\le r+1$;
in fact, \eqtdue\ with $h=k=r+1$ yields
$$
{\de\tilde f_i\over\de w_j}(O)=\cases{\displaystyle{\de g_i\over\de z_j}(O)& for $1\le i\le
	r+1$, $1\le j\le r+1$,\cr
	0& for $1\le i\le r+1$, $r+2\le j\le n$.\cr}
$$
In particular,
$$
{\de\tilde f_{r+1}\over\de w_j}(O)=\cases{0& if $j\neq r+1$,\cr
	\lambda\neq 0& if $j=r+1$.\cr}
$$
Now set $\tilde g_i(w)=g_i(w',w_{r+1},w_{r+1}w_{r+2},\ldots,
w_{r+1}w_n)$, and write again
$$
G(z)=\sum_{l\ge 0}P_{l,z'}(z'')\;,
$$
recalling that $(P_{0,z'})''\equiv O$. For $r+2\le i\le n$ we have
$$
{\de\tilde f_i\over\de w_j}(O)=\lim_{w\to O}{1\over\tilde g_{r+1}(w)}\left[
	{\de\tilde g_i\over\de w_j}(w)-{\tilde g_i(w)\over\tilde g_{r+1}(w)}{\de\tilde g_{r+1}
	\over\de w_j}(w)\right]\;.
\neweq\eqttre
$$
Since
$$
\tilde g_i(w)=\sum_{l\ge 0}(w_{r+1})^l P_{l,w'}(1,w_{r+2},\ldots,w_n)_i\;,
$$
\eqttre\ yields
$$
{\de\tilde f_i\over\de w_j}(O)=\cases{\displaystyle {1\over\lambda}\left[{\de^2 g_i\over\de
z_j\de z_{r+1}}(O)-{1\over\lambda}{\de^2 g_{r+1}\over\de z_j\de z_{r+1}}(O){\de g_i\over
	\de z_{r+1}}(O)\right]&for $r+2\le i\le n$ and $1\le j\le r+1$,\cr
\noalign{\smallskip}
\displaystyle{1\over\lambda}\left[{\de g_i\over\de
	z_j}(O)-{1\over\lambda} {\de g_{r+1}\over\de z_j}(O){\de g_i\over\de z_{r+1}}(O)\right]&for
	$r+2\le i,j\le n$.\cr}
$$
In particular, we find
$$
{\de\tilde f_i\over\de w_j}(O)={1\over\lambda}{\de g_i\over\de z_j}(O)\quad\hbox{for
	$r+s+1\le i\le n$, $r+2\le j\le n$.}
$$
Summing up, we have proved that the Jacobian matrix of $\tilde G$ at the origin is
$$
\tilde{\ca A}=\left|\vcenter{\offinterlineskip
\halign{&\hfil$#$\hfil\cr
\multispan3&\omit\hskip1.5mm\vrule height2mm\hskip 1.5mm&\omit&\omit
	\hskip1.5mm\vrule height2mm width 1pt\hskip 1.5mm\cr
\multispan3\strut\hskip1.5mm\hfil $A$\hfil&\omit\hskip1.5mm\vrule\hskip	1.5mm&*&\omit
	\hskip1.5mm\vrule width 1pt\hskip 1.5mm&\multispan3\hfil$O$\hskip1.5mm\hfil\cr
\multispan3&\omit\hskip1.5mm\vrule height2mm\hskip 1.5mm&\omit&\omit
	\hskip1.5mm\vrule height2mm width 1pt\hskip 1.5mm\cr
\noalign{\hrule}
\multispan3&\omit\hskip1.5mm\vrule height2pt\hskip 1.5mm&\omit&\omit
	\hskip1.5mm\vrule height2pt width 1pt\hskip 1.5mm\cr
\multispan3\strut\hskip1.5mm\hfil $O$\hfil&\omit\hskip1.5mm\vrule\hskip	1.5mm&\lambda&\omit
	\hskip1.5mm\vrule width 1pt\hskip 1.5mm&\multispan3\hfil$O$\hskip1.5mm\hfil\cr
\noalign{\hrule height 1pt}
\multispan3&\omit\hskip1.5mm\vrule height 2pt\hskip 1.5mm&\omit&\omit\hskip1.5mm\vrule
		height 2pt width 1pt\hskip 1.5mm&\multispan3\cr
\strut&&&\omit\hskip1.5mm\vrule \hskip 1.5mm&\omit&\omit\hskip1.5mm\vrule width 1pt\hskip
	1.5mm&\tilde B&\omit\hskip1.5mm\vrule\hskip1.5mm&\hfil$\smash{\raise
	2pt\vbox{\smash{$*$}}}$\hfil\hskip1.5mm\cr
\multispan3\hskip1.5mm\quad\hfil$\smash{\lower 3pt\vbox{\smash{$*$}}}$\quad
	\hfil&\omit\hskip1.5mm\vrule&\hfil$\smash{\lower 3pt\vbox{\smash{$*$}}}$
	\hfil&\omit\hskip1.5mm\vrule width 1pt&\multispan3\hskip-1mm\hrulefill\cr
\multispan3&\omit\hskip1.5mm\vrule height 1pt\hskip 1.5mm&\omit&\omit\hskip1.5mm\vrule
		height 1pt width 1pt\hskip 1.5mm&\omit&\omit\hskip1.5mm\vrule height 1pt\hskip 1.5mm&\cr
\strut&&&\omit\hskip1.5mm\vrule\hskip 1.5mm&\omit&\omit\hskip1.5mm\vrule width 1pt\hskip
	1.5mm&O&\omit\hskip1.5mm\vrule \hskip1.5mm&{1\over\lambda}C\hskip1.5mm\cr}}
\right|\;,
\neweq\eqtqua
$$
where $\tilde B\in M_{s-1,s-1}(\C)$. Now, if we subtract to the $j$-th column of $B$ (for
$j=2,\ldots,s$) the first column of~$B$ multiplied by $\lambda^{-1}\de g_{r+1}/\de
z_{r+j}(O)$ we get
$$
\left|\vcenter{\offinterlineskip
\halign{&\strut\hfil$#$\hfil\cr
\hskip1.5mm\lambda&\omit\hskip1.5mm\vrule\hskip1.5mm&O\hskip1.5mm\cr
\noalign{\hrule}
\omit&\omit\hskip 1.5mm \vrule height2pt\hskip 1.5mm&\omit\cr
\hskip1.5mm*&\omit\hskip1.5mm\vrule\hskip1.5mm&\lambda\tilde B\hskip1.5mm\cr}}
\right|\;.
$$
Since these elementary operations do not change the determinant, we obtain
$\det(B)=\lambda^s\det(\tilde B)$. Therefore
$$
\det(\tilde{\ca A})={1\over\lambda^{n-r-1}}\det(\ca A)\neq 0\;,
$$
and we are done.\qedn

A similar argument yields:

\newthm Proposition \ttre:
Let $M$ be a complex manifold of dimension $n$, and $X\subset M$ a closed submanifold of
dimension~$r\ge 0$. Let $F\in\End(M,X)$ be non-degenerate along $X$, and $\tilde
F\in\End(\tilde M_X,E_X)$ its lifting. Take $p\in X$ and a linear subspace $L\subseteq
E_p$ of dimension~$s-1$ (with $s\ge1$). Assume that
{\medskip
\item{\rm (i)} $\tilde F(L)\subseteq L$, and
\item{\rm (ii)} $dF_p$ is invertible.
\medskip
\noindent Then $\tilde F$ is non-degenerate along $L$, and $d\tilde F_{[v]}$ is
invertible for all~$[v]\in L$.}

\pf
Condition (i) implies that $p$ is a fixed point of $F$, and condition (ii) implies that
$\nu_X(F)=1$. In particular, $\tilde F|_{E_p}$ is induced by the differential of~$F$ at~$p$;
thus $\tilde F|_L$ is injective, and the invertibility of $d\tilde F_{[v]}$ for all~$[v]\in
L$ will imply that $\tilde F$ is non-degenerate along~$L$.

Fix $[v]\in L$, and choose two charts $(V,\phe)$, $(\tilde V,\tilde\phe)$ centered in $p$
adapted to~$X$ such that 
$$
\iota_{p,\phe}([v])=\iota_{p,\tilde\phe}\bigl(\tilde F([v])\bigr)=[1:0:\ldots:0]
$$ 
and
$$
\iota_{p,\phe}(L)=\iota_{p,\tilde\phe}(L)=\{v_{r+s+1}=\cdots=v_n=0\}\;.
$$
Then the charts $(V_{r+1},\chi_{r+1})$ and $(\tilde V_{r+1},
\tilde\chi_{r+1})$ are centered in~$[v]$, respectively in~$\tilde F([v])$, and adapted
to~$L$. The proof then goes on as in the previous
proposition.\qedn

%
%

We are finally ready to prove the main result of this paper:

\newthm Theorem \qtre: (Diagonalization Theorem)
Let $F\in\End(\C^n,O)$ be such that $dF_O$ is invertible and non-diagonalizable. Assume that
$dF_O$ is in Jordan canonical form, with $\rho\ge1$ blocks of lenghts
$\mu_1\ge\cdots\ge\mu_\rho\ge1$ associated respectively to the eigenvalues
$\lambda_1,\ldots,\lambda_\rho\in\C$. Set $\ca
M=\{\mu_1,\ldots,\mu_\rho\}$, and let $(M^0,\ldots,M^{\ell(\ca M)})$ be the sequence of
blow-ups associated to~$\ca M$. Then for $1\le k\le\ell(\ca M)$ there exists a unique $\tilde
F_k\in\End(M^k,E^k)$ such that
$F\circ\pi_k=\pi_k\circ\tilde F_k$, and we have
$\tilde F_k({\bf e}_k)={\bf e}_k$. Furthermore, $d(\tilde F_{\ell(\ca M)})_{{\bf
e}_{\ell(\ca M)}}$ is diagonalizable, with eigenvalues
$\tilde\lambda_1$,~$1$,~$\lambda_2/\lambda_1,\ldots,\lambda_\rho/\lambda_1$ of multiplicity
$1$,~$\mu_1-1$,~$\mu_2,\ldots,\mu_\rho$ respectively, where $\tilde\lambda_1=\lambda_1$
if~$\mu_1>\mu_2$, and $\tilde\lambda_1=\lambda_1^2/\lambda_2$ if $\mu_1=\mu_2$.
More precisely, writing $\chi_{\ell(\ca M)}\circ\tilde F_{\ell(\ca
M)}\circ\chi^{-1}_{\ell(\ca M)}=(\tilde f_1,\ldots,\tilde f_n)$, and denoting by~$a_{11}^j$
the coefficient of $(z_1)^2$ in the power series expansion of~$f_j$, if $\mu_1>\mu_2$ we have
$$
\tilde f_j(w)=\cases{w_1\bigl(\lambda_1-a_{11}^{\mu_1}w_1+2w_2+O(\|w\|^2)\bigr)& if $j=1$,\cr
w_j\bigl(1-{1\over\lambda_1}w_j+{1\over\lambda_1}w_{j+1}+O(\|w\|^2)\bigr)& if $2\le
	j\le\mu_1-1$,\cr
w_{\mu_1}\bigl(1+{a_{11}^{\mu_1}\over\lambda_1}w_1-{1\over\lambda_1}w_{\mu_1}+O(\|w\|^2)\bigr)&
if $j=\mu_1$,\cr
{\lambda_l\over\lambda_1}w_j-{\lambda_l\over\lambda_1^2}w_{j-\nu_l+1}w_j+{1\over\lambda_1}
w_{j-\nu_l+1}w_{j+1}+O(\|w\|^3)& if $j\in\ca P'_{\mu_1,l}\setminus\{\nu_l+\mu_l\}$,
$2\le l\le\rho$,\cr
{\lambda_l\over\lambda_1}w_j-{\lambda_l\over\lambda_1^2}w_{\mu_l+1}w_j+O(\|w\|^3)& if
$j=\nu_l+\mu_l$, $\mu_l<\mu_1-1$,\cr
{\lambda_l\over\lambda_1}w_j+{a_{11}^j\over\lambda_1}w_1w_{\mu_1}
-{\lambda_l\over\lambda_1^2}w_{\mu_1}w_j+O(\|w\|^3)& if $j=\nu_l+\mu_l$,
$\mu_l=\mu_1-1$,\cr}
$$
whereas if $\mu_1=\mu_2$ we have
$$
\tilde f_j(w)=\cases{w_1\bigl({\lambda_1^2\over\lambda_2}-{\lambda_1^2\over\lambda_2^2}
a_{11}^{\nu_2+\mu_2}w_1+{2\lambda_1\over\lambda_2}w_2+O(\|w\|^2)\bigr)& if $j=1$,\cr
w_j\bigl(1-{1\over\lambda_1}w_j+{1\over\lambda_1}w_{j+1}+O(\|w\|^2)\bigr)& if $2\le
j\le\mu_1-1$,\cr 
w_{\mu_1}\bigl(1-{1\over\lambda_1}w_{\mu_1}+O(\|w\|^2)\bigr)& if $j=\mu_1$,\cr
{\lambda_l\over\lambda_1}w_j-{\lambda_l\over\lambda_1^2}w_{j-\nu_l+1}w_j+{1\over\lambda_1}
w_{j-\nu_l+1}w_{j+1}+O(\|w\|^3)& if $j\in\ca P'_{\mu_1,l}\setminus\{\nu_l+\mu_l\}$,
$2\le l\le\rho$,\cr
w_{\mu_2+\nu_2}\bigl({\lambda_2\over\lambda_1}+{a_{11}^{\nu_2+\mu_2}\over\lambda_1}
	w_1+O(\|w\|^2)\bigr)&if $j=\nu_2+\mu_2$,\cr
{\lambda_l\over\lambda_1}w_j+O(\|w\|^3)& if $j=\nu_l+\mu_l$, $\mu_l<\mu_1$,\cr
{\lambda_l\over\lambda_1}w_j+{a_{11}^j\over\lambda_1}w_1w_{\nu_2+\mu_2}+O(\|w\|^3)& if
	$j=\nu_l+\mu_l$, $\mu_l=\mu_1$, $3\le l\le\rho$.\cr}
$$

\pf
Proposition~\tuno\ yields the existence of~$\tilde F_1$; since $\tilde F_1|_{E^1}$ is induced
by the differential of~$F$ at the origin, we see that ${\bf e}_1$ is a fixed point of~$\tilde
F_1$, and more generally that $\tilde F_1(Y^k)=Y^k$ for $k=1,\ldots,\mu_1$. 

By Proposition~\ttre, $d(\tilde F_1)_{[v]}$ is invertible for all $[v]\in
Y^{\mu_1}$. In particular, $\tilde F_1$ is non-degenerate along~$X^1$, and so
Proposition~\tuno\ yields~$\tilde F_2$. Since $d\tilde F_1$ is invertible along~$Y^2$, we can
invoke Proposition~\tdue\ to prove that $d\tilde F_2$ is non-degenerate along~$X^2$, and thus
we get~$\tilde F_3$. Furthermore, being $d\tilde F_2$ invertible along~$X^2$, it is
invertible along the proper transform of~$Y^3$ too, because outside of~$E^2\subset X^2$ it
is given by~$d\tilde F_1$. Then we can again invoke Proposition~\tdue\ to prove that~$\tilde
F_3$ is non-degenerate along~$X^3$, and Proposition~2.1 yields~$\tilde F_4$. Repeating this
procedure we clearly get~$\tilde F_k$ for all~$k$.

To show that ${\bf e}_k$ is a fixed point of~$\tilde F_k$ it suffices to notice that
for $k=2,\ldots,\mu_1$ we have 
$$
\tilde F_1([\de/\de w_k])=[\lambda_1(\de/\de w_k)+(\de/\de w_{k-1})],
$$ 
and $[\de/\de w_{k-1}]\in Y^{k-1}$; analogously, if $\mu_2=\mu_1$ then $\tilde F_1([\de/\de
w_{\nu_2+\mu_2})]=[\lambda_2(\de/\de w_{\nu_2+\mu_2})+(\de/\de w_{\nu_2+\mu_2-1})]$ and
$[\de/\de w_{\nu_2+\mu_2-1}]\in Y^{\mu_1}$. 

We are left to prove that $d(\tilde F_{\ell(\ca M)})_{{\bf e}_{\ell(\ca M)}}$ is
diagonalizable. From
$F\circ\pi_{\ell(\ca M)}=\pi_{\ell(\ca M)}\circ\tilde F_{\ell(\ca M)}$ we easily get
$$
F\circ(\chi_0\circ\pi_{\ell(\ca M)}\circ\chi_{\ell(\ca M)}^{-1})=(\chi_0\circ\pi_{\ell(\ca
	M)}\circ\chi_{\ell(\ca M)}^{-1})\circ\tilde	F\;.
\neweq\eqtcin
$$
Since we know that, writing $F=(f_1,\ldots,f_n)$,
$$
f_j(z)=\cases{\lambda_l z_j+z_{j+1}+\sum\limits_{h,k=1}^n a_{hk}^j z_hz_k+O(\|z\|^3)&if
$\nu_l+1\le j<\nu_l+\mu_l$,\cr
	\lambda_l z_j+\sum\limits_{h,k=1}^n a_{hk}^j z_hz_k+O(\|z\|^3)&if $j=\nu_l+\mu_l$,\cr}
$$
for $1\le l\le\rho$, it is not difficult to check, using \equcin\ and~\equsei, that the
$\tilde f_j$'s have the claimed form, and we are done.\qedn


\smallsect 3. Parabolic curves

From now on we shall assume that $\sp(dF_O)=\{1\}$; in particular, the
Diagonalization Theorem~\qtre\ yields a map tangent to the identity. This allows us
to bring into play Hakim's theory, that we shall now briefly summarize.

Set $\Delta=\{\zeta\in\C\mid |\zeta-1|<1\}$. A
{\it holomorphic curve at the origin} is a holomorphic injective
map~$\phe\colon\Delta\to\C^n\setminus\{O\}$ such that $\phe$ extends continuosly
to~$0\in\de\Delta$ with~$\phe(0)=O$. 

Now take $F\in\End(\C^n,O)$. We shall say that a holomorphic curve at the origin~$\phe$, or
its image~$D=\phe(\Delta)$, is {\it $F$-invariant} if
$F\bigl(\phe(\Delta)\bigr)\subseteq\phe(\Delta)$; that it is {\it stable} if it is
$F$-invariant and $(F|_D)^k\to O$ uniformly on compact subsets of~$D$. A {\it parabolic
curve} is, by definition, a stable holomorphic curve at the origin. Finally, we shall say
that $\phe$ is {\it tangent} to~$v\in\P^{n-1}(\C)$ if
$[\phe(\zeta)]\to v$ as $\zeta\to 0$.

Now let $P_2\colon\C^n\to\C^n$ be a $\C^n$-valued quadratic form. A {\it characteristic
direction} of~$P_2$ is a $v\in\C^n\setminus\{O\}$ such that $P_2(v)=\lambda v$. If
$\lambda=0$ then $v$ is {\it degenerate;} otherwise it is a {\it
non-degenerate} characteristic direction. 

Then (the part we shall need of) Hakim's results can be summarized as follows:

\newthm Theorem \cuno: (Hakim [H2, 3]) 
Let $F\in\End(\C^n,O)$ be such that $dF_O=\id$. Let $P_2\colon\C^n\to\C^n$ be the
quadratic part of the homogeneous expansion of~$F$. If $z^o\in\C^n$, set $z^k=F^k(z^o)$,
and denote by $[z^k]$ its image in $\P^{n-1}(\C)$ when $z^k\neq O$. Then:
{\medskip
\item{\rm (i)}if $z^k\to O$ and $[z^k]\to[v]$ then $v$ is a characteristic direction of
$P_2$;
\item{\rm (ii)}if $v$ is a non-degenerate characteristic direction of $P_2$, then $F$ admits
a parabolic curve tangent to~$[v]$;
\item{\rm (iii)}if $v$ is a non-degenerate characteristic direction of~$P_2$ with
$P_2(v)=\lambda v$ and $D\subset\C^n$ is the parabolic curve given by
part~(ii), then for every $z^o\in D$ and $1\le j\le n$ we have
$$
z^k_j=-{v_j\over\lambda k}+o\left({1\over k}\right)\;.
$$
}

Putting together Theorems~2.4 and~\cuno\ we are able to prove the existence of a parabolic
curve for generic non-diagonalizable maps~$F\in\End(\C^n,O)$ such that $\sp(dF_O)=\{1\}$. In
this context, ``generic'' means $a_{11}^{\mu_1}\neq 0$ and $\mu_2<\mu_1$.

\newthm Corollary \cdue:
Let $F\in\End(\C^n,O)$ be such that $dF_O$ is non-diagonalizable and
$\sp(dF_O)=\{1\}$. Assume without loss of generality that $dF_O$ is in Jordan canonical
form, and let $\ca M$ be the $\rho$-partition of~$n$ induced by the block structure
of~$dF_O$. Assume moreover that $\ell(\ca M)=\mu_1$ and that $a_{11}^{\mu_1}\neq 0$,
where we are using the notations introduced in the previous sections. Then $F$ admits a
parabolic curve~$\phe$ tangent to~$e_1$. Furthermore, if $z^o\in\phe(\Delta)$ and
$z^k=F^k(z^o)$, then 
$$
z^k_j=\cases{(-1)^{\mu_1+j-1}{2\mu_1-1\over a_{11}^{\mu_1}}{2\mu_1-2\choose
\mu_1-1}{(\mu_1+j-2)!\over k^{\mu_1+j-1}}+o\left({1\over k^{\mu_1+j-1}}\right),&
	if $1\le j\le\mu_1$,\cr
o\left({1\over k^{\mu_1+j-\nu_l}}\right),& if $1\le j-\nu_l\le\mu_l<\mu_1-1$,\cr
(-1)^{\mu_1+j-\nu_l}{a_{11}^{\mu_l+\nu_l}(2\mu_1-1)(\mu_l+j-\nu_l)\over
	a_{11}^{\mu_1}}{2\mu_1-2\choose\mu_1-1}{(\mu_1+j-\nu_l-2)!\over
	k^{\mu_1+j-\nu_l}}+o\left({1\over k^{\mu_1+j-\nu_l}}\right),& if
	$1\le j-\nu_l\le\mu_l=\mu_1-1$.\cr}
\neweq\eqcuno
$$

\pf
The idea is to apply Theorem~\cuno\ to the lifting~$\tilde F_{\mu_1}$ of~$F$, and then
use~$\pi_{\mu_1}$ to project the result down to~$F$. Not all the characteristic directions of
the quadratic part of~$\tilde F_{\mu_1}$ at~${\bf e}_{\mu_1}$ are allowable, though. Since
we are working in~$M^{\mu_1}$, characteristic directions tangent to~$\pi^{-1}_{\mu_1}(X^0)$
should be excluded, because the $\tilde F_{\mu_1}$-parabolic curve provided by
Theorem~\cuno.(ii) could be contained in the singular divisor, and thus it would be killed
by~$\pi_{\mu_1}$.  Now, \eqddue\ says that
$\pi_{\mu_1}^{-1}(X^0)$ is given by
$\{w_1=0\}\cup\cdots\cup\{w_{\mu_1}=0\}$; therefore we must look for characteristic
directions~$v$ with~$v_1,\ldots,v_{\mu_1}\neq 0$. Characteristic directions not
tangent to the singular divisor $\pi_k^{-1}(X^0)$ will be called {\it allowable.}

The explicit form of~$\tilde F_{\mu_1}$ given in Theorem~\qtre\ shows that an allowable
characteristic direction~$v$ for~$\tilde F_{\mu_1}$ at~${\bf e}_{\mu_1}$ must satisfy 
$$
\cases{-a^{\mu_1}_{11} v_1+2v_2=\lambda,& for $j=1$,\cr
	-v_j+v_{j+1}=\lambda,& for $2\le j\le \mu_1-1$,\cr
	a^{\mu_1}_{11}v_1-v_{\mu_1}=\lambda,& for $j=\mu_1$,\cr
(-v_j+v_{j+1})v_{j-\nu_l+1}=\lambda v_j,& for $j\in\ca P'_{\mu_1,l}\setminus\{\nu_l+\mu_l\}$,
	$2\le l\le\rho$,\cr
-v_{\mu_l+1}v_j=\lambda v_j,& for $j=\nu_l+\mu_l$, $\mu_l<\mu_1-1$,\cr
a_{11}^jv_1v_{\mu_1}-v_{\mu_1}v_j=\lambda v_j,& for $j=\nu_l+\mu_l$, $\mu_l=\mu_1-1$.\cr}
$$
The unique non-degenerate (i.e., with $\lambda\neq 0$) solution of this system is
$$
v_j=\cases{{1\over a^{\mu_1}_{11}}(2\mu_1-1)\lambda,& for $j=1$,\cr
	(\mu_1+j-2)\lambda,& for $2\le j\le \mu_1$,\cr
	0,& for $j=\nu_l+h$, $1\le h\le\mu_l$, $\mu_l<\mu_1-1$,\cr
{a^{\nu_l+\mu_l}_{11}\over a^{\mu_1}_{11}}(\mu_l+h)\lambda,& for $j=\nu_l+h$, $1\le
h\le\mu_l$, $\mu_l=\mu_1-1$.\cr}
$$
This is an allowable solution; therefore Theorem~\cuno.(ii) yields a $\tilde
F_{\mu_1}$-stable holomorphic curve~$\tilde\phe$ at the origin tangent to~$v$. Since $v$ is
not tangent to~$\pi_{\mu_1}^{-1}(X^0)$, which is invariant under~$\tilde F_{\mu_1}$, the
image of the curve is contained in $M^{\mu_1}\setminus\pi_{\mu_1}^{-1}(X^0)$, which is
exactly the subset of~$M^{\mu_1}$ where~$\pi_{\mu_1}$ is a biholomorphism
with~$\C^n\setminus\{O\}$. Therefore the holomorphic curve
$\phe=\pi_{\mu_1}\circ\tilde\phe$ is a parabolic curve at the origin for~$F$ 
in~$\C^n$, and \eqcuno\ follows from Theorem~\cuno.(iii) and \equsei.\qedn

{\it Remark 3.1:}
Let $\chi\in\Aut(\C^n, O)$ be a (germ of) biholomorphism of~$\C^n$ keeping the origin fixed
and such that the differential of $\hat F=\chi^{-1}\circ F\circ\chi$ is still in Jordan
form; then $\hat a_{11}^{\mu_1}=\alpha\, a_{11}^{\mu_1}$ for a suitable $\alpha\neq 0$, and
thus $F$ is generic iff $\hat F$ is. 
\medbreak

{\it Remark 3.2:}
If $\rho=1$ and $a_{11}^{\mu_1}=0$ but $a_{11}^{\mu_1-1}\neq 0$, it turns out that $d(\tilde
F_{\mu_1-1})_{{\bf e}_{\mu_1-1}}$ is already diagonalizable, and an argument similar
to the one used in the previous proof yields a parabolic curve for $F$ in this case too.
On the other hand, if $\rho\ge 2$ and $\mu_2=\mu_1$ then $\tilde F_{\mu_1+1}$ has {\it no}
allowable non-degenerate characteristic directions at~${\bf e}_{\mu_1+1}$.
\medbreak

{\it Remark 3.3:} We are finally able to explain why diagonalizing simply by blowing-up
points does not work. Indeed, it turns out that in that case the lifted map would have no
allowable characteristic directions; all the relevant dynamics would be inside the singular
divisor, and so one would not easily detect the parabolic curve whose existence is proved
in Corollary~\cdue. 
\medbreak

When $n=2$ (and thus $\rho=1$ and $\mu_1=2$), we are also able to study the non-generic case
$a_{11}^2=0$, obtaining interesting results. For instance, we shall see that (for the first
time, as far as I know) a coefficient of the cubic part of~$F$ enters directly into play even
when the quadratic part of~$F$ is not zero. 

So, assume $n=2$ and $a_{11}^2=0$, and write
$$
\eqalign{f_1(z)&=z_1+z_2+a_{11}^1(z_1)^2+2a_{12}^1 z_1z_2+a_{22}^1(z_2)^2+\cdots,\cr
f_2(z)&=z_2+2a_{12}^2 z_1z_2+a_{22}^2(z_2)^2+a_{111}^2(z_1)^3+\cdots.\cr}
$$
We shall describe our results in terms of the following quantities: 
$$
\eps=a_{11}^1+a_{12}^2,\qquad\hbox{and}\qquad\eta=(a_{11}^1-a_{12}^2)^2+2 a_{111}^2;
$$
they are projective invariants of~$F$ under change of coordinates. More precisely, let again
$\chi\in\Aut(\C^n, O)$ be a (germ of) biholomorphism of~$\C^n$ keeping the origin fixed and
such that the differential of $\hat F=\chi^{-1}\circ F\circ\chi$ is still in Jordan form;
then $\hat a_{11}^2=0$, $\hat\eps=\alpha\eps$ and
$\hat\eta=\alpha^2\eta$ for a suitable $\alpha\neq 0$.

Then:

\newthm Corollary \cqua: Let $F\in\End(\C^2,O)$ be such that $dF_O$ is non-diagonalizable
and $\sp(dF_O)=\{1\}$. Assume that $dF_O$ is in Jordan canonical
form, and that $F$ is non-generic, that is $a_{11}^2=0$. Assume moreover that
$(\eps,\eta)\neq(0,0)$, where $\eps$ and $\eta$ are the invariants just defined. Then:
{\medskip
\item{\rm (i)}if $\eta\neq 0$, $\eps^2$, then $F$ admits two distinct parabolic
curves at the origin;
\item{\rm (ii)}if $\eta=\eps^2\neq 0$, or $\eta=0\neq\eps^2$, then $F$ admits one parabolic
curve at the origin.
\medskip
\noindent In both cases, the parabolic curves are tangent to~$e_1$. Furthermore, if $z^o$
belongs to the image of one of the curves and $z^k=F^k(z^o)$, then
$z^k_1=c_1/k+o(1/k)$ and $z^k_2=c_2/k^2+o(1/k^2)$ for suitable $c_1\neq 0$ and~$c_2\in\C$.}

\pf
The point is that one blow-up is enough to diagonalize such a map; in fact, in this case the
local expansion of $\tilde F_1$ nearby~${\bf e}_1$ is given by
$$
\tilde f_j(w)=\cases{w_1+a_{11}^1(w_1)^2+w_1w_2+O(\|w\|^3),& if $j=1$,\cr
	w_2+a_{111}^2(w_1)^2+(2a_{12}^2-a_{11}^1)w_1w_2-(w_2)^2+O(\|w\|^3),& if $j=2$.\cr}
$$
A direction $[v]\in\P^1(\C)$ is allowable iff $v_1\neq 0$; therefore we can assume
$v_1=1$, and finding the allowable characteristic directions boils down to solving a
quadratic equation whose discriminant is~$\eta$. The allowable characteristic directions
then are multiple of
$$
v_\pm=\left(1,{a_{12}^2-a_{11}^1\pm\sqrt{\eta}\over 2}\right),
$$
and $v_\pm$ is degenerate iff $\eps\pm\sqrt{\eta}=0$. Theorem~\cuno\ thus yields the
assertion, exactly as in the previous corollary.\qedn

{\it Remark 3.4:} 
If $\eps=\eta=0$ several things might happen; we can even have more than
two stable holomorphic curves at the origin. See~[A] and~[CD] for examples.
\medbreak

{\it Remark 3.5:}
A $\C^n$-valued quadratic form $P_2$ on~$\C^n$ induces on the projective space a holomorphic
map $\hat P_2\colon\P^{n-1}(\C)\setminus Z\to\P^{n-1}(\C)$, where $Z$ is the image
in~$\P^{n-1}(\C)$ of the cone $P_2^{-1}(O)\setminus\{O\}\subset\C^n$. If $v\in\C^n$ is a
non-degenerate characteristic direction for~$P_2$, then its image $[v]\in\P^{n-1}(\C)$ is a
fixed point of~$\hat P_2$. In particular, we may then consider the linear map
$$
A_{[v]}=d(\hat P_2)_{[v]}-\id\colon T_{[v]}\bigl(\P^{n-1}(\C)\bigr)\to T_{[v]}\bigl(
	\P^{n-1}(\C)\bigr)\;.
$$
It turns out that this is the same matrix introduced by Hakim~[H2, 3]. She
proved that, under the hypotheses of Theorem~\cuno, if $A_{[v]}$ has $d\ge 0$ eigenvalues
with positive real part then the map actually admits a parabolic holomorphic $(d+1)$-manifold
at the origin. In the case $n=2$, $a_{11}^2=0$ and $(\eps,\eta)\neq(0,0)$, we have
$$
A_{[v_\pm]}=\mp\,2{\sqrt{\eta}\over\eps\pm\sqrt{\eta}}.
$$
In particular, $A_{[v]}=-1$ when $\eta=\eps^2\neq 0$ (where, choosing $\eps$ as 
principal determination of~$\sqrt{\eta}$, the non-degenerate characteristic direction
is~$v_+$), $A_{[v]}=0$ when $\eta=0\neq\eps$, and $\Re A_{[v_\pm]}>0$ iff
$$
\Re\left({\eps\over\pm\sqrt{\eta}}\right)< -1,
$$
when $\eta\neq 0$,~$\eps^2$. In particular,  if $|\Re(\eps/\sqrt{\eta})|>1$ then the map
$F$ admits a parabolic basin of attraction for the origin.
\medbreak

{\it Remark 3.6:}
It is not difficult to compute the matrix $A_{[v]}$ for the allowable characteristic
direction described in the proof of Corollary~\cdue; it is not so easy to compute
the sign of the real part of the eigenvalues, though. For $n\le 20$ we checked that the
matrix~$A_{[v]}$ has no eigenvalue with positive real part, and we suspect that this is true
for all~$n$.
\medbreak

{\it Remark 3.7:}
Hakim~[H3] proved that when $\tilde F\in\End(\C^n,O)$ is a global automorphism of~$\C^n$ with
$d\tilde F_O=\id$, and $v$ is a non-degenerate characteristic direction, then the
set~$\Omega_v$ of orbits $z^k\to O$ such that $[z^k]\to[v]$ is an $\tilde F$-stable
biholomorphic image of~$\C^{d+1}$, where $d\ge 0$ is the number of eigenvalues of~$A_{[v]}$
with positive real part (assuming, for simplicity, that $A_{[v]}$ has no purely imaginary
eigenvalues). This is still true in our situation. Indeed, if our map $F$ is a global
automorphism of~$\C^n$, then its lifting~$\tilde F$ is a global automorphism
of~$M^{\mu_1}\setminus\pi^{-1}_{\mu_1}(X^0)$, which is biholomorphic to~$\C^n\setminus\{O\}$.
Furthermore, if $v$ is an allowable characteristic direction, then~$\Omega_v$ cannot
intersect the singular divisor, because the latter is~$\tilde F$-invariant whereas $v$ is
not tangent to it. This means that we can apply Hakim's result to~$\tilde F$, and
projecting down via~$\pi_{\mu_1}$ we get an $F$-stable $(d+1)$-manifold biholomorphic
to~$\C^{d+1}$. In particular, then, Remark~3.5 yields yet another instance of the
Fatou-Bieberbach phenomenon in~$\C^2$.

\smallsect 4. Regular orbits

In the previous section we have shown that allowable (i.e., not tangent to the
singular divisor) characteristic directions of the lifting of a map~$F$ give
rise to parabolic curves. A priori, other characteristic directions might
also give rise to parabolic curves, or possibly to $F$-orbits converging to the origin. The
aim of this section is to show that this cannot happen, at least in the case $\rho=1$, when
$dF_O$ is the Jordan $n\times n$ block~$J_n$ associated to the eigenvalue~1.

To state more precisely our result, we need some definitions. Let
$\{z^k\}\subset\C^n\setminus\{O\}$ be a sequence converging to the origin. We shall say that
$\{z^k\}$ is {\it $0$-regular} if $\{[z^k]\}$ converges to some $[v]\in\P^{n-1}(\C)$; this
is equivalent to saying that $\pi_1^{-1}(z^k)$ converges to some $[v]\in E^1$. We shall say
that $\{z^k\}$ is {\it $1$-regular} if either~$[v]\neq{\bf e}_1$ (and we shall specify this
case saying that it is $1$-regular {\it of first kind\/}) or $[v]={\bf e}_1$
and~$\{\chi_1\circ\pi_1^{-1}(z^k)\}$ is 0-regular (and then $\{z^k\}$ is $1$-regular {\it of
second kind\/}). Now we proceed by induction. Let~$\{z^k\}$ be $(r-1)$-regular. If it is
$(r-1)$-regular of first kind, we shall also say that it is {\it $r$-regular (of first
kind).} If it is $(r-1)$-regular of second kind, then $\pi_r^{-1}(z^k)$ converges to some
$[v]\in E^r$. We shall say that~$\{z^k\}$ is {\it $r$-regular} if either $[v]\neq{\bf e}_r$
(and we shall again say $r$-regular {\it of first kind\/}) or $[v]={\bf e}_r$ and
$\{\chi_r\circ\pi^{-1}_r(z^k)\}$ is 0-regular (and then $\{z^k\}$ is $r$-regular {\it of
second kind\/}). We stress that we impose no conditions if $[v]\neq{\bf e}_r$; so for most
sequences $r$-regularity is equivalent to $0$-regularity.

Despite its apparent complexity, the condition of $r$-regularity is fairly natural; it is
just a way to say that the different components of the sequence go to zero at comparable
rates. For instance, if for $j=1,\ldots,n$ there are
$a_j\in\C^*$ and $\delta_j>0$ such that
$$
z_j^k={a_j\over k^{\delta_j}}+o\left({1\over k^{\delta_j}}\right)\;,
$$
then $\{z^k\}$ is $r$-regular for every~$r$; and it is easy
to provide examples of much more general $r$-regular sequences.

Now let $F\in\End(\C^n,O)$ be such that $dF_O=J_n$. Assume that $F$ is generic, that is
$a_{11}^n\neq0$, and let~$\tilde F$ be its lifting. We shall say that an $F$-orbit
is {\it regular} if it converges to the origin and it is $n$-regular. A quick look to
\equcin\ and \equsei\ shows that orbits obtained pushing down $0$-regular orbits of $\tilde
F$ tangent to allowable characteristic directions are regular; such orbits are called {\it
standard,} and are the ones described in Corollary~\cdue.
Using this terminology, our aim is to prove that every regular orbit is
standard. To do so, we need a lemma:

\newthm Lemma \suno:
Let $\{w^k\}\subset\C^*$ be a sequence converging to~$0$. Assume there is another sequence
$\{u^k\}\subset\C$ such that $u^k/w^k\to c\in\C$ and
$$
w^{k+1}=w^k(1+u^k)+o\bigl((w^k)^2\bigr)\;.
$$
Then $1/(kw^k)\to -c$. In particular, if $c\neq 0$ we have
$$
w^k=-{1\over ck}+o\left({1\over k}\right)\;.
$$

\pf
Set $\eps^k=w^{k+1}-w^k-u^kw^k$, so that $\eps^k/(w^k)^2\to 0$. We then have
$$
{1\over w^{h+1}}={1\over w^h}-{u^h\over w^h}+{(u^h)^2/w^h+(u^h-1)\eps^h/(w^h)^2\over
	1+u^h+\eps^h/w^h}\;.
$$
Summing this equality for $h=0,\ldots, k-1$ and dividing by~$k$ we find
$$
{1\over kw^k}={1\over kw^0}-{1\over k}\sum_{h=0}^{k-1}{u^h\over w^h}+
	{1\over k}\sum_{h=0}^{k-1}{(u^h)^2/w^h+(u^h-1)\eps^h/(w^h)^2\over
	1+u^h+\eps^h/w^h}\;,
$$
and the assertion follows from the convergence of the averages of a converging sequence.\qedn

Then:

\newthm Theorem \sdue:
Let $F\in\End(\C^n,O)$ be such that $dF_O=J_n$. Assume that $F$ is generic. Then every
regular orbit of $F$ is standard.

\pf
Up to a linear change of coordinates we can assume $a_{11}^n=1$. Let
$\{z^k=F^k(z^o)\}$ be a regular orbit; we first of all want to prove, by induction, that
$\pi_r^{-1}(z^k)\to{\bf e}_r$ for
$r=1,\ldots,n$.

First of all, $0$-regularity yields $[z^k]\to[v]\in\P^{n-1}(\C)$. But then $v$ must be an
eigenvector of~$dF_O$; therefore~$[v]={\bf e}_1$, and thus $\pi_1^{-1}(z^k)\to{\bf e}_1$.
Exactly the same argument shows that $\pi_2^{-1}(z^k)\to{\bf e}_2$.

Now assume that $\pi_r^{-1}(z^k)\to{\bf e}_r$ for some $2\le r\le n-1$, and put
$w^k=\chi_r\circ\pi_r^{-1}(z^k)$. The $0$-regularity of~$\{w^k\}$ implies that
$[w^k]\to[v]\in\P^{n-1}(\C)$; again, $v$ must be (canonically identified to) an
eigenvector of~$d(\tilde F_r)_{{\bf e}_r}$. Now, a computation using \equcin\ and \equsei\
shows that for $1<r<n$ we have
$$
w_j^1=\cases{w_1\bigl(1-a_{11}^r w_1+2w_2-w_{r+1}+O(\|w|^2)\bigr)& if $j=1$,\cr
w_j\bigl(1-w_j+w_{j+1}+O(\|w|^2)\bigr)& if $1<j<r$,\cr
w_r\bigl(1+a_{11}^r w_1-w_r+w_{r+1}+O(\|w|^2)\bigr)& if $j=r$,\cr
a_{11}^jw_1+w_j+w_{j+1}+2a_{12}^jw_1w_2-(a_{11}^r
w_1+w_{r+1})(a_{11}^jw_1+w_j+w_{j+1})+O(\|w|^3)&if $r<j<n$,\cr
a_{11}^n w_1+w_n+2a_{12}^n w_1w_2-(a_{11}^r w_1+w_{r+1})(a_{11}^n w_1+w_n)+O(\|w|^3)&
	if $j=n$.\cr}
\neweq\eqqzero
$$
In particolar, $d(\tilde F_r)_{{\bf e}_r}$ is represented by the matrix
$$
\left|\vcenter{\offinterlineskip
\halign{&\hfil$#$\hfil\cr
\multispan5\strut&\omit\hskip1.5mm\vrule width 1pt\hskip 1.5mm&\cr
\multispan5\strut\hskip1.5mm\hfil $I_r$\hfil&\omit\hskip1.5mm\vrule width 1pt
	\hskip 1.5mm&\multispan3\hfil$O$\hskip1.5mm\hfil\cr
\multispan5\strut&\omit\hskip1.5mm\vrule width1pt\hskip 1.5mm&\cr
\noalign{\hrule height1pt}
\multispan5&\omit\hskip1.5mm\vrule height 2pt
	width 1pt\hskip 1.5mm&\multispan3\cr
\strut a_{11}^{r+1}&\omit\hskip 1mm\vrule\hskip1mm&\multispan3&\omit\hskip1.5mm\vrule 
	width 1pt\hskip 1.5mm&\multispan3\cr
\strut\vdots&\omit\hskip 1mm\vrule\hskip1mm&\quad&\hfil$\smash{\raise
	3pt\vbox{\smash{$O$}}}$\hfil&\quad&\omit\hskip1.5mm\vrule width1pt
	\hskip 1.5mm&\quad&\hfil$\smash{\raise 3pt\vbox{\smash{$J_{n-r}$}}}$\hfil&\quad\cr
\strut a_{11}^n&\omit\hskip 1mm\vrule\hskip1mm&\multispan3&\omit\hskip1.5mm\vrule width 1pt
	\hskip 1.5mm&\multispan3\cr}}
\right|\;.
$$
Therefore $v=(0,v_2,\ldots,v_{r+1},0,\ldots,0)$; to prove that $\pi^{-1}_{r+1}(z^k)\to{\bf
e}_{r+1}$ it suffices to show that $v_{r+1}\neq 0$.

Assume, by contradiction, $v_{r+1}=0$, and let $j_0=\max\{2\le j\le r\mid v_j\neq0\}$. We
know that $w^k_j/w^k_{j_0}\to v_j/v_{j_0}$ for all $j$; in particular, $w^k_j=O(w^k_{j_0})$
if $v_j\neq 0$, and $w^k_j=o(w^k_{j_0})$ if $v_j=0$. Then~\eqqzero\ yields
$$
w^{k+1}_{j_0}=w^k_{j_0}(1-w^k_{j_0})+o\bigl((w^k_{j_0})^2\bigr)\;;
$$
hence using Lemma~\suno\ we find
$$
w^k_{j_0}={1\over k}+o\left({1\over k}\right)\;,
$$
and so
$$
w^k_j={v_j/v_{j_0}\over k}+o\left({1\over k}\right)
\neweq\eqquno
$$
for all $j=1,\ldots,n$. 

We now claim that $v_j/v_{j_0}=j_0-j+1$ for all $j=2,\ldots,j_0$. We argue by induction on
$j_0-j$. Take~$j<j_0$ and assume that $v_{j+1}/v_{j_0}=j_0-j$. Noticing that $w_j^k\neq
0$ for all~$k$ and~$1\le j\le r$ (because $\pi^{-1}_r(z^k)$ does not belong to the
singular divisor), we can write
$$
{w_j^{k+1}\over w_j^k}=1-w_j^k+w_{j+1}^k+O\bigl((w_{j+1}^k)^2\bigr)\;.
$$
If $v_j=0$ we would get
$$
{w_j^{k+1}\over w_j^k}=1+{j_0-j\over k}+o\left({1\over k}\right)\;,
$$
which is impossible because the infinite product $\prod_k(w_j^{k+1}/w_j^k)$ is converging
to~zero. Therefore $v_j\neq 0$; but then applying Lemma~\suno\ to
$$
w^{k+1}_j=w_j^k(1-w^k_j+w^k_{j+1})+o\bigl((w_j^k)^2\bigr)
$$
and recalling \eqquno\ we get $v_j/v_{j_0}=j_0-j+1$, as claimed.

In particular we then have $v_2/v_{j_0}=j_0-1$, and so
$$
{w_1^{k+1}\over w_1^k}=1+{2(j_0-1)\over k}+o\left({1\over k}\right)\;,
$$
which is impossible. The contradiction arises because we assumed $v_{r+1}=0$; therefore we
must have $v_{r+1}\neq 0$, as claimed.

Summing up, we have in particular proved that $\pi_n^{-1}(z^k)\to{\bf e}_n$; set
$w^k=\chi_n\circ\pi^{-1}_n(z^k)$. Notice that, by construction, $w^k_j\neq 0$ for all~$k$
and~$j$. By $0$-regularity, $[w^k]\to[v]\in\P^{n-1}(\C)$; Theorem~\cuno\ then says that $v$
must be a characteristic direction of~$\tilde F_n$ at~${\bf e}_n$, that is a solution of
$$
\cases{-v_1^2+2v_1v_2=\lambda v_1\;,\cr
\noalign{\smallskip}
	-v_j^2+v_jv_{j+1}=\lambda v_j& for $j=2,\ldots,n-1$,\cr
	v_1v_n-v_n^2=\lambda v_n.\cr}
$$
To end the proof we must show that $v$ is allowable, that is that $v_j\ne 0$ for
$j=1,\ldots,n$.

Assume, by contradiction, that there is a $j_0$ such that $v_{j_0}\neq 0$ but $v_{j_0+1}=0$
(where here by $v_{n+1}$ we mean~$v_1$). Then it is easy to prove that
$v_j/v_{j_0}\in\N$ for all $j=1,\ldots,n$; in particular, $v_j/v_{j_0}$ is always
non-negative. Now we have
$$
w^{k+1}_{j_0}=w^k_{j_0}(1-w^k_{j_0})+o\bigl((w^k_{j_0})^2\bigr)\;;
$$
therefore Lemma~\suno\ yields $w_{j_0}^k=1/k+o(1/k)$. Recalling \eqquno\ we then get
$w_j^k=c_j/k+o(1/k)$ with $c_j\ge 0$ for all~$j=1,\ldots,n$. But then arguing exactly as in
the first part of the proof we show that $v_{j_0-1},\ldots,v_1\neq 0$; and then we
get~$v_n\neq 0$, and going up we finally arrive to prove $v_{j_0+1}\neq 0$,
contradiction.\qedn

\setref{CG}
\beginsection References

\pre A M. Abate: Holomorphic dynamical systems with a Jordan fixed point! Preprint!
	1998

\art CD D. Coman, M. Dabija: On the dynamics of some diffeomorphisms of $\C^2$ near
parabolic fixed points! Houston J. Math.! 24 1998 85-96

\art F P. Fatou: Substitutions analytiques et \'equations fonctionnelles \`a deux
variables! Ann. Sc. Ec. Norm. Sup.! {} 1924 67-142

\art H1 M. Hakim: Semi-attractive transformations of $\C^p$! Publ. Math! 38 1994
479-499

\art H2 M. Hakim: Analytic transformations of $(\C^p,0)$ tangent to the identity!
Duke Math. J.! 92 1998 403-428

\pre H3 M. Hakim: Stable pieces of manifolds in transformations tangent to the identity!
Preprint! 1997

\art N Y. Nishimura: Automorphismes analytiques admettant des sousvari\'et\'e de points fixes
attractives dans la direction transversale! J. Math. Kyoto Univ.! 23 1983 289-299

\art R1 L. Reich: Das Typenproblem bei formal-biholomorphen Abbildungen mit anziehendem
Fixpunkt! Math. Ann.! 179 1969 227-250

\art R2 L. Reich: Normalformen biholomorpher Abbildungen mit anziehendem
Fixpunkt! Math. Ann.! 180 1969 233-255

\pre Ri M. Rivi: Stable manifolds for semi-attractive holomorphic germs! Preprint!
1999

\art S S. Sternberg: Local contractions and a theorem of Poincar\'e! Amer. J. Math.! 79
1957 809-824

\art U1 T. Ueda: Local structure of analytic transformations of two complex variables,
I! J. Math. Kyoto Univ.! 26 1986 233-261

\art U2 T. Ueda: Local structure of analytic transformations of two complex variables, 
I\negthinspace I! J. Math. Kyoto Univ.! 31 1991 695-711

\art W B.J. Weickert: Attracting basins for automorphisms of $\C^2$! Invent. Math.!
132 1998 581-605

\art Wu H. Wu: Complex stable manifolds of holomorphic diffeomorphisms! Indiana Univ.
Math. J.! 42 1993 1349-1358

\bye

\smallsect 1. Introduction

In the study of discrete dynamical systems, the description of the asymptotic behavior in a
neighbourhood of a fixed point is of course of paramount importance. The case of a
one-dimensional holomorphic dynamical system generated by a function $f$ with $f(0)=0$ is
fairly well understood (see, e.g., [2,~10]): the behavior depends strongly on the value
of~$\lambda=f'(0)$. 

If $|\lambda|<1$, the origin is attracting, that is there is an open set
$\Omega$ containing~$0$ (the {\it basin of attraction} of the origin) such that every point
of~$\Omega$ converges to~$0$ under iteration of~$f$; furthermore, if $\lambda\ne 0$ then $f$
is linearizable in a neighbourhood of~$0$. If $|\lambda|>1$, the origin is repelling, that is
it is attracting for $f^{-1}$ (which is defined in a neighbourhood of~$0$).

If $|\lambda|=1$ there are three possibilities. If $\lambda$ is a root of unity, then there
is again a basin of attraction~$\Omega$ for the origin (the {\it Leau-Fatou flower\/}), but
this time~$0\in\de\Omega$; furthermore, $f$ is semiconjugate to~$z\mapsto z+1$ on~$\Omega$
(and actually conjugated to it in an open subset of~$\Omega$ containing~$0$ in its boundary).
If $\lambda$ is not a root of unity but satisfies a Diophantine condition (the
Bryuno condition), then Yoccoz (see~[16]) has proved that there is an open set $\Omega$
containing the origin where $f$ is conjugated to an irrational rotation. Finally, if
$\lambda$ neither is a root of unity or satisfies the Bryuno condition then $f$ is not
linearizable in any neighbourhood of~$0$,
and recent results by P\'erez-Marco (see~[11] and references therein) clarify the dynamics in
this case too; in particular, there are no $f$-orbits converging to~$0$ (except the
trivial ones contained in the inverse orbit of the origin).

The situation is not so clear for multidimensional holomorphic dynamical systems. Let $F$ be
a germ of holomorphic self-map of $\C^n$ (with $n\ge 2$) defined in a neighbourhood of the
origin with $F(O)=O$ If no eigenvalue of the differential~$dF_O$ has absolute value one
(the hyperbolic case), then the situation is described by the stable manifold
theorem (see~[15] for a presentation in the complex setting): if $d$ is the number of
eigenvalues (counted with multiplicities) of~$dF_O$ with absolute value less than~1, then
there is are an $F$-invariant $d$-dimensional complex manifold~$W^s$ and an $F$-invariant
$(n-d)$-dimensional complex manifold~$W^u$ such that: $W^s$ and~$W^u$ are transversal at the
origin; $W^s$ is {\it stable} (that is~$F^k(z)\to O$ for all~$z\in W^s$), and $W^u$ is {\it
unstable} (that is stable with respect to~$F^{-1}$). 

In the past few years, a number of results have been obtained in the non-hyperbolic case.
The first one (due originally to Fatou~[4], precised by Ueda~[12, 13] and put into final
form by Hakim~[7]) is about maps~$F$ such that~$dF_O$ has $1$ as eigenvalue of
multiplicity one, and the others eigenvalues have absolute value less than~1. In this case
either $F$ admits a holomorphic curve of fixed points passing through the origin or there
exists a basin of attraction to the origin.

Another situation that has been studied is when $dF_O=\id$. In this case Hakim~[8,~9] (see
also Weickert~[14]) has proved that for $F$ generic there exists an $F$-invariant stable
(i.e., attracted to the origin) holomorphic curve with the origin {\it in its boundary;}
furthermore, we have estimates on the rate of approach of stable orbits to the origin (see
Section~5 for a precise statement of Hakim's results). Notice that, in general, it is not
possible to extend such a stable curve holomorphically through the origin. 

This paper is devoted to the case $dF_O=J$, the Jordan block of order~$n$
associated to the eigenvalue~1 --- that is, to the non-semisimple non-hyperbolic case,
which, as far as I know, up to now has been studied only in very particular examples (see~[1,
3]). Most of our work will be local in nature;
however, as we shall make precise for instance in Remarks~2.1 and~5.6, several results have
global consequences for biholomorphisms of~$\C^n$.

Roughly speaking, we can think of a non-semisimple map (that is, a map with a
non-semisimple differential) as a sort of singularity in the space of all (holomorphic)
maps fixing the origin; indeed, a generic map has semisimple differential. In algebraic
geometry, to resolve a singularity one uses blow-ups; our idea is that the same approach
works for ``resolving'' the non-semisimplicity. In other words, we claim that using a
sequence of blow-ups it is possible to lift the original map to a map which is defined on a
much more complicated space but it has semisimple differential --- and then we are able
to deduce dynamical informations on the original map via the dynamics of the lifted map.
The details are, of course, involved; in particular, it turns out that to extract
dynamically relevant consequences from the construction it is necessary to
blow-up submanifolds, and not only points. Though in Section~3 we shall lay out all the
needed tools and theorems, for the sake of simplicity (and for presenting dynamically
interesting results) in the rest of the paper we shall limit ourselves to the discussion of
the case~$dF_O=J$, describing in detail the sequence of blow-ups needed to ``desingularize''
this type of maps (Theorem~4.3).

The desingularization technique will allow us to apply Hakim's theory to non-semisimple
maps. For instance, we shall be able to prove (Corollaries~5.2, 5.3 and~5.4) that 
generic maps with~$F(O)=O$ and~$dF_O=J$ admits an $F$-invariant stable holomorphic curve with
the origin in its boundary, exactly as in the semisimple case. We shall also be able to
prove that, roughly speaking, this holomorphic curve contains all the stable dynamics of~$F$
in a neighbourhood of the origin (Theorem~6.2).

This paper is organized as follows. In Section~2 we describe a normal form for maps
with~$dF_O=J$; we shall be able to diagonalize the quadratic part of the homogeneous
expansion of~$F$ (see Theorem~2.1 for a precise statement), and to recover a
dynamically meaningful projective invariant (Proposition~2.4). In Section~3 we show how to
lift a map to a blow-up along an invariant submanifold, and we shall describe sufficient
conditions ensuring the possibility to iterate the procedure. In Section~4 we prove the
Desingularization Theorem~4.3 via iterated blow-ups, and we describe the local structure of
the desingularized map. In Section~5 we apply Hakim's theory to the desingularization to get
the existence of stable curves for the original map, and in particular we recover yet
another instance of the Fatou-Bieberbach phenomen (see Remark~5.6). Finally, in Section~6 we
show that, under a weak regularity assumption, the $F$-orbits converging to the origin must
be exactly the ones given by the desingularization.

\smallsect 2. Normal forms

We denote by $\End(\C^n,O)$ the set of germs of holomorphic self-maps of $\C^n$
sending the origin~$O$ to itself; more generally, if $Y$ is a closed set of a complex
manifold~$M$, we shall denote by $\End(M,X)$ the set of germs at $X$ of holomorphic
self-maps of~$M$ sending $X$ into itself. $\Aut(\C^n,O)\subset\End(\C^n,O)$ will denote the
set of germs of biholomorphisms fixing the origin. 

A germ $F\in\End(\C^n,O)$ has a {\it homogeneous expansion} of the form
$$
F(z)=\sum_{j=1}^\infty P_j(z)\;,	
\eqno (2.1)
$$
where $z=(z_1,\ldots,z_n)\in\C^n$, and the $P_j$'s are $n$-uples of homogeneous polynomials
of degree~$j$ in~$z_1,\ldots,z_n$. If~$P_1(z)=Jz$, where
$$
J=J_n=\left|\matrix{1&1&&0\cr &\ddots&\ddots&\cr &&\ddots&1\cr 0&&&1\cr}\right|\in GL_n(\C)
$$
is the standard Jordan block of dimension~$n$ with eigenvalue one, we shall say that $F$ has
a {\it Jordan fixed point} at the origin. 

The aim of this section is the simplification of~$P_2$, the quadratic part of~$F$. The
$\C^n$-valued quadratic form $P_2$ is a $n$-uple $P_2=(\phe_1,\ldots,\phe_n)$ of quadratic
forms; we shall consistently denote by $\Phi_j\colon\C^n\times\C^n\to\C$ the symmetric
bilinear form associated to~$\phe_j$ --- that is $\phe_j(z)=\Phi_j(z,z)$ --- and set
$\Phi_F=(\Phi_1,\ldots,\Phi_n)$. 

We shall say that $F$ is in {\it normal form} if $\phe_1,\ldots,\phe_{n-1}$ are diagonal and
$\phe_n(z)=\eps^{j_0}_n(z_{j_0})^2$ for suitable~$1\le j_0\le n$ and $\eps^{j_0}_n\in\C$ 
($\phe_n\equiv 0$ is allowed). In other words, $F=(f_1,\ldots,f_n)$ is in normal form iff it
can be written as
$$
\cases{\displaystyle f_1(z)=z_1+z_2+\sum_{k=1}^n \eps^k_1 (z_k)^2+O_3\;,\cr
	\displaystyle\hfill\vdots\qquad\hfill\cr
	\displaystyle f_{n-1}(z)=z_{n-1}+z_n+\sum_{k=1}^n \eps^k_{n-1} (z_k)^2+O_3\;,\cr
	\displaystyle f_n(z)=z_n+\eps^{j_0}_n (z_{j_0})^2+O_3\;,\cr}
\eqno (2.2)
$$
where, here and elsewhere, $O_k$ denotes a germ of holomorphic self-map of $\C^n$ whose
homogeneous expansion has no addend of degree less than~$k$. 

The main theorem of this section is that for any $F$ there is a $\chi\in\Aut(\C^n,O)$ such
that $\chi^{-1}\circ F\circ\chi$ is in normal form (the complete statement is even
slightly better; see Theorem~2.1). Before starting with the proof, we record here a couple
of easy computations.

Let $F=\sum_j P_j$, $\tilde F=\sum_j\tilde P_j$ be the homogeneous expansion of two elements
of~$\End(\C^n,O)$, and let 
$$
F\circ\tilde F=\sum_{j=1}^\infty R_j
$$ 
be the homogeneous expansion of the composition. Then it is easy to check that
$$
\eqalign{R_1&=P_1\circ\tilde P_1\;,\cr
	R_2&=P_1\circ\tilde P_2+P_2\circ\tilde P_1\;,\cr
	R_3&=P_1\circ\tilde P_3+2\Phi_F\circ(\tilde P_1\oplus\tilde P_2)+P_3\circ\tilde P_1\;,\cr}
\eqno (2.3)
$$
where $(\tilde P_1\oplus\tilde P_2)(z)=\bigl(\tilde P_1(z),\tilde P_2(z)\bigr)$.

Now let $\chi=\sum_{j=1}^\infty Q_j$ be a (germ of) automorphism of $\C^n$ fixing the
origin. In particular, $Q_1\in GL_n(\C)$,
and Eq.~(2.3) allows us to compute the first few terms
in the homogeneous expansion $\sum_{j=1}^\infty \tilde Q_j$ of $\chi^{-1}$:
$$
\eqalign{\tilde Q_1&=Q_1^{-1}\;,\cr
	\tilde Q_2&=-Q_1^{-1}\circ Q_2\circ Q_1^{-1}\;,\cr
	\tilde Q_3&=-Q_1^{-1}\circ Q_3\circ Q^{-1}_1+2Q^{-1}_1\circ\Phi_\chi\circ\bigl(Q^{-1}_1
	\oplus(Q^{-1}_1\circ Q_2\circ Q^{-1}_1)\bigr)\;,\cr
	\Phi_{\chi^{-1}}&=-Q^{-1}_1\circ\Phi_\chi\circ(Q^{-1}_1\times Q^{-1}_1)\;,\cr}
\eqno (2.4)
$$
where $Q^{-1}_1\times Q^{-1}_1(z,w)=\bigl(Q^{-1}_1(z),Q^{-1}_1(w)\bigr)$.

Putting Eqs.~(2.3) and~(2.4) together, we see that the homogeneous expansion of $\tilde
F=\chi^{-1}\circ F\circ\chi$ starts with
$$
\eqalign{\tilde P_1&=Q^{-1}_1\circ P_1\circ Q_1\;,\cr
\tilde P_2&=Q^{-1}_1\circ P_1\circ Q_2+Q^{-1}_1\circ P_2\circ Q_1-Q^{-1}_1\circ Q_2\circ
	\tilde P_1\;,\cr
\tilde P_3&=Q^{-1}_1\circ\Bigl[P_1\circ Q_3-Q_3\circ\tilde P_1+P_3\circ Q_1+2\Phi_F\circ
	(Q_1\oplus Q_2)-2\Phi_\chi\circ(\tilde P_1\oplus\tilde P_2)\Bigr]\;.\cr}
\eqno (1.5)
$$
We shall be interested only in change of coordinates preserving the linear part of~$F$, that
is such that $\tilde P_1=P_1=J$. Then $Q_1$ must be of the form
$$
Q_1=J_{(\alpha)}=\left|\matrix{\alpha_0&\alpha_1&\cdots&\alpha_{n-1}\cr
	&\ddots&\ddots&\vdots\cr &&\ddots&\alpha_1\cr 0&&&\alpha_0\cr}\right|\in GL_n(\C),
$$
where $\alpha=(\alpha_0,\ldots,\alpha_{n-1})\in\C^*\times\C^{n-1}$, and Eq~(2.5) becomes
$$
\eqalign{\tilde P_2&=Q_1^{-1}\circ\Bigl[P_1\circ Q_2-Q_2\circ P_1+P_2\circ Q_1\Bigr]\;,\cr
	\tilde P_3&=Q^{-1}_1\circ\Bigl[P_1\circ Q_3-Q_3\circ P_1+P_3\circ Q_1+2\Phi_F\circ
	(Q_1\oplus Q_2)-2\Phi_\chi\circ(P_1\oplus\tilde P_2)\Bigr]\;.\cr}
\eqno (2.6)
$$

The main theorem of this section is the following:

\newthm Theorem \duno: Let $F\in\End(\C^n,O)$ be with a Jordan fixed point at the
origin. Then there is $\chi\in\Aut(\C^n,O)$ such that $\tilde F=\chi^{-1}\circ F\circ\chi$ is
in normal form. Furthermore we can also assume that $\eps^k_h=0$ if $h=1,\ldots,n$ and
$k>\lceil n/2\rceil$, and that $j_0\le\lceil n/2\rceil$, where we are using the
notation introduced in Eq.~(2.2).

The proof depends on two lemmas.

\newthm Lemma \ddue: Let $\phe(z)=z^T Az$ be a quadratic form on~$\C^n$, where $A=(a_{hk})\in
M_{n,n}(\C)$ is symmetric. Then there exists a quadratic form $\psi$ on~$\C^n$ such that
setting 
$$
\tilde\phe=\phe-\psi\circ J'-\Psi\circ(I\oplus J')-\Psi\circ(J'\oplus I)\;,
$$
where $\Psi$ is the symmetric bilinear form associated to~$\psi$, $I\in GL_n(\C)$ is the
identity matrix, and $J'=J-I$, we have
$$
\tilde\phe(z)=a_{11}(z_1)^2+\sum_{h=2}^{\lceil n/2 \rceil}\eps_h(z_h)^2
$$
for suitable $\eps_2,\ldots,\eps_{\lceil n/2 \rceil}\in\C$.

\pf 
Let us write $\psi(z)=z^T Bz$, and $\tilde\phe(z)=z^T\tilde A z$, where
$B=(b_{hk})$,
$\tilde A=(\tilde a_{hk})\in M_{n,n}(\C)$ are symmetric. Since
$J'(z)=(z_2,\ldots,z_n,0)$, we have
$$
\cases{\tilde a_{11}=a_{11}\;,\cr
	\tilde a_{1k}=a_{1k}-b_{1,k-1}& for $k=2,\ldots,n$,\cr
	\tilde a_{h1}=a_{h1}-b_{h-1,1}& for $h=2,\ldots,n$,\cr
	\tilde a_{hk}=a_{hk}-b_{h-1,k-1}-b_{h,k-1}-b_{h-1,k}& for $h$,~$k=2,\ldots,n$.\cr}
\eqno (2.7)
$$
So our aim is to find $b_{hk}$ such that $\tilde a_{hk}=0$ if $h\neq k$ or $h=k>\lceil n/2
\rceil$. We shall prove that this is possible working by induction on $h$. More precisely,
we shall prove by induction the following assertion $\ca A_h$: 
\bigskip
{\narrower \noindent``One can choose $b_{ij}$ for
$\min(i,j)<h$ and $b_{hj}$ for $1\le j\le\max(n-h,h)$ so that $\tilde a_{ij}$ is as
desired for all $1\le i\le h$ and $1\le j\le n$''. 
\bigskip
}
\noindent Note that assertion $\ca A_h$ does not
say anything about the values of $b_{hk}$ when $k>\max(n-h,h)$. 

For $h=1$, it suffices to choose $b_{1k}=a_{1,k+1}$ and $b_{h1}=a_{h+1,1}$ for
$h$,~$k=2,\ldots,n$; in particular, $b_{1n}$ and~$b_{n1}$ are momentarily free, and $\ca
A_1$ is proved.

Now assume that $\ca A_{h-1}$ holds, with $h\ge 2$. By symmetry, we have already dealt with
the $\tilde a_{hk}$ for~$k<h$; we must only worry about the $\tilde a_{hk}$ with $k\ge h$.
There are two cases to consider.
\medskip

\item{(a)}{\it $n-h+1\ge h$, that is $h\le\lceil n/2\rceil$} (and
$\max(n-(h-1),h-1)=n-h+1$). In this case, when $k=h$ all addends in Eq.~(2.7) are already
determined, and thus the value of $\tilde a_{hh}$ cannot be anymore altered. If
$h<k\le n-h+1$ we can get $\tilde a_{hk}=0$ by choosing $b_{h,k-1}$ suitably. If $n-h+1<
k\le n$, on the other hand, we get $\tilde a_{hk}=0$ by choosing $b_{h-1,k}$ suitably,
keeping $b_{h,k-1}$ free. Summing up, we have proved $\ca A_h$ in this case, keeping
$b_{hk}$ free for $k\ge n-h=\max(n-h,h)$.

\item{(b)}{\it $n-h+1< h$, that is $h>\lceil n/2\rceil$} (and
$\max(n-(h-1),h-1)=h-1$). In this case, for any $k\ge h$ in Eq.~(2.7) we can choose
$b_{h-1,k}$ suitably so to get $\tilde a_{hk}=0$, even for $k=h$. In particular, we impose
no conditions on $b_{hk}$ for $k\ge h=\max(n-h,h)$, and $\ca A_h$ is proved in this case
too.\qedn

\newthm Lemma \dtre: Let 
$$
\phe(z)=\sum_{j=1}^n \eps_j(z_j)^2
$$
be a diagonal quadratic form on~$\C^n$, and set $j_0=\min\{j\mid \eps_j\neq 0\}$. Then there
exist a quadratic form $\psi$ on~$\C^n$ and $\alpha\in\C^*\times\C^{n-1}$ such that setting 
$$
\tilde\phe=\phe\circ J_{(\alpha)}-\psi\circ J'-\Psi\circ(I\oplus J')-\Psi\circ(J'\oplus I)
$$ 
we have
$$
\tilde\phe(z)=\cases{\alpha_0^2\,\eps_{j_0}(z_{j_0})^2& if $1\le j_0\le\lceil
	n/2\rceil$,\cr
	0&if $\lceil n/2\rceil<j_0\le n$.\cr}
$$

\pf 
Let us write again $\psi(z)=z^T Bz$, and $\tilde\phe(z)=z^T\tilde A z$, where
$B=(b_{hk})$, $\tilde A=(\tilde a_{hk})\in M_{n,n}(\C)$ are symmetric. This time we have
$$
\cases{\tilde a_{11}=\alpha_0^2 \eps_1\;,\cr
	\tilde a_{1k}=\alpha_0\alpha_{k-1}\eps_1-b_{1,k-1}& for $k=2,\ldots,n$,\cr
	\tilde a_{h1}=\alpha_{h-1}\alpha_0\eps_1-b_{h-1,1}& for $h=2,\ldots,n$,\cr
	\displaystyle\tilde a_{hk}=\sum_{j=1}^{\min(h,k)}\!\!\!\alpha_{h-j}\alpha_{k-j}\eps_j  
	-b_{h-1,k-1}-b_{h,k-1}-b_{h-1,k}& for $h$,~$k=2,\ldots,n$.\cr}
\eqno (2.8)
$$
We want to kill all the $\tilde a_{hk}$ except possibly $\tilde a_{j_0,j_0}$ if $j_0\le
\lceil n/2\rceil$. To do so, we shall prove, by induction on $h$, the following
assertion~$\ca B_h$: 
\bigskip
{\narrower \noindent``One can choose $b_{ij}$ for $\min(i,j)<h$, $b_{hk}$ for $1\le
k\le\max(n-h,h)$, and $\alpha_{2j}$ for $1\le j\le h-j_0$ so that $\tilde a_{ij}$ is as
desired for all~$1\le i\le h$ and~$1\le j\le n$. Furthermore, this can be done satisfying
the following requirements:
\medskip
\item{--} if $h<j_0$ then $b_{hk}=0$ for $h\le k\le n-h$;
\item{--} if $h\ge j_0$ and $j_0\le\lceil n/2\rceil$, then $b_{hk}$ for $h\le k\le n-h$ is a
linear combination of terms of the form~$\alpha_i \alpha_j \eps_l$, where the highest among
the indexes of the $\alpha$'s is $h+k-2j_0+1$, and it appears only once in an addend of the
form $\pm\alpha_0\alpha_{h+k-2j_0+1}\eps_{j_0}$.''
\bigskip
}
\noindent  Again, notice that we are imposing no conditions on~$b_{hk}$ for
$k>\max(n-h,h)$, as well as on $\alpha_0$, on~$\alpha_{2j-1}$ for~$1\le j\le h-j_0$ and on
$\alpha_j$ for $2(h-j_0)<j\le n-1$.

For $h=1$, assertion $\ca B_1$ follows trivially from Eq.~(2.8) for all values of~$j_0$.
So assume that $\ca B_{h-1}$ holds (for~$h\ge 2$); as usual, we must only deal with $\tilde
a_{hk}$ for $k\ge h$. There are four cases to consider.
\medskip

\item{(a)} {\it $h\le n-h+1$, that is $h\le\lceil n/2\rceil$, and $h<j_0$.} In this case, if
$h\le k\le n-h+1$ we have $b_{h-1,k}=b_{h-1,k-1}=0$ by $\ca B_{h-1}$, and so $\tilde
a_{hk}=0$ iff $b_{h,k-1}=0$. On the other hand, if $n-h+1<k\le n$ then we can get $\tilde
a_{hk}=0$ by suitably choosing $b_{h-1,k}$, leaving $b_{h,k-1}$ free. So $\ca B_h$ is proved
in this case.

\item{(b)} {\it $h\le n-h+1$, and $h=j_0$ (in particular, $j_0\le\lceil n/2\rceil$).} First
of all we see that $\ca B_{h-1}$ implies $\tilde a_{j_0,j_0}=\alpha_0^2\,\eps_{j_0}$; this
is the only case we cannot kill $\tilde a_{hk}$. If $h=j_0<k\le n-j_0+1$ we get $\tilde
a_{j_0,k}=0$ by choosing $b_{j_0,k-1}=\alpha_0\alpha_{k-j_0}\eps_{j_0}$; finally, if
$n-j_0+1<k\le n$ we get $\tilde a_{j_0,k}=0$ by suitably choosing $b_{j_0-1,k}$. So $\ca
B_{j_0}$ is verified.

\item{(c)} {\it $h\le n-h+1$, and $h>j_0$ (in particular, $j_0\le\lceil n/2\rceil$).} First
of all we have
$$
	\tilde a_{hh}=\sum_{j=j_0}^h \alpha_{h-j}^2\eps_j-b_{h-1,h-1}-2b_{h-1,h}\;.
$$
The inductive hypothesis $\ca B_{h-1}$ says that all the addends in this sum are linear
combinations of terms of the form $\alpha_i\alpha_j\eps_l$, and that the highest index
appearing in some $\alpha$'s is $2(h-j_0)$, in an addend of the form~$\pm\alpha_0
\alpha_{2(h-j_0)}\eps_{j_0}$; therefore we can choose $\alpha_{2(h-j_0)}$ so that $\tilde
a_{hh}=0$. Next, if $h<k\le n-h+1$, we can choose $b_{h,k-1}$ so to get $\tilde a_{hk}=0$;
looking at Eq.~(2.8) and recalling $\ca B_{h-1}$ one immediately sees that $b_{h,k-1}$ is
a linear combination of terms of the form $\alpha_i\alpha_j\eps_l$, and that the highest
index appearing is $h+k-2j_0$ in a single addend of the form $\pm\alpha_0\alpha_{h+k-2j_0}
\eps_{j_0}$ (it comes from $b_{h-1,k}$). Finally, if $n-h+1<k\le n$ we get $\tilde
a_{hk}=0$ by suitably choosing $b_{h-1,k}$, and $\ca B_h$ is again verified.

\item{(d)} {\it $h> n-h+1$, that is $h>\lceil n/2\rceil$.} The argument used in part (b) of
the proof of the previous lemma works in this case too, and we are done.
\qedn

We can now prove Theorem~2.1. To find $Q_1$ and
$Q_2$ such that $\tilde P_2$ is in normal form, write
$Q_1=J_{(\alpha)}$ for $\alpha\in\C^*\times\C^{n-1}$, and $Q_2=(\psi_1,\ldots,\psi_n)$, where
$\psi_j$ is a quadratic form on~$\C^n$ with associated symmetric bilinear form~$\Psi_j$.
Then it is easy to check that
$$
(P_1\circ Q_2-Q_2\circ P_1)_j=\cases{\psi_{j+1}-\psi_j\circ J'-\Psi_j\circ(I\oplus J')
	-\Psi_j\circ(J'\oplus I)& if $1\le j\le n-1$,\cr
	-\psi_n\circ J'-\Psi_n\circ(I\oplus J')-\Psi_n\circ(J'\oplus I)& if $j=n$,\cr}
$$
and that
$$
\tilde\phe_n=\alpha_0^{-1}\bigl[\phe_n\circ J_{(\alpha)}-\psi_n\circ J'-\Psi_n\circ(I\oplus
	J')-\Psi_n\circ(J'\oplus I)\bigr]\;,
$$
where we are using the notations introduced in the statements of Lemmas~2.2 and~2.3.

We now proceed as follows. Using Lemma~2.2, we make a first change of coordinates with
$Q_1=I$ and $\psi_j\equiv 0$ if $j\neq n$, so to put $\tilde\phe_n$ in
diagonal form. Next, using Lemma~2.3, we make a second change of coordinates again with
$\psi_j\equiv 0$ if $j\neq n$, so to put $\tilde\phe_n$ in normal form.
Finally, starting with $h=n-1$ and going up, we apply $n-1$ times Lemma~2.2, making change
of coordinates with $Q_1=I$ and $\psi_j\equiv 0$ if $j\neq h$, so to put $\tilde\phe_h$ in
normal form.\qedn

{\it Remark 2.1:}
A theorem of Weickert [14] and Forstneric [5] implies that for any germ
$\chi\in\Aut(\C^n,O)$ and any $k\in\N$ there is a {\it global} biholomorphism $h$ of $\C^n$
such that $\chi-h=O_k$. In particular, then, Theorem~2.1 holds for global holomorphic
self-maps of~$\C^n$ with a Jordan fixed point at the origin.
\medbreak

The normal form so obtained is not unique; there is still some freedom in the choice of the
$\eps^j_i$. However, we are able to extract from it a projective invariant:

\newthm Proposition \dqua: Let $F\in\End(\C^n,O)$ be with a Jordan fixed point at the
origin. Then the vector $(\eps^1_1,\ldots,\eps^1_n)$ in a normal form of $F$ is uniquely
determined up to a constant factor, and all possible factors do arise.

\pf
To prove the assertion we may assume that both $F$ and $\tilde F=\chi^{-1}\circ F\circ\chi$
are in normal form. Since
$$
(\tilde\eps^1_1,\ldots,\tilde\eps^1_n)=\tilde P_2(e_1)\;,
$$
where $e_1$ is the first vector of the canonical basis of $\C^n$, and analogously for $F$,
it suffices to show that $\tilde P_2(e_1)=\alpha_0 P_2(e_1)$. 

We know that
$$
Q_1\circ\tilde P_2=P_1\circ Q_2-Q_2\circ P_1+P_2\circ Q_1\;,
$$
which is equivalent to
$$
\sum_{h=i}^n \alpha_{h-i}\left(\sum_{j=1}^{\lceil n/2
\rceil}\tilde\eps^j_h(z_j)^2\right)=\psi_{i+1}(z)-\psi_i\bigl(J'(z)\bigr)-\Psi_i\bigl(z,J'(z)\bigr)-\Psi_i\bigl(J'(z),z
	\bigr)+\sum_{h=1}^{\lceil n/2 \rceil} \eps_i^h\left([Q_1(z)]_h\right)^2\;,
\eqno (2.9)
$$
for $i=1,\ldots,n$, where $\psi_{n+1}\equiv0$. Set $\psi_i(z)=z^T B^{(i)} z$, where
$B^{(i)}=(b_{hk}^{(i)})\in M_{n,n}(\C)$ is symmetric, and let~$\ca C_i$ denote the following
statement:
\bigskip
{\narrower\noindent ``We have $\tilde\eps^1_j=\alpha_0\eps^1_j$ for $j=i,\ldots,n$, and
$\displaystyle b_{1j}^{(i)}=\alpha_0\sum_{h=i}^n \alpha_{j+h-i}\eps^1_h $
for $j=1,\ldots,i-1$.''
\bigskip}
\noindent In particular, $\ca C_1$ is our assertion. We shall prove $\ca C_i$ by induction
on~$n-i$.

Putting $i=n$ and $z=e_1$ in Eq.~(2.9) we immediately get $\tilde\eps_n^1=\alpha_0\eps_n^1$.
On the other hand, computing the coefficient of $z_1z_{j+1}$ for $i=n$ and
$j=1,\ldots,n-1$ we get
$$
0=-b_{1j}^{(n)}+\alpha_0\alpha_j\eps^1_n\;,
\eqno (2.10)
$$
and so $\ca C_n$ is proved. Now assume $\ca C_{i+1}$ holds (for $1\le i<n$). Computing in
Eq.~(2.9) the coefficient of $z_1^2$ and using $\ca C_{i+1}$ we get
$$
\alpha_0\tilde\eps^1_i+\alpha_0\sum_{h=i+1}^n\alpha_{h-i}\eps^1_h=\sum_{h=i}^n
	\alpha_{h-i}\tilde\eps^1_h=b_{11}^{(i+1)}+\alpha_0^2\eps_i^1=\alpha_0\sum_{h=i+1}^n
	\alpha_{h-i}\eps_h^1+\alpha_0^2\eps_1^i\;,
$$
and thus $\tilde\eps_i^1=\alpha_0\eps_i^1$. Computing the coefficient of $z_1z_{j+1}$ for
$j=1,\ldots,i-1$ we get
$$
0=b_{1,j+1}^{(i+1)}-b_{1j}^{(i)}+\alpha_0\alpha_j\eps^1_i\;,
$$ 
and recalling $\ca C_{i+1}$ we get the assertion.

Finally, if $\chi(z)=\alpha_0 z$ it is clear that $\tilde P_2(e_1)=\alpha_0 P_2(e_1)$, and
we are done.\qedn

{\it Remark 2.2:}
In this paper our main concern will be maps with $(\eps^1_1,\ldots,\eps^1_n)\neq O$. Anyway,
using a similar argument one can prove that if there is a $j_0\le\lceil n/2\rceil$ such that
$\eps_i^j=0$ for $j<j_0$ and $i=1,\ldots,n$ then the vector $(\eps^{j_0}_{2j_0-1},\ldots,
\eps^{j_0}_n)$ is uniquely determined up to a constant factor. 
\medbreak

{\it Remark 2.3:}
If $\tilde F$ is a normal form of $F$, we have
$$
\tilde P_2(e_1)=Q_1^{-1}\Bigl(J'\bigl(Q_2(e_1)\bigr)+\alpha_0^2\, P_2(e_1)\Bigr)\;,
$$
and we get $\eps_n^1=\tilde\phe_n(e_1)=\alpha_0\phe_n(e_1)$. In particular,
$\eps_n^1\neq 0$ iff $\phe_n(e_1)\neq 0$, which is the generic case.
	
If $\phe_n(e_1)=0$, the same formula yields
$\eps^1_{n-1}=\alpha_0^{-1}\psi_n(e_1)+\alpha_0\phe_{n-1}(e_1)$. Examining the proof of
Lemmas~2.2 and~2.3 it is easy to check that $\psi_n(e_1)=\alpha_0^2 a_{12}^{(n)}$, where we
are writing $\phe_j(z)=z^T A^{(j)} z$. Therefore
$$
\eps^1_{n-1}=\alpha_0(a_{12}^{(n)}+a_{11}^{(n-1)})
$$
if $\eps^1_n=0$.
\medbreak

In Section~5 we shall also need to know the behavior under change of coordinates of part of
the third-order term in the expansion of~$F$ when $n=2$; to be precise, we shall need to
compute~$\tilde P_3(e_1)_2$. 

For the rest of this section we then assume~$n=2$ and set $\eta_2=P_3(e_1)_2$. Let $F$ and
$\tilde F=\chi^{-1}\circ F\circ\chi$ be both in normal form, and with~$\eps^1_2=0$. In
Proposition~2.4 we have already seen that $\tilde\eps^1_1=\alpha_0\eps^1_1$; furthermore,
recalling Eq.~(2.10) we also get~$\psi_2(e_1)=0$, because $\eps^1_2=0$.

On the other hand, Eq.~(2.6) yields
$$
\tilde\eta_2=\alpha_0^2\eta_2-\tilde\eps^1_1{\psi_2(e_1)\over\alpha_0}=\alpha_0^2\eta_2\;.
$$
Thus neither $\eps^1_1$ nor $\eta_2$ are invariant under change of coordinates, but if
$\eps^1_1\neq 0$ then $\eta_2/(\eps^1_1)^2$ is. For reasons that will become clear in
Section~5 we associates to a map with~$\eps^1_2=0$ and~$\eps^1_1\neq 0$ the invariant
$$
\Xi=1+{2\eta_2\over(\eps^1_1)^2}\;.
$$
A natural question is how to compute~$\Xi$. Let then $F$ be not necessarily in normal form,
and write again $\phe_j(z)=z^T A^{(j)} z$, where $A^{(j)}=(a^{(j)}_{hk})\in M_{2,2}(\C)$ is
symmetric. We have already remarked that $\eps^1_2=0$ is equivalent to~$a^{(2)}_{11}=0$. Let
$\tilde F=\chi^{-1}\circ F\circ\chi$, and write $\tilde
P_2=(\tilde\phe_1,\tilde\phe_2)$ and so on. If $\eps^1_2=0$ we get
$$
\tilde P_3(e_1)_2=\alpha^2_0 P_3(e_1)_2+2(a^{(2)}_{12}-a^{(1)}_{11})b^{(2)}_{11}-{2\over
	\alpha_0^2}(b^{(2)}_{11})^2\;,
$$
where $\psi_2(z)=z^T B^{(2)} z$ is the second component of the quadratic part of~$\chi$. 
Since in general
$$
\tilde a^{(2)}_{12}=\alpha_1 a^{(2)}_{11}+\alpha_0
	a^{(2)}_{12}-{1\over\alpha_0}b^{(2)}_{11}\;,
$$
if $\eps^1_2=0$ and $\tilde F$ is in normal form we must have $b^{(2)}_{11}=\alpha_0^2
a^{(2)}_{12}$ and thus
$$
\tilde P_3(e_1)_2=\alpha_0^2\bigl[P_3(e_1)_2-2a^{(1)}_{11}a^{(2)}_{12}\bigr]\;.
$$

Summing up we have proved the following:

\newthm Proposition \dcin: 
Let $F\in\End(\C^2,O)$ be with a Jordan fixed point at the origin. Assume that
$\eps^1_2=0$. Then
$$
\eps^1_1=\alpha_0(a^{(1)}_{11}+a^{(2)}_{12})\qquad\hbox{and}\qquad
	\eta_2=\alpha^2_0\bigl[P_3(e_1)_2-2a^{(1)}_{11}a^{(2)}_{12}\bigr]\;.
$$
In particular, $\eps^1_1\neq 0$ iff $a^{(1)}_{11}+a^{(2)}_{12}\neq 0$, and in that case
$$
\Xi={\bigl(a^{(1)}_{11}-a^{(2)}_{12}\bigr)^2+2P_3(e_1)_2\over
	\bigl(a^{(1)}_{11}+a^{(2)}_{12}\bigr)^2}\;.
$$

\smallsect 3. Blowing up maps along submanifolds

The main tool for the study of maps in normal form will be a suitable sequence of blow-ups,
used to ``resolve the singularity''. In this section we shall describe the general setting.

First of all, we fix a number of notations. If $z=(z_1,\ldots,z_n)\in\C^n$ and $0<r<n$ we
shall write $z'=(z_1,\ldots,z_r)$ and $z''=(z_{r+1},\ldots,z_n)$. If $r=0$, we
set $z''=z$, and $z'$ is empty. Finally, if $V$ is any vector space and $v\in
V\setminus\{O\}$, we denote by $[v]$ the projection of~$v$ in~$\P(V)$.

Let $M$ be a complex manifold of dimension $n\ge 2$, and $X\subset M$ a closed complex
submanifold of dimension $r\ge 0$ (it can be $0$). Let $N_{X/M}$ denote the normal bundle
of~$X$ in~$M$, and let $E_X=\P(N_{X/M})$ be the {\it projective normal bundle,} whose fiber
over $p\in X$ is $E_p=\P(T_pM/T_pX)$. The {\it blow-up of~$M$ along~$X$} is the set
$$
\tilde M_X=(M\setminus X)\cup E_X\;,
$$
endowed with the manifold structure we shall presently describe, together with the
projection $\sigma\colon\tilde M_X\to M$ given by $\sigma|_{M\setminus X}
=\id_{M\setminus X}$ and $\sigma|_{E_p}\equiv\{p\}$ for $p\in X$. The set
$E_X=\sigma^{-1}(X)$ is the {\it exceptional divisor} of the blow-up.

A chart $\phe=(z_1,\ldots,z_n)\colon V\to\C^n$ is {\it adapted to~$X$} if $V\cap
X=\{z_{r+1}=\cdots=z_n=0\}$. Choose a chart $(V,\phe)$ adapted to~$X$, and for
$j=r+1,\ldots,n$ and $q\in V\cap X$ set $X_j=\{z_j=0\}\subset V$,
$L_{j,q}=\P\bigl(\Ker(dz_j)_q /T_qX\bigr)\subset E_q$, $L_j=\bigcup_{q\in V\cap X}L_{j,q}$,
$E_{V\cap X}=\sigma^{-1}(V\cap X)$ and $V_j=(V\setminus X_j)\cup (E_{V\cap X}\setminus L_j)$.
Define $\chi_j\colon V_j\to\C^n$ by
$$
\eqalign{\chi_j(q)=\left(\phe(q)',{z_{r+1}(q)\over z_j(q)},\ldots,z_j(q),\ldots,
	{z_n(q)\over z_j(q)}\right)&\quad\hbox{\rm if $q\in V\setminus X_j\;$,}\cr
	\chi_j([v])=\left(\phe\bigl(\sigma([v])\bigr)',{d(z_{r+1})_{\sigma([v])}(v)\over
	d(z_j)_{\sigma([v])}(v)},\ldots,0,\ldots,{d(z_n)_{\sigma([v])}(v)\over
	d(z_j)_{\sigma([v])}(v)}\right)&\quad \hbox{\rm if $[v]\in E_{V\cap X}\setminus L_j\;$.}\cr}
\eqno (3.1)
$$
Then it is not difficult to check (see [6, pp.~602--604] for details, though the
presentation there is slightly different) that the charts $(V_j,\chi_j)$, together with an
atlas of~$M\setminus X$, endow $\tilde M_X$ with a structure of $n$-dimensional complex
manifold, as claimed, such that the projection~$\sigma$ is holomorphic everywhere.

For future reference, we record here that
$$
\eqalign{\chi_j^{-1}(w)=\phe^{-1}(w',w_jw_{r+1},\ldots,w_j,\ldots,w_jw_n)&\quad\hbox{if
	$w_j\neq 0$,}\cr
	\chi_j^{-1}(w)=\biggl[w_{r+1}\left.{\de\over\de z_{r+1}}\right|_q+\cdots+
	1\cdot\left.{\de\over\de z_j}\right|_q+\cdots +w_n \left.{\de\over\de
	z_n}\right|_q+T_qX\biggr]&\quad\hbox{if $w_j=0$,}\cr}
\eqno (3.2)
$$
where $q=\phe^{-1}(w',0'')$, and that
$$
\phe\circ\sigma\circ\chi_j^{-1}(w)=(w',w_jw_{r+1},\ldots,w_j,\ldots,w_jw_n)\;.
\eqno (3.3)
$$
The fiber~$E_p$ of the exceptional divisor over a point~$p\in X$ is a projective space; the
choice of an adapted chart yields an explicit isomorphism with~$\P^{n-r-1}(\C)$. To wit, we
define the isomorphism $\iota_{p,\phe}\colon E_p\to\P^{n-r-1}(\C)$ by setting
$$
\iota_{p,\phe}\left(\left[v_{r+1}\left.{\de\over\de z_{r+1}}\right|_p+\cdots+v_n
	\left.{\de\over\de z_{r+1}}\right|_p+T_pX\right]\right)=[v_{r+1}:\cdots:v_n]\;.
$$

We are interested to see when a germ $F\in\End(M,X)$ can be extended to the
blow-up~$\tilde M_X$. Take $p\in X$, and choose charts $(V,\phe)$ and $(\tilde
V,\tilde\phe)$ adapted to~$X$ so that $p\in V$ and $F(p)\in\tilde V$. In a neighbourhood
of~$p$ we can write the homogeneous expansion of $G=\tilde\phe\circ F\circ\phe^{-1}$ as
$$
G(z)=\sum_{l\ge 0}P_{l,z'}(z'')\;,
$$
where $P_{l,z'}$ is a $n$-uple of $l$-homogeneous polynomials with coefficients holomorphic
in~$z'$. The condition $F(X)\subseteq X$ then translates to
$$
(P_{0,z'})''\equiv 0\;.
$$
The {\it order of~$F$ at~$p$ along~$X$} is 
$$
\nu_X(F,p)=\min\{l\mid (P_{l,\phe(p)'})''\not\equiv 0\}\ge 1\;;
$$
it is easily checked that $\nu_X(F,p)$ does not depend on the adapted charts chosen. The {\it
order of~$F$ along~$X$} is then given by
$$
\nu_X(F)=\min\{\nu_X(F,p)\mid p\in X\}\;.
$$
Clearly the set $\{p\in X\mid \nu_X(F,p)=\nu_X(F)\}$ is open in~$X$.

We shall say that $F$ is {\it non-degenerate at~$p$ along~$X$} if 
\medskip

\item{(i)} $F^{-1}(p)\subseteq X$,
\item{(ii)} $\nu_X(F,p)=\nu_X(F)$, and
\item{(iii)} $\bigl(P_{l_0,\phe'(p)}(v)\bigr)''=0$ iff $v=O\in\C^{n-r}$, where
$l_0=\nu_X(F)$. 

\medskip
\noindent If $F$ is non-degenerate along~$X$ at all points of~$X$ we shall say that $F$ is
{\it non-degenerate along~$X$.} 

\newthm Theorem \tuno:
Let $M$ be a complex manifold of dimension $n$, and $X\subset M$ a closed submanifold of
dimension~$r\ge 0$. Let $F\in\End(M,X)$ be non-degenerate along $X$. Then there exists a
unique $\tilde F\in\End(\tilde M_X,E_X)$ such that $F\circ\sigma=\sigma\circ\tilde F$.
Furthermore, if $p\in X$ and $(V,\phe)$, $(\tilde V,\tilde\phe)$ are charts adapted to~$X$
with~$p\in V$ and~$F(p)\in\tilde V$, then
$$
\tilde F\bigl([v]\bigr)=(\iota_{F(p),\tilde\phe})^{-1}\left(\left[P_{l_0,\phe(p)'}
	\bigl(\iota_{p,\phe}([v])\bigr)''\right]\right)
\eqno (3.4)
$$
for all $[v]\in E_p$, where $l_0=\nu_X(F)$.

\pf
Since $F^{-1}(X)\subseteq X$, if $q$ does not belong to~$X$ we can safely set $\tilde
F(q)=F(q)$; we are left to define $\tilde F$ on the exceptional divisor.

Choose $p\in X$, and the charts as in the statement of the theorem. For $[v]\in E_p$ choose
$r+1\le j\le n$ so that $[v]\in V_j$; if $\tilde F$ exists, we must have
$$
F\circ\sigma\circ\chi_j^{-1}=\sigma\circ\tilde F\circ\chi_j^{-1}\;.
$$
If $[v]=(\iota_{p,\phe})^{-1}[v_{r+1}:\ldots:v_n]$, we have
$$
[v]=\lim_{\zeta\to 0}\chi_j^{-1}\left(\phe(p)',{v_{r+1}\over v_j},\ldots,\zeta,\ldots,
	{v_n\over v_j}\right)\;,
$$
and so, setting again $G=\tilde\phe\circ F\circ\phe^{-1}$, 
$$
\eqalign{\tilde F([v])&=\lim_{\zeta\to 0}\tilde F\circ\chi_j^{-1}\left(\phe(p)',
	{v_{r+1}\over v_j},\ldots,\zeta,\ldots,{v_n\over v_j}\right)\cr
	&=\lim_{\zeta\to 0}\sigma^{-1}\left(\tilde\phe^{-1}\left(G\left(\phe(p)',{\zeta\over v_j}v
	\right)\right)\right)\;,\cr}
$$
where with a slight abuse of notation we have put $v=(v_{r+1},\ldots,v_n)\in\C^{n-r}$. 

Now it is easy to check that if a sequence $\{q_k\}\subset M\setminus X$ converges
to~$q=F(p)\in X$, then the sequence~$\{\sigma^{-1}(q_k)\}$ converges in~$\tilde M
\setminus X$ iff $\{[\tilde\phe(q_k)'']\}$ converges in $\P^{n-r-1}(\C)$. In our case we have
$$
G\left(\phe(p)',{\zeta\over v_j}v\right)''=\sum_{l\ge l_0}P_{l,\phe(p)'}\left(
	{\zeta\over v_j}v\right)''=\left(\zeta\over v_j\right)^{l_0}\bigl(P_{l_0,\phe(p)'}(v)''+
	\zeta Q(\zeta)\bigr)\;,
$$
for a suitable holomorphic map $Q$. Therefore $[G(\zeta
v/v_j)'']\to[P_{l_0,\phe(p)'}(v)'']$, and thus if $\tilde F$ exists it is given by Eq.~(3.4)
on the exceptional divisor.

To finish the proof we must show that an $\tilde F$ defined by Eq.~(3.4) on the exceptional
divisor and by $F$ elsewhere is holomorphic. Take $[v]\in E_p$, and choose $r+1\le h,k\le n$
so that $[v]\in V_h$ and $\tilde F([v])\in\tilde V_k$; we must show that $\chi_k\circ\tilde
F\circ\chi_h^{-1}$ is holomorphic. We know that
$$
G\circ(\phe\circ\sigma\circ\chi_h^{-1})=(\phe\circ\sigma\circ\chi_k^{-1})\circ
	(\chi_k\circ\tilde F\circ\chi_h^{-1})\;;
$$
so putting $\chi_k\circ\tilde F\circ\chi_h^{-1}=(\tilde f_1,\ldots,\tilde f_n)$ and
recalling Eq.~(3.3) we must have
$$
G(w',w_hw_{r+1},\ldots,w_h,\ldots,w_hw_n)
	=\bigl(\tilde f_1(w),\ldots,\tilde f_r(w),
	\tilde f_k(w)\tilde f_{r+1}(w),\ldots,\tilde f_k(w),\ldots,\tilde f_k(w)\tilde f_n(w)
	\bigr)\;.
$$
Writing $G=(g_1,\ldots,g_n)$ we find that if $w_h\neq 0$ then
$$
\tilde f_i(w)=\cases{g_i(w',w_hw_{r+1},\ldots,w_h,\ldots,w_hw_n)& if $1\le i\le r$ or
	$i=k$,\cr
	\displaystyle {g_i(w',w_hw_{r+1},\ldots,w_h,\ldots,w_hw_n)\over
	g_k(w',w_hw_{r+1},\ldots,w_h,\ldots,w_hw_n)}& if $r+1\le i\neq k\le n$.\cr}
\eqno (3.5)
$$
Since the $g_i$'s are holomorphic and $\{w_h=0\}$ has codimension~1 in $\chi_h(V_h)$, to
end the proof it suffices to show that the quotients in Eq.~(3.5) have a limit when
$w\to\chi_h([v])$. 

Write again $\iota_{p,\phe}([v])=[v_{r+1}:\ldots:v_n]$ and $v=(v_{r+1},\ldots,v_n)$,
and assume then that $w\to\chi_h([v])$. This means that $w'\to\phe(p)'$,
$w_h\to 0$ and $(w_{r+1},\ldots,1,\ldots,w_n)\to v_h^{-1}v$. Now,
$$
\eqalign{G(w',w_hw_{r+1},\ldots,w_h,\ldots,w_hw_n)''
	&=\sum_{l\ge l_0}(w_h)^l P_{l,w'}
	(w_{r+1},\ldots,w_h,\ldots w_n)''\cr
	&=\sum_{l\ge l_0}\left({w_h\over v_h}\right)^l P_{l,w'}(w_{r+1}v_h,\ldots,v_h,\ldots,
	w_n v_h)''\;.\cr}
$$
Since $\tilde F([v])\in\tilde V_k$ we have $P_{l_0,\phe(p)'}(v)_k\neq 0$; therefore
$$
{g_i(w',w_hw_{r+1},\ldots,w_h,\ldots,w_hw_n)\over
	g_k(w',w_hw_{r+1},\ldots,w_h,\ldots,w_hw_n)}\to {P_{l_0,\phe(p)'}(v)_i\over
	P_{l_0,\phe(p)'}(v)_k}\;,
$$
and we are done.\qedn

{\it Remark 3.1:} If $F$ is degenerate along $X$ at some point $p\in X$, then we can still
define~$\tilde F$ by Eq.~(3.4), but the resulting map is not anymore continuous.
\medbreak

As mentioned in the introduction, our construction involves iterated blow-ups; thus we are
interested to know when the map $\tilde F$ is still non-degenerate along a suitable
submanifold of~$\tilde M_X$. We shall limit ourselves to two special cases, which are enough
for our aims. To state the first one, we need the following definition: if
$\sigma\colon\tilde M_X\to M$ is the blow-up of $M$ along~$X$, and $Y\subseteq M$ is a
submanifold of~$M$, then the {\it proper} (or {\it strict\/}) {\it transform} of~$Y$ is
$\tilde Y=\overline{\sigma^{-1}(Y\setminus X)}$. 

\newthm Proposition \tdue:
Let $M$ be a complex manifold of dimension $n$, and $X\subset M$ a closed submanifold of
dimension~$r\ge 0$. Let $F\in\End(M,X)$ be non-degenerate along $X$, and $\tilde
F\in\End(\tilde M_X,E_X)$ its blow-up. Let~$Y\subseteq M$ be a submanifold of~$M$ of
dimension~$r+s$ (with $s\ge1$), and $\tilde Y\subseteq\tilde M$ be its proper transform.
Assume that
\medskip
{
\item{\rm (i)} $Y$ contains properly $X$;
\item{\rm (ii)} $F(Y)\subseteq Y$ and $F^{-1}(Y)\subseteq Y$;
\item{\rm (iii)} $dF_q$ is invertible for all~$q\in Y$.
\medskip
\noindent Then $\tilde F$ is non-degenerate along $\tilde Y$, and $d\tilde F_{\tilde q}$ is
invertible for all~$\tilde q\in\tilde Y$.}

\pf
First of all, notice that if $p\in X$ then $\tilde Y\cap E_p=\P(T_pY/T_pX)$, and that $\tilde
F|_{E_p}$ is induced by $dF_p$. Since, by construction, $\tilde F(\tilde Y)\subseteq Y$ and
$\tilde F^{-1}(\tilde Y\setminus E_X)\subseteq\tilde Y\setminus E_X$, it suffices to prove
that $d\tilde F_{[v]}$ is invertible for all $[v]\in\tilde Y\cap E_X$; 

Fix $p\in X$ and $[v]\in\tilde Y\cap E_p$, and choose two charts $(V,\phe)$ and $(\tilde
V,\tilde\phe)$, centered in~$p$, respectively in~$F(p)$, such that
$$
V\cap X=\{z_{r+1}=\cdots=z_n=0\}\qquad\hbox{\rm and}\qquad V\cap Y
	=\{z_{r+s+1}=\cdots=z_n=0\}\;,
$$
and analogously for $\tilde V$. In particular,
$$
\iota_{p,\phe}(\tilde Y\cap E_p)=\iota_{F(p),\tilde\phe}(\tilde Y\cap E_{F(p)})=
	\{v_{r+s+1}=\cdots=v_n=0\}\;,
$$
and we can also assume that $\iota_{p,\phe}([v])=\iota_{F(p),\tilde\phe}\bigl(\tilde F([v])
\bigr)=[1:0:\cdots:0]$. Then the charts $(V_{r+1},\chi_{r+1})$ and $(\tilde V_{r+1},
\tilde\chi_{r+1})$ are centered in~$[v]$, respectively in~$\tilde F([v])$, and adapted
to~$\tilde Y$ (up to a permutation of the coordinates).

Set $G=\tilde\phe\circ F\circ\phe^{-1}=(g_1,\ldots,g_n)$ and $\tilde G=\tilde\chi_{r+1}
\circ\tilde F\circ\chi_{r+1}^{-1}=(\tilde f_1,\ldots,\tilde f_n)$; the relation between
the~$g_i$'s and the $\tilde f_j$'s is given by Eq.~(3.5). The condition $F(X)\subseteq X$
implies
$$
{\de g_i\over\de z_j}(O)=0\qquad\hbox{for $r+1\le i\le n$ and $1\le j\le r$,}
$$
while the condition $F(Y)\subseteq Y$ implies
$$
{\de g_i\over\de z_j}(O)=0\qquad\hbox{for $r+s+1\le i\le n$ and $r+1\le j\le r+s$.}
\eqno (3.6)
$$
This means that the jacobian matrix of $G$ at the origin is of the form
$$
\ca A=\left|\vcenter{\offinterlineskip
\halign{&\hfil$#$\hfil\cr
\multispan3\strut&\omit\hskip1.5mm\vrule width 1pt\hskip 1.5mm&\cr
\multispan3\strut\hskip1.5mm\hfil $A$\hfil&\omit\hskip1.5mm\vrule width 1pt
	\hskip 1.5mm&\multispan3\hfil$*$\hskip1.5mm\hfil\cr
\multispan3\strut&\omit\hskip1.5mm\vrule width 1pt\hskip 1.5mm&\cr
\noalign{\hrule height 1pt}
\multispan3&\omit\hskip1.5mm\vrule height 1pt width 1pt\hskip 1.5mm&\multispan3\cr
\strut&&&\omit\hskip1.5mm\vrule width 1pt\hskip 1.5mm&B&\omit\hskip1.5mm
	\vrule\hskip1.5mm&*\hskip1.5mm\cr
\multispan3\hskip1.5mm\quad\hfil$\smash{\lower 3pt\vbox{\smash{$O$}}}$\quad
	\hfil&\omit\hskip1.5mm\vrule width 1pt&\multispan3\hskip-1mm\hrulefill\cr
\multispan3&\omit\hskip1.5mm\vrule height 1pt width 1pt\hskip 1.5mm&\omit&\omit
\hskip1.5mm\vrule height 1pt\hskip1.5mm&\omit\cr
\strut&&&\omit\hskip1.5mm\vrule width 1pt\hskip 1.5mm&O&\omit\hskip1.5mm\vrule
	\hskip1.5mm&C\hskip1.5mm\cr}}
\right|\;,
$$
with $A\in M_{r,r}(\C)$, $B\in M_{s,s}(\C)$ and $C\in M_{n-r-s,n-r-s}(\C)$. Since, by
assumption, $dF_p$ is invertible, we have
$$
\det(\ca A)=\det(A)\det(B)\det(C)\neq 0\;.
$$
Finally, $\tilde F([v])\in\tilde V_{r+1}$ translates in
$$
\lambda={\de g_{r+1}\over\de z_{r+1}}(O)\neq 0\;.
$$
Our aim is to compute $\de\tilde f_i/\de w_j$ at $w=O$. This is easy when $1\le i\le r+1$;
in fact, Eq.~(3.5) with $h=k=r+1$ yields
$$
{\de\tilde f_i\over\de w_j}(O)=\cases{\displaystyle{\de g_i\over\de z_j}(O)& for $1\le i\le
	r+1$, $1\le j\le r+1$,\cr
	0& for $1\le i\le r+1$, $r+2\le j\le n$.\cr}
$$
In particular,
$$
{\de\tilde f_{r+1}\over\de w_j}(O)=\cases{0& if $j\neq r+1$,\cr
	\lambda\neq 0& if $j=r+1$.\cr}
$$
Now set $\tilde g_i(w)=g_i(w',w_{r+1},w_{r+1}w_{r+2},\ldots,
w_{r+1}w_n)$, and write again
$$
G(z)=\sum_{l\ge 0}P_{l,z'}(z'')\;,
$$
recalling that $(P_{0,z'})''\equiv O$. For $r+2\le i\le n$ we have
$$
{\de\tilde f_i\over\de w_j}(O)=\lim_{w\to O}{1\over\tilde g_{r+1}(w)}\left[
	{\de\tilde g_i\over\de w_j}(w)-{\tilde g_i(w)\over\tilde g_{r+1}(w)}{\de\tilde g_{r+1}
	\over\de w_j}(w)\right]\;.
\eqno (3.7)
$$
Since
$$
\tilde g_i(w)=\sum_{l\ge 0}(w_{r+1})^l P_{l,w'}(1,w_{r+2},\ldots,w_n)_i\;,
$$
Eq.~(3.7) yields
$$
\eqalign{{\de\tilde f_i\over\de w_j}(O)={1\over\lambda}\left[{\de^2 g_i\over\de z_j
	\de z_{r+1}}(O)-{1\over\lambda}{\de^2 g_{r+1}\over\de z_j\de z_{r+1}}(O){\de g_i\over
	\de z_{r+1}}(O)\right]&\quad\hbox{for $r+2\le i\le n$ and $1\le j\le r+1$,}\cr
{\de\tilde f_i\over\de w_j}(O)={1\over\lambda}\left[{\de g_i\over\de z_j}(O)-{1\over\lambda}
	{\de g_{r+1}\over\de z_j}(O){\de g_i\over\de z_{r+1}}(O)\right]&\quad\hbox{for $r+2\le
	i,j\le n$.}\cr}
\eqno (3.8)
$$
In particular, Eq.~(3.6) yields
$$
{\de\tilde f_i\over\de w_j}(O)={1\over\lambda}{\de g_i\over\de z_j}(O)\quad\hbox{for
	$r+s+1\le i\le n$, $r+2\le j\le n$.}
$$
We have proved that the Jacobian matrix of $\tilde G$ at the origin is
$$
\tilde{\ca A}=\left|\vcenter{\offinterlineskip
\halign{&\hfil$#$\hfil\cr
\multispan3&\omit\hskip1.5mm\vrule height2mm\hskip 1.5mm&\omit&\omit
	\hskip1.5mm\vrule height2mm width 1pt\hskip 1.5mm\cr
\multispan3\strut\hskip1.5mm\hfil $A$\hfil&\omit\hskip1.5mm\vrule\hskip	1.5mm&*&\omit
	\hskip1.5mm\vrule width 1pt\hskip 1.5mm&\multispan3\hfil$O$\hskip1.5mm\hfil\cr
\multispan3&\omit\hskip1.5mm\vrule height2mm\hskip 1.5mm&\omit&\omit
	\hskip1.5mm\vrule height2mm width 1pt\hskip 1.5mm\cr
\noalign{\hrule}
\multispan3&\omit\hskip1.5mm\vrule height2pt\hskip 1.5mm&\omit&\omit
	\hskip1.5mm\vrule height2pt width 1pt\hskip 1.5mm\cr
\multispan3\strut\hskip1.5mm\hfil $O$\hfil&\omit\hskip1.5mm\vrule\hskip	1.5mm&\lambda&\omit
	\hskip1.5mm\vrule width 1pt\hskip 1.5mm&\multispan3\hfil$O$\hskip1.5mm\hfil\cr
\noalign{\hrule height 1pt}
\multispan3&\omit\hskip1.5mm\vrule height 2pt\hskip 1.5mm&\omit&\omit\hskip1.5mm\vrule
		height 2pt width 1pt\hskip 1.5mm&\multispan3\cr
\strut&&&\omit\hskip1.5mm\vrule \hskip 1.5mm&\omit&\omit\hskip1.5mm\vrule width 1pt\hskip
	1.5mm&\tilde B&\omit\hskip1.5mm\vrule\hskip1.5mm&\hfil$\smash{\raise
	2pt\vbox{\smash{$*$}}}$\hfil\hskip1.5mm\cr
\multispan3\hskip1.5mm\quad\hfil$\smash{\lower 3pt\vbox{\smash{$*$}}}$\quad
	\hfil&\omit\hskip1.5mm\vrule&\hfil$\smash{\lower 3pt\vbox{\smash{$*$}}}$
	\hfil&\omit\hskip1.5mm\vrule width 1pt&\multispan3\hskip-1mm\hrulefill\cr
\multispan3&\omit\hskip1.5mm\vrule height 1pt\hskip 1.5mm&\omit&\omit\hskip1.5mm\vrule
		height 1pt width 1pt\hskip 1.5mm&\omit&\omit\hskip1.5mm\vrule height 1pt\hskip 1.5mm&\cr
\strut&&&\omit\hskip1.5mm\vrule\hskip 1.5mm&\omit&\omit\hskip1.5mm\vrule width 1pt\hskip
	1.5mm&O&\omit\hskip1.5mm\vrule \hskip1.5mm&{1\over\lambda}C\hskip1.5mm\cr}}
\right|\;,
\eqno (3.9)
$$
where $\tilde B\in M_{s-1,s-1}(\C)$. Now, if we subtract to the $j$-th column of $B$ (for
$j=2,\ldots,s$) the first column of~$B$ multiplied by $\lambda^{-1}\de g_{r+1}/\de
z_{r+j}(O)$ we get
$$
\left|\vcenter{\offinterlineskip
\halign{&\strut\hfil$#$\hfil\cr
\hskip1.5mm\lambda&\omit\hskip1.5mm\vrule\hskip1.5mm&O\hskip1.5mm\cr
\noalign{\hrule}
\omit&\omit\hskip 1.5mm \vrule height2pt\hskip 1.5mm&\omit\cr
\hskip1.5mm*&\omit\hskip1.5mm\vrule\hskip1.5mm&\lambda\tilde B\hskip1.5mm\cr}}
\right|\;.
$$
Since these elementary operations do not change the determinant, we obtain
$$
\det(B)=\lambda^2\det(\tilde B)\;.
$$
Therefore
$$
\det(\tilde{\ca A})={1\over\lambda^2}\det(\ca A)\neq 0\;,
$$
and we are done.\qedn

A similar argument yields:

\newthm Proposition \ttre:
Let $M$ be a complex manifold of dimension $n$, and $X\subset M$ a closed submanifold of
dimension~$r\ge 0$. Let $F\in\End(M,X)$ be non-degenerate along $X$, and $\tilde
F\in\End(\tilde M_X,E_X)$ its blow-up. Take $p\in X$ and a linear subspace $L\subseteq E_p$
of dimension~$s-1$ (with $s\ge1$). Assume that
{\medskip
\item{\rm (i)} $\tilde F(L)\subseteq L$, and
\item{\rm (ii)} $dF_p$ is invertible.
\medskip
\noindent Then $\tilde F$ is non-degenerate along $L$, and $d\tilde F_{[v]}$ is
invertible for all~$[v]\in L$.}

\pf
Condition (i) implies that $p$ is a fixed point of $F$, and condition (ii) implies that
$\nu_X(F)=1$. In particular, $\tilde F|_{E_p}$ is induced by the differential of~$F$ at~$p$;
thus $\tilde F|_L$ is injective, and the invertibility of $d\tilde F_{[v]}$ for all~$[v]\in
L$ will imply that $\tilde F$ is non-degenerate along~$L$.

Fix $[v]\in L$, and choose two charts $(V,\phe)$, $(\tilde V,\tilde\phe)$ centered in $p$
adapted to~$X$ such that 
$$
\iota_{p,\phe}([v])=\iota_{p,\tilde\phe}\bigl(\tilde F([v])\bigr)=[1:0:\ldots:0]
$$ 
and
$$
\iota_{p,\phe}(L)=\iota_{p,\tilde\phe}(L)=\{v_{r+s+1}=\cdots=v_n=0\}\;.
$$
Then the charts $(V_{r+1},\chi_{r+1})$ and $(\tilde V_{r+1},
\tilde\chi_{r+1})$ are centered in~$[v]$, respectively in~$\tilde F([v])$, and adapted
to~$L$ (up to a permutation of the coordinates). The proof then goes on as in the previous
proposition.\qedn

In particular, the latter proposition applies when $L$ reduces to a single point, that must
be an eigenvector of the differential of~$F$ at~$p$. This case is important for dynamical
reasons:

\newthm Lemma \tqua:
Let $M$ be a complex manifold of dimension~$n$, and $p\in M$. Let $F\in\End(M,p)$ be
non-degenerate at~$p$, and $\tilde F\in\End(\tilde M_p,E_p)$ its
blow-up. Let $\{z^k=F^k(z^0)\}\subset M\setminus\{p\}$ be an $F$-orbit converging to~$p$.
Assume there is $v\in T_pM$ such that $\sigma^{-1}(z^k)\to[v]$. Then $\tilde F([v])=[v]$.

\pf
Indeed, since $\sigma^{-1}(z^k)=\sigma^{-1}\bigl(F^k(z^0)\bigr)=\tilde
F^k\bigl(\sigma^{-1}(z^0)
\bigr)$, we have
$$
\tilde F([v])=\lim_{k\to\infty}\tilde F\bigl(\tilde F^k\bigl(\sigma^{-1}(z^0)\bigr)\bigr)
	=\lim_{k\to\infty}\tilde F^{k+1}\bigl(\sigma^{-1}(z^0)\bigr)=\lim_{k\to\infty}\sigma^{-1}
	(z^{k+1})=[v]\;.
$$
\qedn

{\it Remark 3.2:}
In local coordinates, the hypothesis $\sigma^{-1}(z^k)\to[v]$ translates to
$[w^k]\to\iota_{p,\phe}([v])\in\P^{n-1}(\C)$, where $w^k=\phe(z^k)$ and $(V,\phe)$ is a
chart centered at~$p$. The assertion then becomes $P_{l_0}(v)=\lambda v$ for a suitable
$v\in\C^*$, where $l_0=\nu_p(F)$ and $P_{l_0}$ is the first non-zero term in the homogeneous
expansion of $\phe\circ F\circ\phe^{-1}$ at~$p$.

\smallsect 4. The desingularization theorem

We can now describe the sequence of blow-ups we shall need. Let $M^0=\C^n$,
$\chi_0=\id_{\C^n}$, ${\bf e}_0=O$ and~$X^0=\{O\}$. We start by blowing up the origin, and
thus by setting $M^1=\tilde M^0_{X^0}$ and $\pi_1=\sigma_1\colon M^1\to M^0$. Since
$M^0=\C^n$ has a canonical chart adapted to $X^0$ (that is, centered at the origin), the
exceptional divisor~$E^1=\pi_1^{-1}(X^0)$ is canonically isomorphic to $\P^{n-1}(\C)$. This
allows us to define a distinguished point~${\bf e}_1\in E^1$, corresponding
to~$[e_1]\in\P^{n-1}(\C)$, and also distinguished linear subspaces $Y^k\subset E^1$ for
$k=2,\ldots,n-1$, corresponding to the linear subspaces generated by
$\{[e_1],\ldots,[e_k]\}\subset\P^{n-1}(\C)$.

Now put $X^1=\{{\bf e}_1\}$ and set $M^2=\tilde M^1_{X^1}$. Let $X^2\subset M^2$ be the
proper transform of~$Y^2$, and set~$M^3=\tilde M^2_{X^2}$. Next, let $X^3\subset M^3$ be the
proper transform (with respect to $\sigma_3\colon M^3\to M^2$) of the proper transform (with
respect to $\sigma_2\colon M^2\to M^1$) of~$Y^3$, and put $M^4=\tilde M^3_{X^3}$. Proceeding
in this way, we define for $k=2,\ldots,n-1$ the manifold~$M^{k+1}$ as the blow-up of~$M^k$
along the iterated proper transform~$X^k$ of~$Y^k$, and we denote by $\sigma_{k+1}\colon
M^{k+1}\to M^k$ the associated projection, and by $E^k\subset M^k$ the exceptional divisor.

\medbreak
{\it Remark 4.1:} As we shall see in the proof of Theorem~4.3, the choice of this sequence of
blow-ups is dictated by the geometry of the maps we want to study. To be more precise, one
has to blow-up along the iterated proper transforms of a flag $Y^1\subset\cdots\subset
Y^{n-1}$ invariant under the action of the differential of the map.
\medbreak

At each stage of this construction there are canonical charts adapted to the submanifolds
involved:

\newthm Lemma \quno:
For $1\le k\le n$ we can find a distinguished point ${\bf e}_k\in M^k$ and a canonical
chart $(V_k,\chi_k)$ centered in~${\bf e}_k$ such that
$V_k\cap X^k=\chi_k^{-1}(\{w_1=w_{k+1}=\cdots=w_n=0\})$, that is
$(V_k,\chi_k)$ is adapted to~$X^k$ (up to a permutation of the coordinates). Furthermore,
$$
V_k\cap\pi_k^{-1}(X^0)=\chi_k^{-1}\left(\bigcup_{h=1}^k\{w_h=0\}\right)\supset V_k\cap X^k\;,
\eqno (4.1)
$$
where $\pi_k=\sigma_1\circ\cdots\circ\sigma_k$, and
$$
\chi_{k-1}\circ\sigma_k\circ\chi_k^{-1}(w)=
	\cases{(w_1,w_1w_2,\ldots,w_1w_n)& if $k=1$,\cr
	(w_1w_k,w_2,\ldots,w_k,w_kw_{k+1},\ldots,w_kw_n)& if $2\le k\le n$.\cr}
\eqno (4.2)
$$
Finally, for $h=k+1,\ldots,n-1$ the intersection of the iterated proper transform
of~$Y^h$ with~$V_k$ is 
$$
\chi_k^{-1}(\{w_1=w_{h+1}=\cdots=w_n=0\})\;.$$

\pf
For $k=1$, the existence of a canonical chart adapted to~$X^0$ yields, via the
construction described in the previous section, a canonical chart~$(V_1,\chi_1)$ adapted
to~$X^1$ --- and in turn this yields a canonical basis $\{\de/\de w_1,\ldots,\de/\de w_n\}$
of~$T_{{\bf e}_1}M^1$. Furthermore, it is easy to check that
$$
V_1\cap\pi_1^{-1}(X^0)=\chi_1^{-1}(\{w_1=0\})\supset V_1\cap X^1\;,
$$
and
$$
\chi_0\circ\sigma_1\circ\chi_1^{-1}(w)=(w_1,w_1w_2,\ldots,w_1w_n)\;.
$$
So the lemma is proved for $k=1$.

Assume, by induction, that the lemma holds for $k-1$. In particular, we have a
ditinguished point~${\bf e}_{k-1}$ and a canonical chart $(V_{k-1},\chi_{k-1})$ centered
at~${\bf e}_{k-1}$ and adapted to~$X^{k-1}$ (up to a permutation of the coordinates). We thus
have a canonical basis $\{\de/\de w_1,\ldots,\de/\de w_n\}$ of $T_{{\bf e}_{k-1}}M^{k-1}$
such that
$\{\de/\de w_2,\ldots,\de/\de w_{k-1}\}$ spans $T_{{\bf e}_{k-1}}X^{k-1}$. Put
$$
{\bf e}_k=\left[{\de\over\de w_k}+T_{{\bf e}_{k-1}}X^{k-1}\right]\in\sigma_k^{-1}({\bf
		e}_{k-1})\;,
$$
and let $(V_k,\chi_k)$ be the canonical chart centered in~${\bf e}_k$ constructed, as before,
via~$(V_{k-1},\chi_{k-1})$. Then it is not too difficult to check using the inductive
hypothesis that $(V_k,\chi_k)$ is as desired.\qedn

The set $\pi_k^{-1}(X^0)\subset M^k$ will be called the {\it singular divisor} of~$M^k$.

To prove the Desingularization Theorem~4.3 we need still another local computation:

\newthm Lemma \qdue:
Fix $1\le k\le n$, and write $z=\chi_0\circ\pi_k\circ\chi_k^{-1}(w)$. Then
$$
z_j=\cases{w_1\prod\limits_{h=2}^j (w_h)^2\,\prod\limits_{h=j+1}^k w_h& if $1\le j\le k$,\cr
	w_1\prod\limits_{h=2}^k (w_h)^2\, w_j&if $k+1\le j\le n$.\cr}
\eqno (4.3)
$$
Furthermore, if $z_1,\ldots,z_k\ne0$ then
$$
w_j=\cases{(z_1)^2/z_k& if $j=1$,\cr
	z_j/z_{j-1}& if $2\le j\le k$,\cr
	z_j/z_k& if $k+1\le j\le n$.\cr}
\eqno (4.4)
$$

\pf
We shall prove Eq.~(4.3) by induction on~$k$. For $k=1$ it follows from Eq.~(4.2). Assume it
holds for~$k-1$, and write $u=(\chi_{k-1}\circ\sigma_k\circ\chi_k^{-1})(w)$. Then
$z=(\chi_0\circ\pi_{k-1}\circ\chi_{k-1}^{-1})(u)$; therefore the inductive hypothesis
together with Eq.~(4.2) yield
$$
z_j=\cases{w_1w_k\prod\limits_{h=2}^j (w_h)^2\,\prod\limits_{h=j+1}^{k-1} w_h& if $1\le j\le
k-1$,\cr
	w_1w_k\prod\limits_{h=2}^{k-1} (w_h)^2\, w_k& if $j=k$,\cr
	w_1w_k\prod\limits_{h=2}^{k-1}(w_h)^2\,w_kw_j& if $k+1\le j\le n$,\cr}
$$
and Eq.~(4.3) follows. Finally, Eq.~(4.4) is a trivial consequence of Eq.~(4.3).\qedn

At last we are ready to prove the main result of this paper:

\newthm Theorem \qtre: (Desingularization Theorem)
Let $F\in\End(\C^n,O)$ be with a Jordan fixed point at the origin. Assume that $F$ is in
normal form. Then for $1\le k\le n$ there exists a unique $\tilde F_k\in\End(M^k,E^k)$ such
that $F\circ\pi_k=\pi_k\circ\tilde F_k$. Furthermore, ${\bf e}_k$ is a fixed point of~$\tilde
F_k$, and $\tilde G_k=(\tilde g_1,\ldots,\tilde g_n)=\chi_k\circ\tilde F_k\circ\chi_k^{-1}$
in a neighbourhood of the origin is given by
$$
\tilde g_j(w)=\cases{w_1[1+\eps^1_1w_1+w_2+O_2]& if $j=1$,\cr
	\eps^1_jw_1+w_j+w_{j+1}\cr
	\phantom{\eps^1_jw_1}-(\eps^1_1w_1+w_2)(\eps^1_jw_1+w_j+w_{j+1})+\eta_j w^2_1+O_3&
	if $2\le j\le n$,\cr}
\eqno (4.5)
$$
when $k=1$ and by
$$
\tilde g_j(w)=\cases{w_1[1-\eps^1_kw_1+2w_2-w_{k+1}+O_2]& if $j=1$,\cr
	w_j[1-w_j+w_{j+1}+O_2]& if $2\le j\le k-1$,\cr
	w_k[1+\eps^1_kw_1-w_k+w_{k+1}+O_2]& if $j=k$,\cr
	\eps_j^1w_1+w_j+w_{j+1}\cr
	\phantom{\eps_j^1w_1}-(\eps^1_kw_1+w_{k+1})(\eps^1_jw_1+w_j+w_{j+1})+O_3& 
	if $k+1\le j\le n$,\cr}
\eqno (4.6)
$$
when $2\le k\le n$, where $\eta_j$ is the third-order term of the homogeneous
expansion of~$f_j$ evaluated in~$e_1$, and we have set $w_{n+1}\equiv 0$.
\hfill\break\indent
In particular, $d(\tilde F_n)_{{\bf e}_n}=\id$ always, and $d(\tilde F_{n-1})_{{\bf
e}_{n-1}}=\id$ if $\eps^1_n=0$.

\pf
Theorem~3.1 yields the existence of~$\tilde F_1$; since $\tilde F_1|_{E^1}$ is induced by
the differential of~$F$ at the origin, we see that ${\bf e}_1$ is a fixed point of~$\tilde
F_1$, and more generally that $\tilde F_1(Y^k)=Y^k$ for $k=2,\ldots,n$. 

We claim that $d(\tilde F_1)_{[v]}$ is invertible for all $[v]\in Y^{n-1}$. Indeed, fix such
a $[v]$, and choose $2\le k\le n-1$ such that $[v]\in Y^k\setminus Y^{k-1}$ (the case
$[v]={\bf e}_1$ has already been dealt with in Proposition~3.3). In particular, both $[v]$
and $\tilde F([v])$ belong to the domain of the $k$-th canonical chart of~$M^1$, which in
this proof we shall denote by~$(V,\chi)$. In particular, $[v]=\chi^{-1}(w^o)$, where
$w^o=(v_1/v_k,\ldots,v_{k-1}/v_k,0,\ldots,0)$. We now proceed as in the proof of
Proposition~3.2, using Eq.~(3.5) with $h=k$ and~$r=0$, and Eq.~(3.7) with $k$ instead
of~$r+1$ and $w^o$ instead of~$O$. The Jacobian matrix of $\chi\circ\tilde
F_1\circ\chi^{-1}$ at~$w^o$ turns out to be
$$
\left|\vcenter{\offinterlineskip
\halign{&\hfil$#$\hfil\cr
\multispan3\strut&\omit\hskip1.5mm\vrule width 1pt\hskip 1.5mm&\cr
\multispan3\strut\hskip1.5mm\hfil $J_{k-1}$\hfil&\omit\hskip1.5mm\vrule width 1pt
	\hskip 1.5mm&\multispan5\hfil$*$\hskip1.5mm\hfil\cr
\multispan3\strut&\omit\hskip1.5mm\vrule width 1pt\hskip 1.5mm&\cr
\noalign{\hrule height 1pt}
\multispan3&\omit\hskip1.5mm\vrule height 2pt width 1pt\hskip 1.5mm&\multispan5\cr
\strut&&&\omit\hskip1.5mm\vrule width 1pt\hskip 1.5mm&1&\omit\hskip1.5mm
	\vrule\hskip1.5mm&\quad&O&\quad\cr
\multispan3&\omit\hskip1.5mm\vrule width1pt&\multispan5\hskip-1mm\hrulefill\cr
\multispan3&\omit\hskip1.5mm\vrule width1pt&\strut&\omit\hskip1.5mm
	\vrule\hskip1.5mm&\cr
\multispan3\hskip1.5mm\quad\hfil$\smash{\vbox{\smash{$O$}}}$\quad
	\hfil&\omit\hskip1.5mm\vrule width 1pt&*&\omit\hskip1.5mm
	\vrule\hskip1.5mm&\quad&J_{n-k}&\quad\cr
\multispan3&\omit\hskip1.5mm\vrule width1pt&&\omit\hskip1.5mm
	\vrule\hskip1.5mm&\strut&\cr}}
\right|\;;
$$
in particular, the determinant is~1, and $d(\tilde F_1)_{[v]}$ is invertible.

By Proposition~3.3, $\tilde F_1$ is non-degenerate along $X^1=\{{\bf e}_1\}$, and so
Theorem~3.1 yields~$\tilde F_2$. Since $d\tilde F_1$ is invertible along~$Y^2$, we can
invoke Proposition~3.2 to prove that $d\tilde F_2$ is non-degenerate along~$X^2$, and thus
we get~$\tilde F_3$. Furthermore, being $d\tilde F_2$ invertible along~$X^2$, it is
invertible along the proper transform of~$Y^3$ too, because outside of~$E^2\subset X^2$ it
is given by~$d\tilde F_1$. Then we can again invoke Proposition~3.2 to prove that~$\tilde
F_3$ is non-degenerate along~$X^3$, and Theorem~3.1 yields~$\tilde F_4$. Repeating this
procedure we clearly get~$\tilde F_k$ for all~$k$.

Now we prove that ${\bf e}_k$ is a fixed point of~$\tilde F_k$. For $k=1$ it has already
been noticed. The action of~$\tilde F_2$ on~${\bf e}_2$ is given by~$d(\tilde F_1)_{{\bf
e}_1}$; then Eq.~(3.9) with $r=0$, $s=1$ and $\lambda=1$ clearly shows (because~$C=J_{n-1}$)
that~$\tilde F_2({\bf e}_2)={\bf e}_2$. Then we can repeat the argument, and Eq.~(3.9) with
$r=k-2$, $s=1$, $\lambda=1$ and~$C=J_{n-k+1}$ yields~$\tilde F_k({\bf e}_k)={\bf e}_k$. 

We are left to prove Eqs.~(4.5) and~(4.6). The main tool here is the equality
$$
F\circ(\chi_0\circ\pi_k\circ\chi_k^{-1})=(\chi_0\circ\pi_k\circ\chi_k^{-1})\circ\tilde G_k\;,
\eqno (4.7)
$$
which follows easily from $F\circ\pi_k=\pi_k\circ\tilde F_k$. Writing
$F=(f_1,\ldots,f_n)$ and recalling Lemma~4.2, Eq.~(4.7) yields
$$
\tilde g_j(w)=\cases{f_1(z)^2/f_k(z)& if $j=1$,\cr
	f_j(z)/f_{j-1}(z)& if $2\le j\le k$,\cr
	f_j(z)/f_k(z)& if $k+1\le j\le n$,\cr}
\eqno (4.8)
$$
where $z$ is given by Eq.~(4.3). Now to get Eqs.~(4.5) and~(4.6) from Eq.~(4.8) is just a
matter of computations. For instance, if $2\le j\le k-1$ we have
$$
\eqalign{{f_j(z)\over f_{j-1}(z)}&={\displaystyle w_1\prod_{h=2}^j(w_h)^2\prod_{h=j+1}^k
	w_h\left[1+w_{j+1}+\eps^1_jw_1\prod_{h=j+1}^k w_h+O_3\right]\over\displaystyle
	w_1\prod_{h=2}^{j-1}(w_h)^2\prod_{h=j}^k w_h\,\Bigl[1+w_j+O_3\Bigr]}\cr
	&=w_j\left[1+w_{j+1}+\eps^1_jw_1\prod_{h=j+1}^k w_h+O_3\right]\Bigl[1-w_j+w_j^2+O_3
	\Bigr]\cr
	&=w_j\left[1-w_j+w_{j+1}+\eps^1_jw_1\prod_{h=j+1}^k w_h+w_j^2-w_jw_{j+1}+O_3\right]\;,\cr}
$$
and similar computations yield the other components (showing by the way that it is easy to
get further terms in the homogeneous expansion of~$\tilde F_k$).\qedn

The map $\tilde F_n$ (or $\tilde F_{n-1}$ if $\eps^1_n=0$) is the {\it desingularization}
of~$F$. Using the canonical coordinates centered in~${\bf e}_n$ (respectively, in ${\bf
e}_{n-1}$) we have
$$
\tilde g_j(w)=\cases{w_1[1-\eps^1_n w_1+2w_2+O_2]& for $j=1$,\cr
	w_j[1-w_j+w_{j+1}+O_2]& for $2\le j\le n-1$,\cr
	w_n[1+\eps^1_n w_1-w_n+O_2]& for $j=n$,\cr}
$$
when $\eps^1_n\neq 0$;
$$
\tilde g_j(w)=\cases{w_1[1-\eps^1_{n-1} w_1+2w_2-w_n+O_2]& for $j=1$,\cr
	w_j[1-w_j+w_{j+1}+O_2]& for $2\le j\le n-2$,\cr
	w_{n-1}[1+\eps^1_{n-1}w_1-w_{n-1}+w_n+O_2]& for $j=n-1$,\cr
	w_n-(\eps^1_{n-1}w_1+w_n)w_n+O_3& for $j=n$,\cr}
$$
when $n\ge3$ and $\eps^1_n=0$;
$$
\tilde g_j(w)=\cases{w_1[1+\eps^1_1 w_1+w_2+O_2]& for $j=1$,\cr
	w_2+\eta_2(w_1)^2-(\eps^1_1w_1+w_2)w_2+O_3& for $j=2$,\cr}
$$
when $n=2$ and $\eps^1_2=0$.

\smallsect 5. Stable holomorphic curves

The main consequence of the Desingularization Theorem~4.3 is that after blowing up we have
diagonalized the differential of the map. This allows us to bring into play Hakim's theory,
that we shall now briefly summarize (though with a slightly different terminology).

Set $\Delta=\{\zeta\in\C\mid |\zeta-1|<1\}$. A
{\it holomorphic curve at the origin} is a holomorphic injective
map~$\phe\colon\Delta\to\C^n\setminus\{O\}$ such that $\phe$ extends continuosly
to~$0\in\de\Delta$ with~$\phe(0)=O$. 

Now take $F\in\End(\C^n,O)$. We shall say that a holomorphic curve at the origin~$\phe$, or
its image~$D=\phe(\Delta)$, is {\it $F$-invariant} if
$F\bigl(\phe(\Delta)\bigr)\subseteq\phe(\Delta)$; that it is {\it stable} if it is
$F$-invariant and $(F|_D)^k\to O$ uniformly on compact subsets of~$D$.
Finally, we shall say that $\phe$ is {\it tangent} to~$v\in\P^{n-1}(\C)$ if $[\phe(\zeta)]\to
v$ as $\zeta\to 0$.

Now let $P_2\colon\C^n\to\C^n$ be a $\C^n$-valued quadratic form. A {\it characteristic
direction} of~$P_2$ is a $v\in\C^n\setminus\{O\}$ such that $P_2(v)=\lambda v$. If
$\lambda=0$ then $v$ is {\it degenerate;} otherwise it is a {\it
non-degenerate} characteristic direction. 

Then Hakim's results can be summarized as follows:

\newthm Theorem \cuno: (Hakim [8, 9]) 
Let $F\in\End(\C^n,O)$ be such that $dF_O=\id$. Let $P_2\colon\C^n\to\C^n$ be the
quadratic part of the homogeneous expansion of~$F$. If $z^o\in\C^n$, set $z^k=F^k(z^o)$,
and denote by $[z^k]$ its image in $\P^{n-1}(\C)$ when $z^k\neq 0$. Then:
{\medskip
\item{\rm (i)}if $z^k\to O$ and $[z^k]\to[v]$ then $v$ is a characteristic direction of
$P_2$;
\item{\rm (ii)}if $v$ is a non-degenerate characteristic direction of $P_2$, then $F$ admits
a stable holomorphic curve at~$O$ tangent to~$[v]$;
\item{\rm (iii)}if $v$ is a non-degenerate characteristic direction of~$P_2$ with
$P_2(v)=\lambda v$ and $D\subset\C^n$ is the stable
holomorphic curve at the origin given by part~(ii), then for every $z^o\in D$ and $1\le
j\le n$ we have
$$
z^k_j=-{v_j\over\lambda k}+o\left({1\over k}\right)\;.
$$
}

Putting together Theorems~4.3 and~5.1 we are able to prove the existence of stable
holomorphic curve at the origin for generic maps with a Jordan fixed point at the origin.
We first consider the case $n\ge 3$:

\newthm Corollary \cdue:
Let $n\ge 3$, and $F\in\End(\C^n,O)$ be with a Jordan fixed point at the origin. Then:
{\medskip
\item{\rm (i)}If $\eps^1_n\neq 0$, then $F$ admits a stable holomorphic curve at the
origin~$\phe$ tangent to~$e_1$. Furthermore, if $F$ is in normal form, $z^o\in\phe(\Delta)$
and $z^k=F^k(z^o)$, then
$$
z^k_j=(-1)^{n+j-1}{2n-1\over \eps^1_n}{2n-2\choose n-1}{(n+j-2)!\over k^{n+j-1}}+o\left(
	{1\over k^{n+j-1}}\right)
\eqno (5.1)
$$
for $j=1,\ldots,n$.
\item{\rm (ii)}If $\eps^1_n=0$ but $\eps^1_{n-1}\neq 0$, then $F$ admits a stable
holomorphic curve at the origin~$\phe$ tangent to~$e_1$. Furthermore, if $F$ is in normal
form, $z^o\in\phe(\Delta)$ and $z^k=F^k(z^o)$, then
$$
z^k_j=\cases{\displaystyle(-1)^{n+j}{2n-3\over\eps^1_{n-1}}{2n-4\choose n-2}{(n+j-3)!
	\over k^{n+j-2}}+o\left({1\over k^{n+j-2}}\right)&for $j=1,\ldots, n-1$,\cr
	 \displaystyle o\left({1\over k^{2n-2}}\right)& for $j=n$.\cr}
\eqno (5.2)
$$
}

{\it Remark 5.1:}
Lemma~3.4 implies that if $O$ is a Jordan fixed point then a stable holomorphic curve at
the origin, if it is tangent to some vector, it must be tangent to~$e_1$.
\medbreak

\pf
We can clearly assume that $F$ is in normal form. The idea is to apply Theorem~5.1 to the
desingularization~$\tilde F_k$ of~$F$ (where
$k=n$ or~$n-1$), and then use~$\pi_k$ to project the result down to~$F$. Not all the
characteristic directions of the quadratic part at~${\bf e}_k$ of~$\tilde F_k$ are
allowable, though. Since we are working in~$M^k$, characteristic directions tangent
to~$\pi^{-1}_k(X^0)$ should be excluded, because the $\tilde F_k$-stable holomorphic curve
provided by Theorem~5.1.(ii) could be contained in the singular divisor, and thus it
would be killed by~$\pi_k$.  Equation~(4.1) says that $\pi_k^{-1}(X^0)$ is given by
$\{w_1=0\}\cup\cdots\cup\{w_k=0\}$; therefore we must look for characteristic directions~$v$
with~$v_1,\ldots,v_k\neq 0$. In general, characteristic directions not tangent to the
singular divisor $\pi_k^{-1}(X^0)$ will be called {\it allowable.}

We start by considering case (i). Then, since we are assuming $v_1,\ldots,v_n\neq 0$, an
allowable characteristic direction for~$\tilde F_n$ at~${\bf e}_n$ must satisfy 
$$
\cases{-\eps^1_nv_1+2v_2=\lambda\;,\cr
	-v_j+v_{j+1}=\lambda,& for $2\le j\le n-1$,\cr
	\eps^1_nv_1-v_n=\lambda.\cr}
$$
The unique solution of this system is
$$
v_j=\cases{(2n-1)\lambda/\eps^1_n& for $j=1$,\cr
		(n+j-2)\lambda& for $2\le j\le n$.\cr}
$$
This is an allowable solution; therefore Theorem~5.1.(ii) yields a $\tilde F_k$-stable
holomorphic curve~$\tilde\phe$ at the origin tangent to~$v$. Since $v$ is not tangent
to~$\pi_n^{-1}(X^0)$, which is invariant under~$\tilde F_n$, the image of the curve is
contained in $M^n\setminus\pi_n^{-1}(X^0)$, which is exactly the subset of~$M^n$
where~$\pi_n$ is a biholomorphism with~$\C^n\setminus\{O\}$. Therefore the holomorphic curve
$\phe=\pi_n\circ\tilde\phe$ is an $F$-stable holomorphic curve at the origin in~$\C^n$, and
Eq.~(5.1) follows from Lemma~4.2 and Theorem~5.1.(iii).

Now let us consider case (ii). This time we are assuming $v_1,\ldots,v_{n-1}\neq 0$; then an
allowable characteristic direction for~$\tilde F_{n-1}$ at~${\bf e}_{n-1}$ must satisfy
$$
\cases{-\eps^1_{n-1}v_1+2v_2-v_n=\lambda\;,\cr
	-v_j+v_{j+1}=\lambda,& for $2\le j\le n-2$,\cr
	\eps^1_{n-1}v_1-v_{n-1}+v_n=\lambda\;,\cr
	-v_n(\eps^1_{n-1}v_1+v_n)=\lambda v_n\;.\cr}
$$
If $v_n\neq 0$, the first and the last equations yield $v_2=0$, which is excluded. So
$v_n=0$, and the system reduces to the one studied in case (i). Therefore the unique
solution is
$$
v_j=\cases{(2n-3)\lambda/\eps^1_{n-1}& for $j=1$,\cr
		(n+j-3)\lambda& for $2\le j\le n-1$,\cr
	0&for $j=n$,\cr}
$$
and we get the assertion exactly as before.\qedn

{\it Remark 5.2:}
A different desingularization of the map~$F$ can be obtained simply by blowing up the
points~${\bf e}_0,\ldots,{\bf e}_{n-1}$. But it turns out that such a procedure is too rough:
we lose informations, and the desingularized map we get has {\it no} allowable
characteristic directions.
\medbreak

When $n=2$ and $\eps^1_2\neq 0$, arguing exactly as in part~(i) of the previous corollary we
obtain

\newthm Corollary \ctre:
Let $F\in\End(\C^2,O)$ be with a Jordan fixed point at the origin. Assume that
$\eps^1_2\neq 0$. Then $F$ admits a stable holomorphic curve at the
origin~$\phe$ tangent to~$e_1$. Furthermore, if $F$ is in normal form, $z^o\in\phe(\Delta)$
and $z^k=F^k(z^o)$, then
$$
\cases{\displaystyle z^k_1={6\over \eps^1_2}{1\over k^2}+o\left({1\over k^2}\right)\;,\cr
	\displaystyle z^k_2=-{12\over \eps^1_2}{1\over k^3}+o\left({1\over k^3}\right)\;.\cr}
$$

We now shall deal with the case $n=2$, $\eps^1_2=0$. We point out that for
the first time a coefficient of the cubic part of~$F$ enters directly into play even when
the quadratic part of~$F$ is not zero. 

Let $F$ be in normal form with $n=2$ and $\eps^1_2=0$. If $\eps^1_1\neq 0$, we set
$\Xi=1+2\eta_2/(\eps^1_1)^2$; Proposition~2.5 shows how to compute $\eps^1_1$ and~$\Xi$ for
maps not necessarily in normal form.

\newthm Corollary \cqua: Let $F\in\End(\C^2,O)$ be with a Jordan fixed point at the
origin. Assume that $\eps^1_2=0$. Then:
{\medskip
\item{\rm (i)}if $\eps^1_1\neq 0$ and $\Xi\neq 0$,~$1$, or if~$\eps^1_1=0$ and~$\eta_2\neq
0$, then $F$ admits two distinct stable holomorphic curves at the origin.
\item{\rm (ii)}if $\eps^1_1\neq 0$ and $\Xi=0$,~$1$ then $F$ admits one stable 
holomorphic curve at the origin.
\medskip
\noindent In both cases, the stable curves are tangent to~$e_1$. Furthermore, if $F$ is in
normal form, $z^o$ belongs to the image of the curve and $z^k=F^k(z^o)$, then
$z^k_1=c_1/k+o(1/k)$ and $z^k_2=c_2/k^2+o(1/k^2)$ for suitable $c_1\neq 0$ and~$c_2\in\C$.}

{\it Remark 5.3:} 
If $\eps^1_1=\eta_2=0$ several things might happen; we can even have more than
two stable holomorphic curves at the origin. See~[1] and~[3] for examples.
\medbreak

\pf
Let $\tilde F_1$ be the desingularization of~$F$. An easy computation shows that a
vector
$v$ is an allowable non-degenerate characteristic direction iff $v_1\neq 0$ and
$$
\eta_2\left({v_1\over\lambda}\right)^2+2\eps^1_1{v_1\over\lambda}-2=0\quad\hbox{and}
	\quad {v_2\over\lambda}=1-\eps^1_1{v_1\over\lambda}\;.
$$
Therefore, taking $v_1=1$ for the sake of definiteness, we find that
\medskip
\item{--}if $\eps^1_1\neq 0$ and $\Xi\neq 0$ then we get two allowable characteristic
directions $v_{\pm}=(1,u_{\pm})$, where
$$
u_{\pm}={-1\pm\sqrt{\Xi}\over 2}\eps^1_1\;,
\eqno (5.3)
$$
and $+\sqrt{\Xi}$ denotes the square root with positive real part (or with positive
imaginary part if~$\Xi$ is a negative number). The direction~$v_+$ is always
non-degenerate, whereas~$v_-$ is non-degenerate iff~$\Xi\neq 1$.

\item{--}if $\eps^1_1\neq 0$ and $\Xi=0$ then we get one allowable characteristic
direction $v_1=(1,-\eps^1_1/2)$ which is non-degenerate.

\item{--}if $\eps^1_1=0$ and $\eta_2\ne0$ then we get two allowable characteristic
directions $v_\pm=(1,\pm\sqrt{\eta_2/2})$, both non-degenerate.

\item{--}if $\eps^1_1=\eta_2=0$ then we get one allowable characteristic
direction $v_1=(1,0)$, but it is degenerate.
\medskip
\noindent The assertion then follows exactly as in the previous corollaries.\qedn

{\it Remark 5.4:}
A $\C^n$-valued quadratic form $P_2$ on~$\C^n$ induces on the projective space a holomorphic
map $\hat P_2\colon\P^{n-1}(\C)\setminus Z\to\P^{n-1}(\C)$, where $Z$ is the image
in~$\P^{n-1}(\C)$ of the cone $P_2^{_1}(O)\setminus\{O\}\subset\C^n$. If $v\in\C^n$ is a
non-degenerate characteristic direction for~$P_2$, then its image $[v]\in\P^{n-1}(\C)$ is a
fixed point of~$\hat P_2$. In particular, we may then consider the linear map
$$
A_{[v]}=d(\hat P_2)_{[v]}-\id\colon T_{[v]}\bigl(\P^{n-1}(\C)\bigr)\to T_{[v]}\bigl(
	\P^{n-1}(\C)\bigr)\;.
$$
An easy computation shows that this is the same matrix introduced by Hakim~[8, 9]. She
proved that, under the hypotheses of Theorem~5.1, if $A_{[v]}$ has $d\ge 0$ eigenvalues with
positive real part, then the map actually admits a stable holomorphic $(d+1)$-manifold
at the origin. Applying this fact to the case $n=2$, $\eps^1_2=0$ and~$\eps^1_1\neq 0$ of
the previous corollary, it is not difficult to check that, writing $\Xi=\rho e^{i\theta}$,
$A_{[v_-]}$ has positive real part iff
$$
0<\rho<{1+\cos\theta\over 2}=\cos^2\left({\theta\over2}\right)
$$
(and $A_{[v_+]}$ has always negative real part). Then in this case we have found a basin of
attraction for the origin.
\medbreak

{\it Remark 5.5:}
It is not difficult to compute the matrix $A_{[v]}$ for the allowable characteristic
directions described in the proof of Corollaries~5.2 and 5.3; it is not so easy to compute
the sign of the real part of the eigenvalues, though. For $n\le 7$ we checked that the
matrix~$A_{[v]}$ has no eigenvalue with positive real part, and we suspect that this is true
for all~$n$.
\medbreak

{\it Remark 5.6:}
Hakim~[9] proved that when $\tilde F\in\End(\C^n,O)$ is a global automorphism of~$\C^n$ with
$d\tilde F_O=\id$, and $v$ is a non-degenerate characteristic direction, then the
set~$\Omega_v$ of orbits $z^k\to O$ with $[z^k]\to[v]$ is an $\tilde F$-stable biholomorphic
image of~$\C^{d+1}$, where $d\ge 0$ is the number of eigenvalues of~$A_{[v]}$ with positive
real part (assuming, for simplicity, that $A_{[v]}$ has no purely imaginary eigenvalues).
This is still true in our situation. Indeed, if our map $F$ is a global automorphism
of~$\C^n$, then its desingularization~$\tilde F$ is a global automorphism
of~$M^k\setminus\pi^{-1}_k(X^0)$, which is biholomorphic to~$\C^n\setminus\{O\}$.
Furthermore, if $v$ is an allowable characteristic direction, then~$\Omega_v$ cannot
intersect the singular divisor, because the latter is~$\tilde F$-invariant whereas $v$ is
not tangent to it. This means that we can apply Hakim's result to~$\tilde F$, and
projecting down via~$\pi_k$ we get an $F$-stable $(d+1)$-manifold biholomorphic
to~$\C^{d+1}$. In particular, then, Remark~5.4 yields yet another instance of the
Fatou-Bieberbach phenomenon. 
\medbreak

\smallsect 6. Regular orbits

In the previous section we have shown that allowable (i.e., not tangent to the
singular divisor) characteristic directions of the desingularization of a map~$F$ give
rise to $F$-stable curves at the origin. A priori, other characteristic directions might
also give rise to $F$-stable curves, or at least to $F$-orbits converging to the origin. The
aim of this section is to show that in the generic case $\eps^1_n\neq0$ this cannot happen
--- and thus, roughly speaking, the local stable dynamics of $F$ is completely described by
Corollary~5.2.

To state more precisely our result, we need some definitions. Let
$\{z^k\}\subset\C^n\setminus\{O\}$ be a sequence converging to the origin. We shall say that
$\{z_k\}$ is {\it $0$-regular} if $\{[z^k]\}$ converges to some $[v]\in\P^{n-1}(\C)$; this
is equivalent to saying that $\pi_1^{-1}(z^k)$ converges to some $[v]\in E^1$. We shall say
that $\{z^k\}$ is {\it $1$-regular} if either~$[v]\neq{\bf e}_1$ (and we shall specify this
case saying that it is $1$-regular {\it of first kind\/}) or $[v]={\bf e}_1$
and~$\{\chi_1\circ\pi_1^{-1}(z^k)\}$ is 0-regular (and then $\{z^k\}$ is $1$-regular {\it of
second kind\/}). Now we proceed by induction. Let~$\{z^k\}$ be $(r-1)$-regular. If it
$(r-1)$-regular of first kind, we shall also say that it is {\it $r$-regular (of first
kind).} If it is $(r-1)$-regular of second kind, then $\pi_r^{-1}(z^k)$ converges to some
$[v]\in E^r$. We shall say that~$\{z^k\}$ is {\it $r$-regular} if either $[v]\neq{\bf e}_r$
(and then again $r$-regular of first kind) or $[v]={\bf e}_r$ and
$\{\chi_r\circ\pi^{-1}_r(z^k)\}$ is 0-regular (and then $\{z^k\}$ is $r$-regular {\it of
second kind\/}). We stress that we impose no conditions if $[v]\neq{\bf e_r}$; so for most
sequences $r$-regularity is equivalent to $0$-regularity.

The condition of $r$-regularity is just a way to require that the different components of
the sequence go to zero at comparable rates. For instance, if there
are
$a_j\in\C^*$ and $\delta_j>0$ such that
$$
z_j^k={a_j\over k^{\delta_j}}+o\left({1\over k^{\delta_j}}\right)\;,
$$
then Lemma~4.2 shows that $\{z^k\}$ is $r$-regular for every~$r$; but of course it is easy
to provide examples of much more general $r$-regular sequences.

Now let $F\in\End(\C^n,O)$ be with a Jordan fixed point at the origin. Assume that $F$ is in
normal form with $\eps^1_n\neq 0$, and let $\tilde F$ be its desingularization. We shall say
that an $F$-orbit is {\it regular} if it converges to the origin and it is $n$-regular.
A quick look to Lemma~4.2 shows
that orbits obtained pushing down $0$-regular orbits of $\tilde F$ tangent to allowable
characteristic directions are regular; such orbits are called {\it standard,} and are the
ones described in Corollaries~5.2.(i) and 5.3 (the orbits described in Corollary~5.2.(ii)
are easily seen to be $(n-1)$-regular).

Using this terminology, our aim is to prove that every regular orbit is
standard. To do so, we need a lemma:

\newthm Lemma \suno:
Let $\{w^k\}\subset\C^*$ be a sequence converging to~$0$. Assume there is another sequence
$\{u^k\}\subset\C$ such that $u^k/w^k\to c\in\C$ and
$$
w^{k+1}=w^k(1+u^k)+o\bigl((w^k)^2\bigr)\;.
$$
Then $1/(kw^k)\to -c$. In particular, if $c\neq 0$ we have
$$
w^k=-{1\over ck}+o\left({1\over k}\right)\;.
$$

\pf
Set $\eps^k=w^{k+1}-w^k-u^kw^k$, so that $\eps^k/(w^k)^2\to 0$. We then have
$$
{1\over w^{h+1}}={1\over w^h}-{u^h\over w^h}+{(u^h)^2/w^h+(u^h-1)\eps^h/(w^h)^2\over
	1+u^h+\eps^h/w^h}\;.
$$
Summing this equality for $h=0,\ldots, k-1$ and dividing by~$k$ we find
$$
{1\over kw^k}={1\over kw^0}-{1\over k}\sum_{h=0}^{k-1}{u^h\over w^h}+
	{1\over k}\sum_{h=0}^{k-1}{(u^h)^2/w^h+(u^h-1)\eps^h/(w^h)^2\over
	1+u^h+\eps^h/w^h}\;,
$$
and the assertion follows from the convergence of the averages of a converging sequence.\qedn

Then:

\newthm Theorem \sdue:
Let $F\in\End(\C^n,O)$ be with a Jordan fixed point at the origin. Assume that $F$ is in
normal form with $\eps^1_n\neq 0$. Then every regular orbit of $F$ is standard.

\pf
We can clearly assume $\eps^1_n=1$. Let $\{z^k=F^k(z^0)\}$ be a regular orbit; we first of
all want to prove, by induction, that $\pi_r^{-1}(z^k)\to{\bf e}_r$ for $r=1,\ldots,n$.

First of all, $0$-regularity yields $[z^k]\to[v]\in\P^{n-1}(\C)$; but then Lemma~3.4 and
Remark~3.2 imply that~$[v]={\bf e}_1$, and thus $\pi_1^{-1}(z^k)\to{\bf e}_1$. Exactly the
same argument shows that $\pi_2^{-1}(z^k)\to{\bf e}_2$.

Now assume that $\pi_r^{-1}(z^k)\to{\bf e}_r$ for some $2\le r\le n-1$, and put
$w^k=\chi_r\circ\pi_r^{-1}(z^k)$. The $0$-regularity of~$\{w^k\}$ implies that
$[w^k]\to[v]\in\P^{n-1}(\C)$; but then Lemma~3.4 and Remark~3.2 (applied to the blow-up
of~${\bf e}_r$ in~$M^r$) show that $v$ must be (canonically identified with) an eigenvector
of~$d(\tilde F_r)_{{\bf e}_r}$. Now, by Theorem~4.3 we know that, in the canonical
coordinates, $d(\tilde F_r)_{{\bf e}_r}$ is represented by the matrix
$$
\left|\vcenter{\offinterlineskip
\halign{&\hfil$#$\hfil\cr
\multispan5\strut&\omit\hskip1.5mm\vrule width 1pt\hskip 1.5mm&\cr
\multispan5\strut\hskip1.5mm\hfil $I_r$\hfil&\omit\hskip1.5mm\vrule width 1pt
	\hskip 1.5mm&\multispan3\hfil$O$\hskip1.5mm\hfil\cr
\multispan5\strut&\omit\hskip1.5mm\vrule width1pt\hskip 1.5mm&\cr
\noalign{\hrule height1pt}
\multispan5&\omit\hskip1.5mm\vrule height 2pt
	width 1pt\hskip 1.5mm&\multispan3\cr
\strut\eps^1_{r+1}&\omit\hskip 1mm\vrule\hskip1mm&\multispan3&\omit\hskip1.5mm\vrule 
	width 1pt\hskip 1.5mm&\multispan3\cr
\strut\vdots&\omit\hskip 1mm\vrule\hskip1mm&\quad&\hfil$\smash{\raise
	3pt\vbox{\smash{$O$}}}$\hfil&\quad&\omit\hskip1.5mm\vrule width1pt
	\hskip 1.5mm&\quad&\hfil$\smash{\raise 3pt\vbox{\smash{$J_{n-r}$}}}$\hfil&\quad\cr
\strut\eps^1_n&\omit\hskip 1mm\vrule\hskip1mm&\multispan3&\omit\hskip1.5mm\vrule width 1pt
	\hskip 1.5mm&\multispan3\cr}}
\right|\;.
$$
Therefore $v=(0,v_2,\ldots,v_{r+1},0,\ldots,0)$; to prove that $\pi^{-1}_{r+1}(z^k)\to{\bf
e}_{r+1}$ it suffices to show that $v_{r+1}\neq 0$.

Assume, by contradiction, $v_{r+1}=0$, and let $j_0=\max\{2\le j\le r\mid v_j\neq0\}$. We
know that $w^k_j/w^k_{j_0}\to v_j/v_{j_0}$ for all $j$; in particular, $w^k_j=O(w^k_{j_0})$
if $v_j\neq 0$, and $w^k_j=o(w^k_{j_0})$ if $v_j=0$. Then Eq.~(4.6) yields
$$
w^{k+1}_{j_0}=w^k_{j_0}(1-w^k_{j_0})+o\bigl((w^k_{j_0})^2\bigr)\;;
$$
hence using Lemma~6.1 we find
$$
w^k_{j_0}={1\over k}+o\left({1\over k}\right)\;,
$$
and so
$$
w^k_j={v_j/v_{j_0}\over k}+o\left({1\over k}\right)
\eqno (6.1)
$$
for all $j=1,\ldots,n$. 

We now claim that $v_j/v_{j_0}=j_0-j+1$ for all $j=2,\ldots,j_0$. We argue by induction on
$j_0-j$. Take~$j<j_0$ and assume that $v_{j+1}/v_{j_0}=j_0-j$. Noticing that $w_j^k\neq
0$ for all~$k$ and~$1\le j\le r$ (because $\pi^{-1}_r(z^k)$ does not belong to the
singular divisor), we can write
$$
{w_j^{k+1}\over w_j^k}=1-w_j^k+w_{j+1}^k+O\bigl((w_{j+1}^k)^2\bigr)\;.
$$
If $v_j=0$ we would get
$$
{w_j^{k+1}\over w_j^k}=1+{j_0-j\over k}+o\left({1\over k}\right)\;,
$$
which is impossible because the infinite product $\prod_k(w_j^{k+1}/w_j^k)$ is converging
to~zero. Therefore $v_j\neq 0$; but then applying Lemma~6.1 to
$$
w^{k+1}_j=w_j^k(1-w^k_j+w^k_{j+1})+o\bigl((w_j^k)^2\bigr)
$$
and recalling Eq.~(6.1) we get $v_j/v_{j_0}=j_0-j+1$, as claimed.

In particular we then have $v_2/v_{j_0}=j_0-1$, and so
$$
{w_1^{k+1}\over w_1^k}=1+{2(j_0-1)\over k}+o\left({1\over k}\right)\;,
$$
which is impossible. The contradiction arises because we assumed $v_{r+1}=0$; therefore we
must have $v_{r+1}\neq 0$, as claimed.

Summing up, we have in particular proved that $\pi_n^{-1}(z^k)\to{\bf e}_n$; set
$w^k=\chi_n\circ\pi^{-1}_n(z^k)$. Notice that, by construction, $w^k_j\neq 0$ for all~$k$
and~$j$. By $0$-regularity, $[w^k]\to[v]\in\P^{n-1}(\C)$; Theorem~5.1 then says that $v$
must be a characteristic direction of~$\tilde F_n$ at~${\bf e}_n$, that is a solution of
$$
\cases{-v_1^2+2v_1v_2=\lambda v_1\;,\cr
\noalign{\smallskip}
	-v_j^2+v_jv_{j+1}=\lambda v_j& for $j=2,\ldots,n-1$,\cr
	v_1v_n-v_n^2=\lambda v_n.\cr}
$$
To end the proof we must show that $v$ is allowable, that is that $v_j\ne 0$ for
$j=1,\ldots,n$.

Assume, by contradiction, that there is a $j_0$ such that $v_{j_0}\neq 0$ but $v_{j_0+1}=0$
(where by $v_{n+1}$ we mean~$v_1$, of course). Then it is easy to prove that
$v_j/v_{j_0}\in\N$ for all $j=1,\ldots,n$; in particular, $v_j/v_{j_0}$ is always
non-negative. Now we have
$$
w^{k+1}_{j_0}=w^k_{j_0}(1-w^k_{j_0})+o\bigl((w^k_{j_0})^2\bigr)\;;
$$
therefore Lemma~6.1 yields $w_{j_0}^k=1/k+o(1/k)$. Recalling Eq.~(6.1) we then get
$w_j^k=c_j/k+o(1/k)$ with $c_j\ge 0$ for all~$j=1,\ldots,n$. But then arguing exactly as in
the first part of the proof we show that $v_{j_0-1},\ldots,v_1\neq 0$; and then we
get~$v_n\neq 0$, and going up we finally arrive to prove $v_{j_0+1}\neq 0$,
contradiction.\qedn

\setref{10}
\beginsection References

\pre 1 M. Abate: Holomorphic dynamical systems with a Jordan fixed point, I! Preprint!
	1998

\book 2 L. Carleson, T.W. Gamelin: Complex dynamics! Springer-Verlag, Berlin,
1993

\pre 3 D. Coman, M. Dabija: On the dynamics of some diffeomorphisms of $\C^2$ near
parabolic fixed points! Preprint! 1996

\art 4 P. Fatou: Substitutions analytiques et \'equations fonctionnelles \`a deux
variables! Ann. Sc. Ec. Norm. Sup.! {} 1924 67-142

\pre 5 F. Forstneri\v c: Interpolation by holomorphic automorphisms and embeddings in
$\C^n$! Preprint! 1996

\book 6 P. Griffiths, J. Harris: Principles of algebraic geometry! Wyley, New York, 1978

\art 7 M. Hakim: Semi-attractive transformations of $\C^p$! Publ. Math! 38 1994
479-499

\art 8 M. Hakim: Analytic transformations of $(\C^p,0)$ tangent to the identity!
Duke Math. J.! 92 1998 403-428

\pre 9 M. Hakim: Stable pieces of manifolds in transformations tangent to the identity!
Preprint! 1997

\pre 10 J. Milnor: Dynamics in one complex variable: introductory
lectures! Preprint \#1990/5, SUNY StonyBrook, New York! 1990

\art 11 R. P\'erez-Marco: Fixed points and circle maps! Acta Math.! 179 1997
243-294

\art 12 T. Ueda: Local structure of analytic transformations of two complex variables,
I! J. Math. Kyoto Univ.! 26 1986 233-261

\art 13 T. Ueda: Local structure of analytic transformations of two complex variables, 
I\negthinspace I! J. Math. Kyoto Univ.! 31 1991 695-711

\art 14 B.J. Weickert: Attracting basins for automorphisms of $\C^2$! Invent. Math.!
132 1998 581-605

\art 15 H. Wu: Complex stable manifolds of holomorphic diffeomorphisms! Indiana Univ.
Math. J.! 42 1993 1349-1358

\art 16 J.-C. Yoccoz: Th\'eor\`eme de Siegel, nombres de Bryuno et polyn\^omes
quadratiques! Ast\'e\-ris\-que! 231 1995 3-88

\bye